\documentclass[A4j,11pt]{article}
\usepackage{amsfonts,amsmath,amscd}
\usepackage{graphics}
\usepackage{amssymb}
\begin{document}

\title{{\bf Volume-preserving mean curvature flow\\
for tubes in rank one symmetric spaces\\
of non-compact type}}
\author{{\bf Naoyuki Koike}}
\date{}
\maketitle

\begin{abstract}
First we investigate the evolutions of the radius function and its gradient along the volume-preserving mean 
curvature flow starting from a tube (of nonconstant radius) over a compact closed domain of a reflective 
submanifold in a symmetric space under certain condition for the radius function.  
Next, we prove that the tubeness is preserved along the flow in the case where the ambient space is a rank one 
symmetric space of non-compact type, the reflective submanifold is 
an invariant submanifold and the radius function of the initial tube is radial.  
Furthermore, in this case, we prove that the flow reaches to the invariant submanifold or 
it exists in infinite time and converges to another tube of constant mean curvature 
in the $C^{\infty}$-topology in infinite time.  
\end{abstract}

\vspace{0.5truecm}






\section{Introduction} 
Let $f_t$'s ($t\in[0,T)$) be a one-parameter $C^{\infty}$-family of immersions 
of an $n$-dimensional compact manifold $M$ into an $(n+1)$-dimensional Riemannian manifold $\overline M$, 
where $T$ is a positive constant or $T=\infty$.  Define a map 
$\widetilde f:M\times[0,T)\to\overline M$ by $\widetilde f(x,t)=f_t(x)$ 
($(x,t)\in M\times[0,T)$).  
Denote by $\pi_M$ the natural projection of $M\times[0,T)$ onto $M$.  
For a vector bundle $E$ over $M$, denote by $\pi_M^{\ast}E$ the induced bundle of $E$ 
by $\pi_M$.  
Also, denote by $H_t, g_t$ and $N_t$ the mean curvature, the induced metric and the outward unit normal vector of $f_t$, respectively.  Define the function $H$ over $M\times[0,T)$ by 
$H_{(x,t)}:=(H_t)_x$ ($(x,t)\in M\times [0,T)$), 
the section $g$ of $\pi_M^{\ast}(T^{(0,2)}M)$ by 
$g_{(x,t)}:=(g_t)_x$ ($(x,t)\in M\times [0,T)$) and the section $N$ of 
$\widetilde f^{\ast}(T\overline M)$ by $N_{(x,t)}:=(N_t)_x$ ($(x,t)\in M\times [0,T)$), where 
$T^{(0,2)}M$ is the tensor bundle of degree $(0,2)$ of $M$ and $T\overline M$ is the tangent 
bundle of $\overline M$.  
The average mean curvature $\overline H(:[0,T)\to{\Bbb R})$ is defined by 
$$\overline H_t:=\frac{\int_MH_tdv_{g_t}}{\int_Mdv_{g_t}},\leqno{(1.1)}$$
where $dv_{g_t}$ is the volume element of $g_t$.  
The flow $f_t$'s ($t\in[0,T)$) is called a {\it volume-preserving mean curvature flow} if 
it satisfies 
$$\widetilde f_{\ast}\left(\frac{\partial}{\partial t}\right)=(\overline H-H)N.\leqno{(1.2)}$$
In particular, if $f_t$'s are embeddings, then we call $M_t:=f_t(M)$'s 
$(0\in[0,T))$ rather than $f_t$'s $(0\in[0,T))$ a volume-preserving mean curvature flow.  
Note that, if $M$ has no boundary and if $f$ is an embedding, then, along this flow, 
the volume of $(M,g_t)$ decreases but the volume of the domain $D_t$ sorrounded by $f_t(M)$ is preserved 
invariantly.  

First we shall recall the result by M. Athanassenas ([A1,2]).  
Let $P_i$ ($i=1,2$) be affine hyperplanes in the $(n+1)$-dimensional Euclidean 
space ${\Bbb R}^{n+1}$ meeting a affine line ${\it l}$ orthogonally 
and $E$ a closed domain of ${\Bbb R}^{n+1}$ with $\partial E=P_1\cup P_2$.  
Also, let $M$ be a hypersurface of revolution in ${\Bbb R}^{n+1}$ such that 
$M\subset E$, $\partial\,M\subset P_1\cup P_2$ and that 
$M$ meets $P_1$ and $P_2$ orthogonally.  
Let $D$ be the closed domain surrouded by $P_1,P_2$ and $M$, and $d$ 
the distance between $P_1$ and $P_2$.  She ([A1,2]) proved the following fact.  

\vspace{0.3truecm}

\noindent
{\bf Fact 1.} {\sl Let $M_t$ ($t\in[0,T)$) be the volume-preserving mean curvature flow 
starting from $M$ such that $M_t$ meets $P_1$ and $P_2$ orthogonally for all $t\in[0,T)$.  
Then the following statements {\rm(i)} and {\rm (ii)} hold:

{\rm (i)} $M_t$ ($t\in[0,T)$) remain to be hypersurfaces of revolution.  

{\rm (ii)} If ${\rm Vol}(M)\leq\frac{{\rm Vol}(D)}{d}$ holds, then $T=\infty$ and 
as $t\to\infty$, the flow $M_t$ converges to the cylinder $C$ such that 
the volume of the closed domain surrounded by $P_1,P_2$ and $C$ is equal to 
${\rm Vol}(D)$.}

\vspace{0.3truecm}

E. Cabezas-Rivas and V. Miquel ([CM1,2,3]) proved the similar result in certain kinds of 
rotationally symmetric spaces.  
Let $\overline M$ be an $(n+1)$-dimensional rotationally symmetric space 
(i.e., $SO(n)$ acts on $\overline M$ isometrically and its fixed point set is a one-dimensional 
submanifold).  
Note that real space forms are rotationally symmetric spaces.  
Denote by ${\it l}$ the fixed point set of the action, which is an 
one-dimensional totally geodesic submanifold in $\overline M$.  
Let $P_i$ ($i=1,2$) be totally geodesic hypersurfaces (or equidistant hypersurfaces) in 
$\overline M$ meeting ${\it l}$ orthogonally and $E$ a closed domain of $\overline M$ with 
$\partial E=P_1\cup P_2$, where we note that they treat the case where $P_i$ ($i=1,2$) are totally geodesic 
hypersurfaces (resp. equidistant hypersurfaces) in [CM1,2] (resp. [CM3]).  
An embedded hypersurface $M$ in $\overline M$ is called 
a {\it hypersurface of revolution} if it is invariant with respect to the $SO(n)$-action.  
Let $M$ be a hypersurface of revolution in $\overline M$ such that 
$M\subset E$, $\partial\,M\subset P_1\cup P_2$ and that $M$ meets $P_1$ and $P_2$ 
orthogonally.  
Let $D$ be the closed domain surrouded by $P_1,P_2$ and $M$, and $d$ 
the distance between $P_1$ and $P_2$.  
They ([CM1,2,3]) proved the following fact.  

\vspace{0.3truecm}

\noindent
{\bf Fact 2.} {\sl Assume that ${\rm Sec}(v,w)<0$ for any $v\in T{\it l}$ and $w\in T^{\perp}{\it l}$ and that 
${\rm Sec}(w_1,w_2)\leq 0$ for any $w_1,\,w_2\in T^{\perp}{\it l}$, where 
${\rm Sec}(\cdot,\bullet)$ denotes the sectional curvature of the $2$-plane spanned by $\cdot$ and $\bullet$.  
Let $M_t$ ($t\in[0,T)$) be the volume-preserving mean curvature flow starting from $M$ such that $M_t$ meets 
$P_1$ and $P_2$ orthogonally for all $t\in[0,T)$.  Then the following statements {\rm(i)} and {\rm (ii)} hold:

{\rm (i)} $M_t$ ($t\in[0,T)$) remain to be hypersurfaces of revolution.  

{\rm (ii)}  If ${\rm Vol}(M)\leq C$ holds, where $C$ is a constant depending on 
${\rm Vol}(D)$ and $d$, then $T=\infty$ and, as $t\to\infty$, the flow $M_t$ ($t\in[0,T)$) converges to 
a hypersurface of revolution $C$ of constant mean curvature such that the volume of the closed domain surrounded 
by $P_1,P_2$ and $C$ is equal to ${\rm Vol}(D)$.}

\vspace{0.3truecm}

A {\it symmetric space of compact type} (resp. {\it non-compact type}) is a naturally 
reductive Riemannian homogeneous space $\overline M$ such that, for each point $p$ of $\overline M$, 
there exists an isometry of $\overline M$ having $p$ as an isolated fixed point and that the 
isometry group of $\overline M$ is a semi-simple Lie group each of whose irreducible factors is 
compact (resp. not compact) (see [He]).  
Note that symmetric spaces of compact type other than a sphere 
and symmetric spaces of non-compact type other than a (real) hyperbolic space are not 
rotationally symmetric.  
An {\it equifocal submanifold} in a (general) symmetric space is a compact submanifold (without boundary) 
satisfying the following conditions:

\vspace{0.2truecm}

\noindent
(E-i) the normal holonomy group of $M$ is trivial,

\vspace{0.1truecm}

\noindent
(E-ii) $M$ has a flat section, that is, for each $x\in M$, 
$\Sigma_x:=\exp^{\perp}(T^{\perp}_xM)$ is totally geodesic and the induced 
metric on $\Sigma_x$ is flat, where $T^{\perp}_xM$ is the normal space of $M$ 
at $x$ and $\exp^{\perp}$ is the normal exponential map of $M$.  

\vspace{0.1truecm}

\noindent
(E-iii) for each parallel normal vector field $v$ of $M$, the focal radii of $M$ 
along the normal geodesic $\gamma_{v_x}$ (with $\gamma'_{v_x}(0)=v_x$) are 
independent of the choice of $x\in M$, where $\gamma'_{v_x}(0)$ is the 
velocity vector of $\gamma_{v_x}$ at $0$.  

\vspace{0.2truecm}


\noindent
In [Ko2], we showed that the mean curvature flow starting from an equifocal submanifold in 
a symmetric space of compact type collapses to one of its focal submanifolds in finite 
time.  
In [Ko3], we showed that the mean curvature flow starting from a certain kind of 
(not necessarily compact) submanifold satisfying the above conditions (E-i), (E-ii) and 
(E-iii) in a symmetric space of non-compact type collapses to one of its focal 
submanifolds in finite time.  
The following question arise naturally:

\vspace{0.3truecm}

\noindent
{\it Question.} {\sl In what case, does the volume-preserving mean curvature flow starting 
from a submanifold in a symmetric space of compact type (or non-compact type) 
converges to a submanifold satisfying the above conditions (E-i), (E-ii) and (E-iii)?}

\vspace{0.3truecm}

\noindent
Let $M$ be an equifocal hypersurface in a rank ${\it l}(\geq 2)$ symmetric space $\overline M$ of compact type 
or non-compact type.  Then it admits a reflective focal submanifold $F$ 
and it is a tube (of constant radius) over $F$, 
where the ``{\it reflectivity}" means that the submanifold is a connected 
component of the fixed point set of an involutive isometry of $\overline M$ and a 
``{\it tube of constant radius} $r(>0)$ {\it over} $F$" means the image of 
$t_r(F):=\{\xi\in T^{\perp}F\,\vert\,\vert\vert\xi\vert\vert=r\}$ by the normal exponential 
map $\exp^{\perp}$ of $F$ under the assumption that the restriction 
$\exp^{\perp}\vert_{t_{r}(F)}$ of 
$\exp^{\perp}$ to $t_{r}(F)$ is an embedding, where $T^{\perp}F$ is the normal bundle of 
$F$ and $\vert\vert\cdot\vert\vert$ is the norm of $(\cdot)$.  
See [L1,2] about the classification of reflective submanifolds in symmetric spaces.  
Any reflective submanifold in a symmetric space $\overline M$ of compact type or non-compact 
type is a singular orbit of a Hermann action (i.e., the action of the symmetric subgroup 
of the isometry group of $\overline M$) (see [KT]).  
Note that even if T. Kimura and M. Tanaka ([KT]) proved this fact in compact type case, the proof is valid 
in non-compact type case.  From this fact, it is shown that $M$ is 
curvature-adapted, where ``the curvature-adpatedness" means that, for any point $x\in M$ and 
any normal vector $\xi$ of $M$ at $x$, $R(\cdot,\xi)\xi$ preserves the tangent space $T_xM$ of $M$ at 
$x$ invariantly, and that the restriction $R(\cdot,\xi)\xi\vert_{T_xM}$ of $R(\cdot,\xi)\xi$ to $T_xM$ 
and the shape operator $A_{\xi}$ commute to each other 
($R\,:\,$ the curvature tensor of $\overline M$).  
The notion of the curvature-adaptedness was introduced in [BV].  
For a non-constant positive-valued function $r$ over $F$, the image of 
$t_r(F):=\{\xi\in T^{\perp}F\,\vert\,\vert\vert\xi\vert\vert=r(\pi(\xi))\}$ by 
$\exp^{\perp}$ is called the {\it tube of non-constant radius} $r$ {\it over} $F$ 
in the case where the restriction 
$\exp^{\perp}\vert_{t_r(F)}$ of $\exp^{\perp}$ to $t_r(F)$ is an embedding, 
where $\pi$ is the bundle projection of $T^{\perp}F$.  
Note that $\exp^{\perp}\vert_{t_r(F)}$ is an embedding for 
a non-constant positive-valued function $r$ over $F$ such that $\max\,r$ is sufficiently small 
because $F$ is homogeneous.  
Since $F$ is reflective, so is also the normal umbrella 
$F^{\perp}_x:=\exp^{\perp}(T^{\perp}_xF)$ of $F$ at $x$ and hence $F^{\perp}_x$ is a 
symmetric space.  

\vspace{0.3truecm}

\noindent
{\bf Motivation.} 
If $F^{\perp}_x$ is a rank one symmetric space, then tubes over $F$ 
of constant radius satisfies the above conditions (E-i), (E-ii) and (E-iii).  
Hence, when $F^{\perp}_x$ is of rank one, 
it is very interesting to invesitigate in what case the volume-preserving mean curvature 
flow starting from a tube of non-constant radius over $F$ converges to a tube of constant radius 
over $F$.  

\vspace{0.3truecm}

Under this motivation, we try to derive a result similar to those of 
M. Athanassenas ([A1,2]) and E. Cabezas-Rivas and V. Miquel ([CM1,2]) in this paper.  

Let $\gamma:[0,\infty)\to \overline M$ be any normal geodesic of $F$.  
Denote by $r_{co}(\gamma)$ the first conjugate radius along the geodesic $\gamma$ in $F^{\perp}_x$, 
$r_{fo}(\gamma)$ the first focal radius of $F$ along $\gamma$.  
$r_{\gamma}:=\min\{r_{co}(\gamma),r_{fo}(\gamma)\}$.  
It is shown that, if $F^{\perp}_x$ also is of rank one, then $r_{\gamma}$ is independent of the choice of 
$\gamma$.  Hence we denote $r_{\gamma}$ by $r_F$ in this case.  
The setting in this paper is as follows.  

\vspace{0.25truecm}

\noindent
{\bf Setting (S).} Let $F$ be a reflective submanifold in a symmetric space 
$\overline M$ of compact type or non-compact type and $B$ be a compact closed domain in $F$ with smooth boundary 
which is star-shaped with respect to some $x_{\ast}\in B$ and does not intersect with the cut locus of $x_{\ast}$ in $F$.  
Assume that the normal umbrellas of $F$ are rank one symmetric spaces.  
Set $\displaystyle{P:=\mathop{\cup}_{x\in\partial B}F^{\perp}_x}$ and denote by 
$E$ the closed domain in $\overline M$ surrounded by $P$.  
Let $M:=t_{r_0}(B)$ and $f:=\exp^{\perp}\vert_{t_{r_0}(B)}$, where $r_0$ is a non-constant positive 
$C^{\infty}$-function over $B$ with $r_0<r_F$ such that ${\rm grad}\,r_0=0$ holds along $\partial B$.  
Denote by $D$ the closed domain surrouded by $P$ and $f(M)$.  
See Figure 1 about this setting.  In the case where $\overline M$ is of compact type, see also Figure 4.  

\vspace{0.3truecm}

\centerline{
\unitlength 0.1in
%
\hspace{2.65truecm}}

\vspace{0.5truecm}

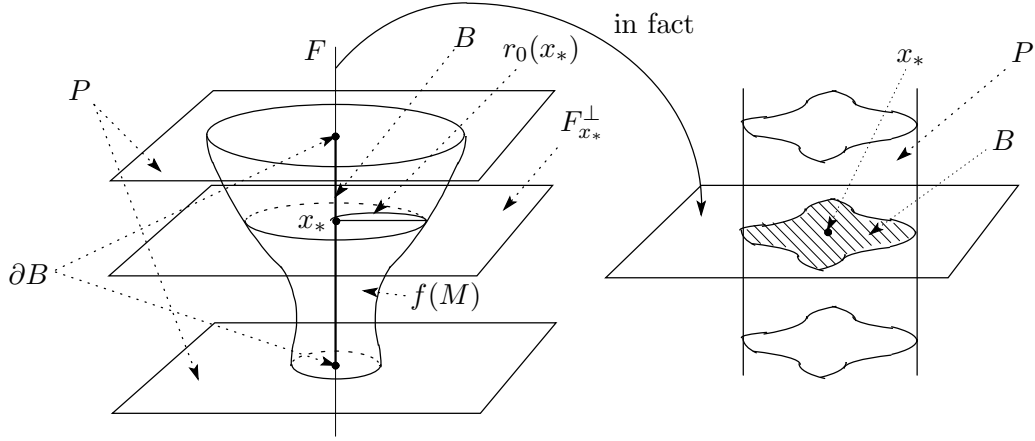
\centerline{{\bf Figure 1.}}

\vspace{0.3truecm}

\noindent
{\it Remark 1.1.} (i) At least one of singular orbits of any Hermann action of cohomogeneity one 
on any symmetric space $\overline M$ without Euclidean part other than a sphere and a (real) hyperbolic 
space is a reflective submanifold whose normal umbrellas are rank one symmetric spaces and 
tubes of constant radius over the reflective singular orbit satisfy the above conditions 
(E-i), (E-ii) and (E-iii) (i.e., of constant mean curvature).  
Note that, when $\overline M$ is of compact type, the Hermann action has exactly two singular 
orbits and, 
when $\overline M$ is of non-compact type, the Hermann action has the only one singular orbit.  
Hermann actions of cohomogeneity one on irreducible symmetric spaces of compact type or 
non-compact type are classified in [BT].  

(ii) At least one of singular orbits of any Hermann action of cohomogeneity greater than one 
on any symmetric space $\overline M$ of compact type or non-compact type is a reflective 
submanifold, but the normal umbrellas are higher rank symmetric spaces and hence tubes of constant radius over 
the reflective singular orbit do not satisfy the above conditions (E-i), (E-ii) and (E-iii) 
(i.e., not of constant mean curvature).  

\vspace{0.3truecm}

Here we shall state the difference of the above setting (S) (which includes the setting in [CM1,2]) 
from the setting in [CM3] (see $(1.3)$ in [CM3]).  
Note that the setting in [A1,2] is included in both settings because the equidistant hypersurfaces of 
a totally geodesic hypersurface in a Euclidean space are totally geodesic.  
The setting in [CM2] is as follows.  

\vspace{0.25truecm}

\noindent
{\bf Setting([CM2]).} Let $\overline M=J\times D^n(R)$ and 
$\overline g=f(r)^2dz^2+dr^2+h(r)^2g_{S^{n-1}}$, where 
$J$ are an interval, $D^n(R)$ is the ball of radius $R$ centered at the origin $0$ in ${\Bbb R}^n$, 
$S^{n-1}$ is the $(n-1)$-dimensional sphere of radius $1$ centered at $O$ in ${\Bbb R}^n$, $z$ is the natural 
coordinate of $J$, $r$ is the radius function from $0$ in ${\Bbb R}^n$, 
$g_{S^{n-1}}$ is the standard metric of $S^{n-1}$ and $f,h$ are positive functions.  
Let $J'$ be an closed subinterval of $J$.  
E. Cabezas-Rivas and V. Miquel ([CM2]) considered the volume-preserving mean curvature flow starting from 
a tube over the one-dimensional totally geodesic submanifold ${\it l}:=J'\times\{0\}$ in 
$(\overline M,\overline g)$ which is orthogonal to the totally geodesic hypersurfaces 
$\partial J'\times D^n(R)$.  
See Figure 2 about this setting in the case where $\overline M$ is of positive curvature.  

\vspace{0.3truecm}

\centerline{
\unitlength 0.1in
\begin{picture}( 39.6000, 33.1000)( 11.8000,-39.7000)
%
\special{pn 4}%
\special{pa 2490 1550}%
\special{pa 2490 3100}%
\special{fp}%
\put(16.3000,-14.2000){\makebox(0,0)[lb]{$S^{n-1}$}}%
\put(18.7000,-19.3000){\makebox(0,0)[lb]{$0$}}%
%
\special{pn 8}%
\special{pa 2090 2400}%
\special{pa 3000 2400}%
\special{fp}%
%
\special{pn 8}%
\special{ar 2090 2410 910 260  0.0000000 6.2831853}%
%
\special{pn 8}%
\special{ar 2090 2400 400 126  0.0000000 6.2831853}%
%
\special{pn 8}%
\special{pa 1920 1960}%
\special{pa 2086 2410}%
\special{dt 0.045}%
\special{sh 1}%
\special{pa 2086 2410}%
\special{pa 2082 2340}%
\special{pa 2068 2360}%
\special{pa 2044 2354}%
\special{pa 2086 2410}%
\special{fp}%
%
\special{pn 8}%
\special{ar 4370 1790 600 600  0.0000000 6.2831853}%
%
\special{pn 20}%
\special{ar 4370 1790 160 600  1.5707963 4.7123890}%
%
\special{pn 8}%
\special{ar 4370 1790 600 150  6.2831853 6.2831853}%
\special{ar 4370 1790 600 150  0.0000000 1.8202947}%
%
\special{pn 8}%
\special{ar 4240 1840 740 210  4.6940844 6.2831853}%
%
\special{pn 20}%
\special{ar 4250 1440 830 770  0.5112180 1.5571546}%
%
\special{pn 8}%
\special{ar 4360 3370 600 600  0.0000000 6.2831853}%
%
\special{pn 8}%
\special{ar 4360 3370 160 600  1.5707963 4.7123890}%
%
\special{pn 8}%
\special{ar 4230 3420 740 210  4.6940844 6.2831853}%
%
\special{pn 8}%
\special{ar 4240 3020 830 770  0.5112180 1.5571546}%
%
\special{pn 8}%
\special{pa 4010 1150}%
\special{pa 4260 1340}%
\special{dt 0.045}%
\special{sh 1}%
\special{pa 4260 1340}%
\special{pa 4220 1284}%
\special{pa 4218 1308}%
\special{pa 4196 1316}%
\special{pa 4260 1340}%
\special{fp}%
\put(39.8000,-11.2000){\makebox(0,0)[rb]{$J$}}%
\put(48.1000,-12.2000){\makebox(0,0)[lb]{$\{z_0\}\times D^n(R)$}}%
\put(50.4000,-20.2000){\makebox(0,0)[lt]{$(\overline M,\overline g)$}}%
%
\special{pn 8}%
\special{ar 4330 3370 160 600  4.7123890 4.7439679}%
\special{ar 4330 3370 160 600  4.8387048 4.8702837}%
\special{ar 4330 3370 160 600  4.9650206 4.9965995}%
\special{ar 4330 3370 160 600  5.0913363 5.1229153}%
\special{ar 4330 3370 160 600  5.2176521 5.2492311}%
\special{ar 4330 3370 160 600  5.3439679 5.3755469}%
\special{ar 4330 3370 160 600  5.4702837 5.5018627}%
\special{ar 4330 3370 160 600  5.5965995 5.6281785}%
\special{ar 4330 3370 160 600  5.7229153 5.7544942}%
\special{ar 4330 3370 160 600  5.8492311 5.8808100}%
\special{ar 4330 3370 160 600  5.9755469 6.0071258}%
\special{ar 4330 3370 160 600  6.1018627 6.1334416}%
\special{ar 4330 3370 160 600  6.2281785 6.2597574}%
\special{ar 4330 3370 160 600  6.3544942 6.3860732}%
\special{ar 4330 3370 160 600  6.4808100 6.5123890}%
\special{ar 4330 3370 160 600  6.6071258 6.6387048}%
\special{ar 4330 3370 160 600  6.7334416 6.7650206}%
\special{ar 4330 3370 160 600  6.8597574 6.8913363}%
\special{ar 4330 3370 160 600  6.9860732 7.0176521}%
\special{ar 4330 3370 160 600  7.1123890 7.1439679}%
\special{ar 4330 3370 160 600  7.2387048 7.2702837}%
\special{ar 4330 3370 160 600  7.3650206 7.3965995}%
\special{ar 4330 3370 160 600  7.4913363 7.5229153}%
\special{ar 4330 3370 160 600  7.6176521 7.6492311}%
\special{ar 4330 3370 160 600  7.7439679 7.7755469}%
\put(44.8000,-27.0000){\makebox(0,0)[rb]{$S^{n-1}$}}%
%
\special{pn 8}%
\special{ar 4340 3200 880 290  0.7878061 2.2857549}%
\put(49.3000,-37.3000){\makebox(0,0)[lt]{$D^n(R)$}}%
\put(50.0000,-27.4000){\makebox(0,0)[lb]{$r$-curves}}%
%
\special{pn 8}%
\special{pa 5138 2780}%
\special{pa 4678 3460}%
\special{dt 0.045}%
\special{sh 1}%
\special{pa 4678 3460}%
\special{pa 4732 3416}%
\special{pa 4708 3416}%
\special{pa 4700 3394}%
\special{pa 4678 3460}%
\special{fp}%
%
\special{pn 8}%
\special{pa 5138 2780}%
\special{pa 4788 3600}%
\special{dt 0.045}%
\special{sh 1}%
\special{pa 4788 3600}%
\special{pa 4834 3548}%
\special{pa 4810 3552}%
\special{pa 4796 3532}%
\special{pa 4788 3600}%
\special{fp}%
%
\special{pn 4}%
\special{pa 2090 1560}%
\special{pa 2090 2530}%
\special{fp}%
%
\special{pn 4}%
\special{pa 2090 2610}%
\special{pa 2090 3110}%
\special{fp}%
\put(22.3000,-14.4000){\makebox(0,0)[lb]{$J$}}%
%
\special{pn 20}%
\special{sh 1}%
\special{ar 2090 2410 10 10 0  6.28318530717959E+0000}%
\special{sh 1}%
\special{ar 2090 2410 10 10 0  6.28318530717959E+0000}%
%
\special{pn 8}%
\special{pa 3000 1560}%
\special{pa 3000 3120}%
\special{fp}%
%
\special{pn 8}%
\special{pa 1180 1570}%
\special{pa 1180 3130}%
\special{fp}%
%
\special{pn 4}%
\special{pa 1690 1560}%
\special{pa 1690 3120}%
\special{fp}%
%
\special{pn 8}%
\special{pa 2220 1440}%
\special{pa 2090 1600}%
\special{dt 0.045}%
\special{sh 1}%
\special{pa 2090 1600}%
\special{pa 2148 1562}%
\special{pa 2124 1560}%
\special{pa 2118 1536}%
\special{pa 2090 1600}%
\special{fp}%
%
\special{pn 8}%
\special{pa 1740 1450}%
\special{pa 1858 2288}%
\special{dt 0.045}%
\special{sh 1}%
\special{pa 1858 2288}%
\special{pa 1870 2220}%
\special{pa 1852 2236}%
\special{pa 1830 2226}%
\special{pa 1858 2288}%
\special{fp}%
%
\special{pn 8}%
\special{ar 4220 1830 440 200  3.1415927 4.7123890}%
%
\special{pn 8}%
\special{ar 4210 1820 440 120  1.5560590 3.1415927}%
%
\special{pn 20}%
\special{ar 4250 1830 480 380  1.5561346 3.1415927}%
%
\special{pn 8}%
\special{ar 4210 3410 440 200  3.1415927 4.7123890}%
%
\special{pn 8}%
\special{ar 4240 3410 480 380  1.5561346 3.1415927}%
%
\special{pn 8}%
\special{pa 4210 1940}%
\special{pa 4260 1940}%
\special{pa 4260 1890}%
\special{pa 4210 1890}%
\special{pa 4210 1940}%
\special{fp}%
%
\special{pn 8}%
\special{pa 4220 1620}%
\special{pa 4270 1620}%
\special{pa 4270 1580}%
\special{pa 4220 1580}%
\special{pa 4220 1620}%
\special{fp}%
%
\special{pn 8}%
\special{pa 4250 2210}%
\special{pa 4290 2210}%
\special{pa 4290 2160}%
\special{pa 4250 2160}%
\special{pa 4250 2210}%
\special{fp}%
%
\special{pn 8}%
\special{pa 4210 3210}%
\special{pa 4260 3210}%
\special{pa 4260 3160}%
\special{pa 4210 3160}%
\special{pa 4210 3210}%
\special{fp}%
%
\special{pn 8}%
\special{pa 4200 3490}%
\special{pa 4250 3490}%
\special{pa 4250 3430}%
\special{pa 4200 3430}%
\special{pa 4200 3490}%
\special{fp}%
%
\special{pn 8}%
\special{pa 4240 3790}%
\special{pa 4290 3790}%
\special{pa 4290 3740}%
\special{pa 4240 3740}%
\special{pa 4240 3790}%
\special{fp}%
%
\special{pn 8}%
\special{ar 4630 1790 100 536  1.5707963 4.7123890}%
%
\special{pn 8}%
\special{ar 4100 1790 170 530  1.6902253 4.5975124}%
%
\special{pn 8}%
\special{ar 4710 3370 140 560  1.8897140 4.3506597}%
%
\special{pn 8}%
\special{pa 4470 2720}%
\special{pa 4620 2920}%
\special{dt 0.045}%
\special{sh 1}%
\special{pa 4620 2920}%
\special{pa 4596 2856}%
\special{pa 4588 2878}%
\special{pa 4564 2880}%
\special{pa 4620 2920}%
\special{fp}%
%
\special{pn 8}%
\special{pa 5140 2780}%
\special{pa 4680 3250}%
\special{dt 0.045}%
\special{sh 1}%
\special{pa 4680 3250}%
\special{pa 4742 3216}%
\special{pa 4718 3212}%
\special{pa 4712 3188}%
\special{pa 4680 3250}%
\special{fp}%
%
\special{pn 20}%
\special{ar 4380 1970 630 580  3.3973145 6.0090179}%
%
\special{pn 8}%
\special{pa 4530 1930}%
\special{pa 4580 1930}%
\special{pa 4580 1890}%
\special{pa 4530 1890}%
\special{pa 4530 1930}%
\special{fp}%
%
\special{pn 8}%
\special{pa 4540 1650}%
\special{pa 4590 1650}%
\special{pa 4590 1610}%
\special{pa 4540 1610}%
\special{pa 4540 1650}%
\special{fp}%
%
\special{pn 8}%
\special{pa 3940 1920}%
\special{pa 3990 1920}%
\special{pa 3990 1880}%
\special{pa 3940 1880}%
\special{pa 3940 1920}%
\special{fp}%
%
\special{pn 20}%
\special{sh 1}%
\special{ar 4200 1930 10 10 0  6.28318530717959E+0000}%
\special{sh 1}%
\special{ar 4200 1930 10 10 0  6.28318530717959E+0000}%
\put(41.6000,-19.0000){\makebox(0,0)[rb]{$0$}}%
%
\special{pn 20}%
\special{sh 1}%
\special{ar 4210 1630 10 10 0  6.28318530717959E+0000}%
\special{sh 1}%
\special{ar 4210 1630 10 10 0  6.28318530717959E+0000}%
%
\special{pn 8}%
\special{pa 4630 880}%
\special{pa 4390 1380}%
\special{dt 0.045}%
\special{sh 1}%
\special{pa 4390 1380}%
\special{pa 4438 1330}%
\special{pa 4414 1332}%
\special{pa 4402 1312}%
\special{pa 4390 1380}%
\special{fp}%
%
\special{pn 8}%
\special{pa 4640 880}%
\special{pa 4730 2070}%
\special{dt 0.045}%
\special{sh 1}%
\special{pa 4730 2070}%
\special{pa 4746 2002}%
\special{pa 4726 2018}%
\special{pa 4706 2006}%
\special{pa 4730 2070}%
\special{fp}%
\put(44.0000,-8.3000){\makebox(0,0)[lb]{$\partial J'\times D^n(R)$}}%
\put(41.9000,-16.3000){\makebox(0,0)[rb]{$z_0$}}%
%
\special{pn 8}%
\special{pa 4970 1260}%
\special{pa 4770 1690}%
\special{dt 0.045}%
\special{sh 1}%
\special{pa 4770 1690}%
\special{pa 4816 1638}%
\special{pa 4792 1642}%
\special{pa 4780 1622}%
\special{pa 4770 1690}%
\special{fp}%
%
\special{pn 8}%
\special{pa 4560 1420}%
\special{pa 4550 1470}%
\special{pa 4610 1500}%
\special{pa 4620 1450}%
\special{pa 4620 1450}%
\special{pa 4560 1420}%
\special{fp}%
%
\special{pn 8}%
\special{pa 3950 1550}%
\special{pa 3950 1610}%
\special{pa 4000 1590}%
\special{pa 4000 1520}%
\special{pa 4000 1520}%
\special{pa 3950 1550}%
\special{fp}%
%
\special{pn 8}%
\special{pa 4250 1450}%
\special{pa 4300 1450}%
\special{pa 4300 1410}%
\special{pa 4250 1410}%
\special{pa 4250 1450}%
\special{fp}%
%
\special{pn 8}%
\special{pa 3970 2130}%
\special{pa 4000 2140}%
\special{pa 4000 2100}%
\special{pa 3960 2080}%
\special{pa 3960 2080}%
\special{pa 3970 2130}%
\special{fp}%
%
\special{pn 8}%
\special{pa 4560 2160}%
\special{pa 4610 2130}%
\special{pa 4610 2060}%
\special{pa 4550 2090}%
\special{pa 4550 2090}%
\special{pa 4560 2160}%
\special{fp}%
%
\special{pn 8}%
\special{pa 1370 1220}%
\special{pa 1500 2390}%
\special{dt 0.045}%
\special{sh 1}%
\special{pa 1500 2390}%
\special{pa 1514 2322}%
\special{pa 1494 2338}%
\special{pa 1474 2326}%
\special{pa 1500 2390}%
\special{fp}%
\put(12.6000,-11.7000){\makebox(0,0)[lb]{$D^n(R)$}}%
%
\special{pn 8}%
\special{pa 2930 1230}%
\special{pa 2770 2400}%
\special{dt 0.045}%
\special{sh 1}%
\special{pa 2770 2400}%
\special{pa 2800 2338}%
\special{pa 2778 2348}%
\special{pa 2760 2332}%
\special{pa 2770 2400}%
\special{fp}%
\special{pa 2770 2400}%
\special{pa 2770 2400}%
\special{dt 0.045}%
\put(27.2000,-11.8000){\makebox(0,0)[lb]{$r$-curve}}%
\end{picture}%
\hspace{1.55truecm}}

\vspace{0.3truecm}

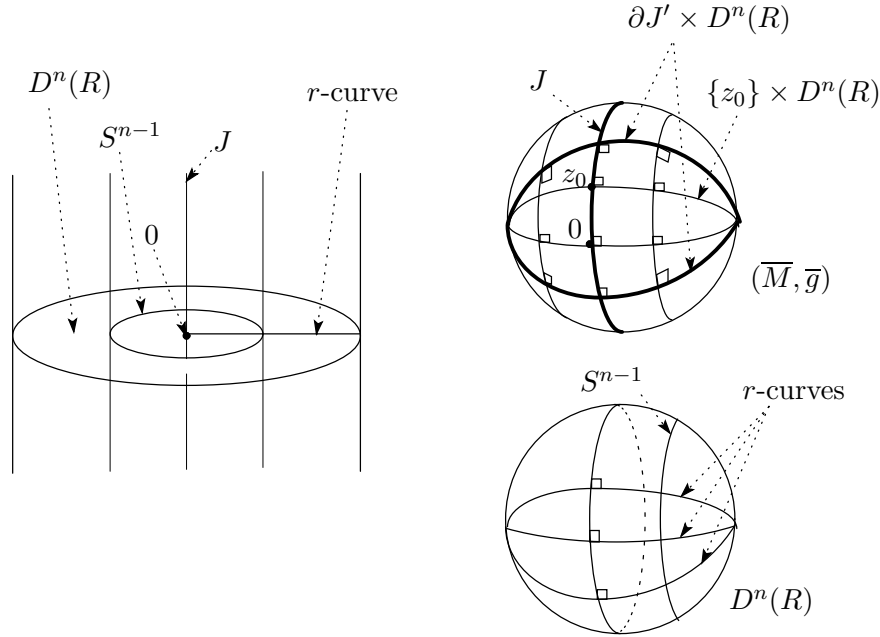
\centerline{{\bf Figure 2.}}

\vspace{0.25truecm}

The setting in [CM3] is as follows.  

\vspace{0.25truecm}

\noindent
{\bf Setting([CM3]).} 
Let $\overline M=J\times D^n(R)$ and 
$\overline g=dz^2+f(z)^2dr^2+f(z)^2h(r)^2g_{S^{n-1}}$, where 
$J,D^n(R),S^{n-1},z,r$ and $g_{S^{n-1}}$ are as above and $f,h$ are positive functions.  
Let $J'$ be an closed subinterval of $J$.  
E. Cabezas-Rivas and V. Miquel ([CM3]) considered the volume-preserving mean curvature flow starting from 
a tube over the one-dimensional totally geodesic submanifold ${\it l}:=J'\times\{0\}$ in 
$(\overline M,\overline g)$ which is orthogonal to the equidistant hypersurfaces 
$\partial J'\times D^n(R)$ (which are not totally geodesic).  
See Figure 3 about this setting in the case where $\overline M$ is of positive curvature.  

\vspace{0.5truecm}

\centerline{
\unitlength 0.1in
\begin{picture}( 43.9000, 32.5000)(  7.5000,-39.7000)
%
\special{pn 20}%
\special{sh 1}%
\special{ar 2680 2410 10 10 0  6.28318530717959E+0000}%
\special{sh 1}%
\special{ar 2680 2410 10 10 0  6.28318530717959E+0000}%
%
\special{pn 4}%
\special{pa 2550 1550}%
\special{pa 2550 3100}%
\special{fp}%
\put(17.7000,-12.1000){\makebox(0,0)[rb]{$S^{n-1}$}}%
\put(18.8000,-19.5000){\makebox(0,0)[lb]{$0$}}%
%
\special{pn 8}%
\special{pa 2090 2410}%
\special{pa 3000 2410}%
\special{fp}%
%
\special{pn 8}%
\special{ar 2090 2410 910 260  0.0000000 6.2831853}%
%
\special{pn 8}%
\special{ar 2090 2410 460 136  0.0000000 6.2831853}%
%
\special{pn 8}%
\special{pa 1930 1970}%
\special{pa 2082 2424}%
\special{dt 0.045}%
\special{sh 1}%
\special{pa 2082 2424}%
\special{pa 2080 2354}%
\special{pa 2064 2374}%
\special{pa 2042 2368}%
\special{pa 2082 2424}%
\special{fp}%
%
\special{pn 8}%
\special{ar 4370 1790 600 600  0.0000000 6.2831853}%
%
\special{pn 8}%
\special{ar 4370 1790 160 600  1.5707963 4.7123890}%
%
\special{pn 8}%
\special{ar 4370 1790 600 150  6.2831853 6.2831853}%
\special{ar 4370 1790 600 150  0.0000000 1.8202947}%
%
\special{pn 8}%
\special{ar 4360 3370 600 600  0.0000000 6.2831853}%
%
\special{pn 8}%
\special{ar 4360 3370 160 600  1.5707963 4.7123890}%
%
\special{pn 8}%
\special{ar 4230 3420 740 210  4.6940844 6.2831853}%
%
\special{pn 8}%
\special{ar 4240 3020 830 770  0.5112180 1.5571546}%
\put(36.5000,-15.5000){\makebox(0,0)[rb]{$J$}}%
\put(49.7000,-20.8000){\makebox(0,0)[lt]{$(\overline M,\overline g)$}}%
%
\special{pn 8}%
\special{ar 4330 3370 160 600  4.7123890 4.7439679}%
\special{ar 4330 3370 160 600  4.8387048 4.8702837}%
\special{ar 4330 3370 160 600  4.9650206 4.9965995}%
\special{ar 4330 3370 160 600  5.0913363 5.1229153}%
\special{ar 4330 3370 160 600  5.2176521 5.2492311}%
\special{ar 4330 3370 160 600  5.3439679 5.3755469}%
\special{ar 4330 3370 160 600  5.4702837 5.5018627}%
\special{ar 4330 3370 160 600  5.5965995 5.6281785}%
\special{ar 4330 3370 160 600  5.7229153 5.7544942}%
\special{ar 4330 3370 160 600  5.8492311 5.8808100}%
\special{ar 4330 3370 160 600  5.9755469 6.0071258}%
\special{ar 4330 3370 160 600  6.1018627 6.1334416}%
\special{ar 4330 3370 160 600  6.2281785 6.2597574}%
\special{ar 4330 3370 160 600  6.3544942 6.3860732}%
\special{ar 4330 3370 160 600  6.4808100 6.5123890}%
\special{ar 4330 3370 160 600  6.6071258 6.6387048}%
\special{ar 4330 3370 160 600  6.7334416 6.7650206}%
\special{ar 4330 3370 160 600  6.8597574 6.8913363}%
\special{ar 4330 3370 160 600  6.9860732 7.0176521}%
\special{ar 4330 3370 160 600  7.1123890 7.1439679}%
\special{ar 4330 3370 160 600  7.2387048 7.2702837}%
\special{ar 4330 3370 160 600  7.3650206 7.3965995}%
\special{ar 4330 3370 160 600  7.4913363 7.5229153}%
\special{ar 4330 3370 160 600  7.6176521 7.6492311}%
\special{ar 4330 3370 160 600  7.7439679 7.7755469}%
\put(44.8000,-27.0000){\makebox(0,0)[rb]{$S^{n-1}$}}%
%
\special{pn 8}%
\special{ar 4340 3200 880 290  0.7878061 2.2857549}%
\put(49.6000,-37.0000){\makebox(0,0)[lt]{$D^n(R)$}}%
\put(48.6000,-27.4000){\makebox(0,0)[lb]{$r$-curves}}%
%
\special{pn 8}%
\special{pa 5138 2780}%
\special{pa 4678 3460}%
\special{dt 0.045}%
\special{sh 1}%
\special{pa 4678 3460}%
\special{pa 4732 3416}%
\special{pa 4708 3416}%
\special{pa 4700 3394}%
\special{pa 4678 3460}%
\special{fp}%
%
\special{pn 8}%
\special{pa 5138 2780}%
\special{pa 4788 3600}%
\special{dt 0.045}%
\special{sh 1}%
\special{pa 4788 3600}%
\special{pa 4834 3548}%
\special{pa 4810 3552}%
\special{pa 4796 3532}%
\special{pa 4788 3600}%
\special{fp}%
%
\special{pn 4}%
\special{pa 2090 1560}%
\special{pa 2090 2530}%
\special{fp}%
%
\special{pn 4}%
\special{pa 2090 2610}%
\special{pa 2090 3110}%
\special{fp}%
\put(22.3000,-14.4000){\makebox(0,0)[lb]{$J$}}%
%
\special{pn 20}%
\special{sh 1}%
\special{ar 2090 2410 10 10 0  6.28318530717959E+0000}%
\special{sh 1}%
\special{ar 2090 2410 10 10 0  6.28318530717959E+0000}%
\put(43.4000,-8.9000){\makebox(0,0)[rb]{$\partial J'\times D^n(R)$}}%
%
\special{pn 8}%
\special{pa 3000 1560}%
\special{pa 3000 3120}%
\special{fp}%
%
\special{pn 8}%
\special{pa 1180 1570}%
\special{pa 1180 3130}%
\special{fp}%
%
\special{pn 4}%
\special{pa 1630 1550}%
\special{pa 1630 3110}%
\special{fp}%
%
\special{pn 8}%
\special{pa 2220 1440}%
\special{pa 2090 1600}%
\special{dt 0.045}%
\special{sh 1}%
\special{pa 2090 1600}%
\special{pa 2148 1562}%
\special{pa 2124 1560}%
\special{pa 2118 1536}%
\special{pa 2090 1600}%
\special{fp}%
%
\special{pn 8}%
\special{ar 4210 1820 440 120  1.5560590 3.1415927}%
%
\special{pn 8}%
\special{ar 4210 3410 440 200  3.1415927 4.7123890}%
%
\special{pn 8}%
\special{ar 4240 3410 480 380  1.5561346 3.1415927}%
%
\special{pn 8}%
\special{pa 4210 1940}%
\special{pa 4260 1940}%
\special{pa 4260 1890}%
\special{pa 4210 1890}%
\special{pa 4210 1940}%
\special{fp}%
%
\special{pn 8}%
\special{pa 4220 1610}%
\special{pa 4270 1610}%
\special{pa 4270 1570}%
\special{pa 4220 1570}%
\special{pa 4220 1610}%
\special{fp}%
%
\special{pn 8}%
\special{pa 4250 2210}%
\special{pa 4290 2210}%
\special{pa 4290 2160}%
\special{pa 4250 2160}%
\special{pa 4250 2210}%
\special{fp}%
%
\special{pn 8}%
\special{pa 4210 3210}%
\special{pa 4260 3210}%
\special{pa 4260 3160}%
\special{pa 4210 3160}%
\special{pa 4210 3210}%
\special{fp}%
%
\special{pn 8}%
\special{pa 4200 3490}%
\special{pa 4250 3490}%
\special{pa 4250 3430}%
\special{pa 4200 3430}%
\special{pa 4200 3490}%
\special{fp}%
%
\special{pn 8}%
\special{pa 4240 3790}%
\special{pa 4290 3790}%
\special{pa 4290 3740}%
\special{pa 4240 3740}%
\special{pa 4240 3790}%
\special{fp}%
%
\special{pn 8}%
\special{ar 4630 1790 100 536  1.5707963 4.7123890}%
\put(46.0000,-10.6000){\makebox(0,0)[lb]{$\{z_0\}\times D^n(R)$}}%
%
\special{pn 8}%
\special{ar 4100 1790 170 530  1.6902253 4.5975124}%
%
\special{pn 8}%
\special{ar 4710 3370 140 560  1.8897140 4.3506597}%
%
\special{pn 8}%
\special{pa 4470 2720}%
\special{pa 4620 2920}%
\special{dt 0.045}%
\special{sh 1}%
\special{pa 4620 2920}%
\special{pa 4596 2856}%
\special{pa 4588 2878}%
\special{pa 4564 2880}%
\special{pa 4620 2920}%
\special{fp}%
%
\special{pn 8}%
\special{pa 5140 2780}%
\special{pa 4680 3250}%
\special{dt 0.045}%
\special{sh 1}%
\special{pa 4680 3250}%
\special{pa 4742 3216}%
\special{pa 4718 3212}%
\special{pa 4712 3188}%
\special{pa 4680 3250}%
\special{fp}%
%
\special{pn 20}%
\special{ar 4370 1510 530 60  6.2831853 6.2831853}%
\special{ar 4370 1510 530 60  0.0000000 3.1415927}%
%
\special{pn 20}%
\special{ar 4360 2054 642 156  0.7915646 2.3598777}%
%
\special{pn 20}%
\special{ar 4370 1820 160 570  2.3680987 3.5995097}%
%
\special{pn 8}%
\special{ar 4370 1500 530 60  3.1415927 3.1822706}%
\special{ar 4370 1500 530 60  3.3043045 3.3449825}%
\special{ar 4370 1500 530 60  3.4670164 3.5076943}%
\special{ar 4370 1500 530 60  3.6297282 3.6704062}%
\special{ar 4370 1500 530 60  3.7924401 3.8331181}%
\special{ar 4370 1500 530 60  3.9551520 3.9958299}%
\special{ar 4370 1500 530 60  4.1178638 4.1585418}%
\special{ar 4370 1500 530 60  4.2805757 4.3212537}%
\special{ar 4370 1500 530 60  4.4432876 4.4839655}%
\special{ar 4370 1500 530 60  4.6059994 4.6466774}%
\special{ar 4370 1500 530 60  4.7687113 4.8093893}%
\special{ar 4370 1500 530 60  4.9314232 4.9721011}%
\special{ar 4370 1500 530 60  5.0941350 5.1348130}%
\special{ar 4370 1500 530 60  5.2568469 5.2975249}%
\special{ar 4370 1500 530 60  5.4195588 5.4602367}%
\special{ar 4370 1500 530 60  5.5822706 5.6229486}%
\special{ar 4370 1500 530 60  5.7449825 5.7856605}%
\special{ar 4370 1500 530 60  5.9076943 5.9483723}%
\special{ar 4370 1500 530 60  6.0704062 6.1110842}%
\special{ar 4370 1500 530 60  6.2331181 6.2737960}%
%
\special{pn 8}%
\special{ar 4370 2250 642 156  3.9233076 3.9533828}%
\special{ar 4370 2250 642 156  4.0436083 4.0736835}%
\special{ar 4370 2250 642 156  4.1639091 4.1939843}%
\special{ar 4370 2250 642 156  4.2842098 4.3142850}%
\special{ar 4370 2250 642 156  4.4045106 4.4345858}%
\special{ar 4370 2250 642 156  4.5248114 4.5548865}%
\special{ar 4370 2250 642 156  4.6451121 4.6751873}%
\special{ar 4370 2250 642 156  4.7654129 4.7954880}%
\special{ar 4370 2250 642 156  4.8857136 4.9157888}%
\special{ar 4370 2250 642 156  5.0060144 5.0360895}%
\special{ar 4370 2250 642 156  5.1263151 5.1563903}%
\special{ar 4370 2250 642 156  5.2466159 5.2766911}%
\special{ar 4370 2250 642 156  5.3669166 5.3969918}%
\special{ar 4370 2250 642 156  5.4872174 5.4916207}%
\put(41.8000,-17.4000){\makebox(0,0)[rb]{$z_0$}}%
%
\special{pn 20}%
\special{sh 1}%
\special{ar 4210 1940 10 10 0  6.28318530717959E+0000}%
\special{sh 1}%
\special{ar 4210 1940 10 10 0  6.28318530717959E+0000}%
%
\special{pn 8}%
\special{pa 3950 1590}%
\special{pa 4000 1590}%
\special{pa 4000 1550}%
\special{pa 3950 1550}%
\special{pa 3950 1590}%
\special{fp}%
%
\special{pn 8}%
\special{pa 4540 1610}%
\special{pa 4590 1610}%
\special{pa 4590 1570}%
\special{pa 4540 1570}%
\special{pa 4540 1610}%
\special{fp}%
%
\special{pn 8}%
\special{pa 3970 2170}%
\special{pa 4010 2170}%
\special{pa 4010 2120}%
\special{pa 3970 2120}%
\special{pa 3970 2170}%
\special{fp}%
%
\special{pn 8}%
\special{pa 4560 2200}%
\special{pa 4600 2200}%
\special{pa 4600 2150}%
\special{pa 4560 2150}%
\special{pa 4560 2200}%
\special{fp}%
%
\special{pn 8}%
\special{pa 4530 1930}%
\special{pa 4580 1930}%
\special{pa 4580 1880}%
\special{pa 4530 1880}%
\special{pa 4530 1930}%
\special{fp}%
%
\special{pn 8}%
\special{pa 3940 1920}%
\special{pa 3990 1920}%
\special{pa 3990 1870}%
\special{pa 3940 1870}%
\special{pa 3940 1920}%
\special{fp}%
%
\special{pn 8}%
\special{ar 4370 1690 590 80  6.2831853 6.2831853}%
\special{ar 4370 1690 590 80  0.0000000 3.1415927}%
%
\special{pn 20}%
\special{sh 1}%
\special{ar 4200 1770 10 10 0  6.28318530717959E+0000}%
\special{sh 1}%
\special{ar 4200 1770 10 10 0  6.28318530717959E+0000}%
%
\special{pn 8}%
\special{pa 4880 1100}%
\special{pa 4740 1760}%
\special{dt 0.045}%
\special{sh 1}%
\special{pa 4740 1760}%
\special{pa 4774 1700}%
\special{pa 4752 1708}%
\special{pa 4734 1692}%
\special{pa 4740 1760}%
\special{fp}%
%
\special{pn 8}%
\special{pa 4570 3240}%
\special{pa 4620 3240}%
\special{pa 4620 3190}%
\special{pa 4570 3190}%
\special{pa 4570 3240}%
\special{fp}%
%
\special{pn 8}%
\special{pa 4570 3470}%
\special{pa 4620 3470}%
\special{pa 4620 3420}%
\special{pa 4570 3420}%
\special{pa 4570 3470}%
\special{fp}%
%
\special{pn 8}%
\special{pa 4590 3700}%
\special{pa 4640 3700}%
\special{pa 4640 3650}%
\special{pa 4590 3650}%
\special{pa 4590 3700}%
\special{fp}%
%
\special{pn 8}%
\special{pa 3670 1510}%
\special{pa 4200 1840}%
\special{dt 0.045}%
\special{sh 1}%
\special{pa 4200 1840}%
\special{pa 4154 1788}%
\special{pa 4156 1812}%
\special{pa 4134 1822}%
\special{pa 4200 1840}%
\special{fp}%
%
\special{pn 8}%
\special{pa 4070 940}%
\special{pa 4120 2190}%
\special{dt 0.045}%
\special{sh 1}%
\special{pa 4120 2190}%
\special{pa 4138 2124}%
\special{pa 4118 2138}%
\special{pa 4098 2124}%
\special{pa 4120 2190}%
\special{fp}%
%
\special{pn 8}%
\special{pa 4070 940}%
\special{pa 4310 1570}%
\special{dt 0.045}%
\special{sh 1}%
\special{pa 4310 1570}%
\special{pa 4306 1502}%
\special{pa 4292 1520}%
\special{pa 4268 1516}%
\special{pa 4310 1570}%
\special{fp}%
\put(14.7000,-14.8000){\makebox(0,0)[rb]{$D^n(R)$}}%
%
\special{pn 8}%
\special{pa 1330 1530}%
\special{pa 1390 2370}%
\special{dt 0.045}%
\special{sh 1}%
\special{pa 1390 2370}%
\special{pa 1406 2302}%
\special{pa 1386 2318}%
\special{pa 1366 2306}%
\special{pa 1390 2370}%
\special{fp}%
%
\special{pn 8}%
\special{pa 1640 1240}%
\special{pa 1770 2310}%
\special{dt 0.045}%
\special{sh 1}%
\special{pa 1770 2310}%
\special{pa 1782 2242}%
\special{pa 1764 2258}%
\special{pa 1742 2246}%
\special{pa 1770 2310}%
\special{fp}%
%
\special{pn 8}%
\special{pa 2890 1440}%
\special{pa 2780 2410}%
\special{dt 0.045}%
\special{sh 1}%
\special{pa 2780 2410}%
\special{pa 2808 2346}%
\special{pa 2786 2358}%
\special{pa 2768 2342}%
\special{pa 2780 2410}%
\special{fp}%
\put(26.1000,-13.7000){\makebox(0,0)[lb]{$r$-curve}}%
%
\special{pn 8}%
\special{ar 4620 3370 140 560  5.0741182 5.1084039}%
\special{ar 4620 3370 140 560  5.2112611 5.2455468}%
\special{ar 4620 3370 140 560  5.3484039 5.3826897}%
\special{ar 4620 3370 140 560  5.4855468 5.5198325}%
\special{ar 4620 3370 140 560  5.6226897 5.6569754}%
\special{ar 4620 3370 140 560  5.7598325 5.7941182}%
\special{ar 4620 3370 140 560  5.8969754 5.9312611}%
\special{ar 4620 3370 140 560  6.0341182 6.0684039}%
\special{ar 4620 3370 140 560  6.1712611 6.2055468}%
\special{ar 4620 3370 140 560  6.3084039 6.3426897}%
\special{ar 4620 3370 140 560  6.4455468 6.4798325}%
\special{ar 4620 3370 140 560  6.5826897 6.6169754}%
\special{ar 4620 3370 140 560  6.7198325 6.7541182}%
\special{ar 4620 3370 140 560  6.8569754 6.8912611}%
\special{ar 4620 3370 140 560  6.9941182 7.0284039}%
\special{ar 4620 3370 140 560  7.1312611 7.1655468}%
\special{ar 4620 3370 140 560  7.2684039 7.3026897}%
\special{ar 4620 3370 140 560  7.4055468 7.4398325}%
\end{picture}%
\hspace{2.55truecm}}

\vspace{0.5truecm}

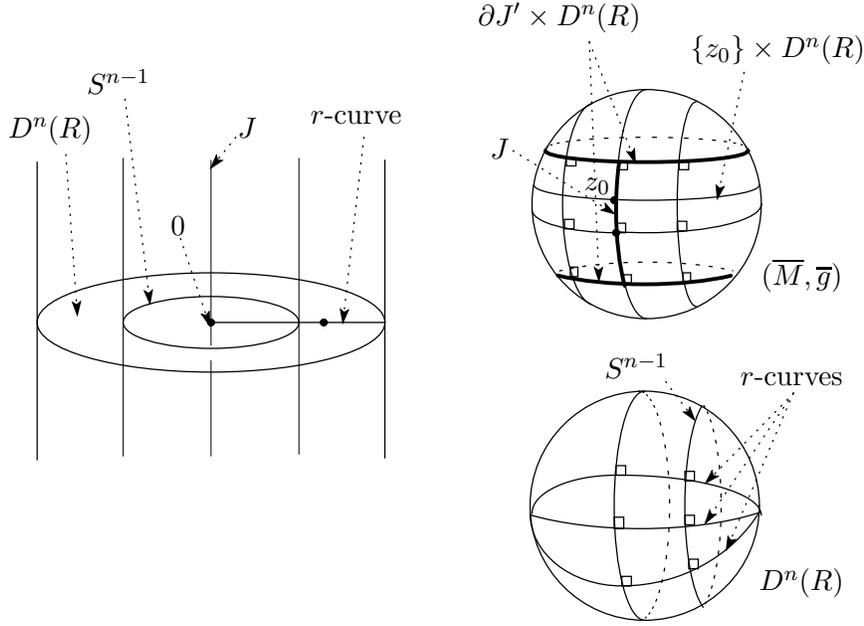
\centerline{{\bf Figure 3.}}

\vspace{0.5truecm}

Let $F$ and $F^{\perp}_{x_{\ast}}$ be as in the above setting $(S)$, 
$I$ the radial geodesic in $F^{\perp}_{x_{\ast}}$ starting from $x_{\ast}$ removed 
the closure of a neighborhood of $x_{\ast}$ and $S^{n-1}$ the unit geodesic sphere centered at $x_{\ast}$ in 
$F^{\perp}_{x_{\ast}}$.  Then $\overline M$ in the setting (S) is 
locally diffeomorphic to $B\times I\times S^{n-1}$ and $\overline g$ in the setting (S) is locally described as 
$\overline g_{(x,r,v)}=A_{r,x}^{\ast}(g_F)_x+dr^2+{\widehat A}_{r,v}^{\ast}(g_{S^{n-1}})_v$ 
($(x,r,v)\in B\times I\times S^{n-1}$), where $g_F,\,g_{S^{n-1}}$ are the induced metrics 
on $F$ and $S^{n-1}$ (where we note that $g_{S^{n-1}}$ is not the canonical metric) and 
$A_{r,x}$ (resp. $\widehat A_{r,v}$) is a $(1,1)$-tensor of $T_xF$ (resp. $T_vS^{n-1}$) depending only on $r$ 
(which are independent of $x$ and $v$ essentially).  
We consider the volume-preserving mean curvature flow starting from a tube over $B$.  
Thus our setting (S) includes Setting([CM2]) (see Figure 4).  
Under the above setting (S), 
we consider the volume-preserving mean curvature flow $f_t$ ($t\in[0,T)$) starting from $f$ and satisfying 
the following boundary condition:

\vspace{0.15truecm}

{\rm (C1)} ${\rm grad}\,r_t=0$ holds along $\partial B$ for all $t\in[0,T)$, where $r_t$ is the radius function 
of $M_t:=f_t(M)$ (i.e., $M_t=\exp^{\perp}(t_{r_t}(B))$), where $r_t$ is possible to be multi-valued.  

\vspace{0.15truecm}

\noindent
Furthermore, assume that the this flow satisfies the following condition:

\vspace{0.15truecm}

{\rm (C2)} $({\rm grad}\,r_t)_x$ belongs to a common eigenspace of the family 
$\{R(\cdot,\xi)\xi\}_{\xi\in T^{\perp}_xB}$ for all $(x,t)\in B\times[0,T)$.  

\vspace{0.5truecm}

\centerline{
\unitlength 0.1in
%
\hspace{1.55truecm}}

\vspace{0.4truecm}

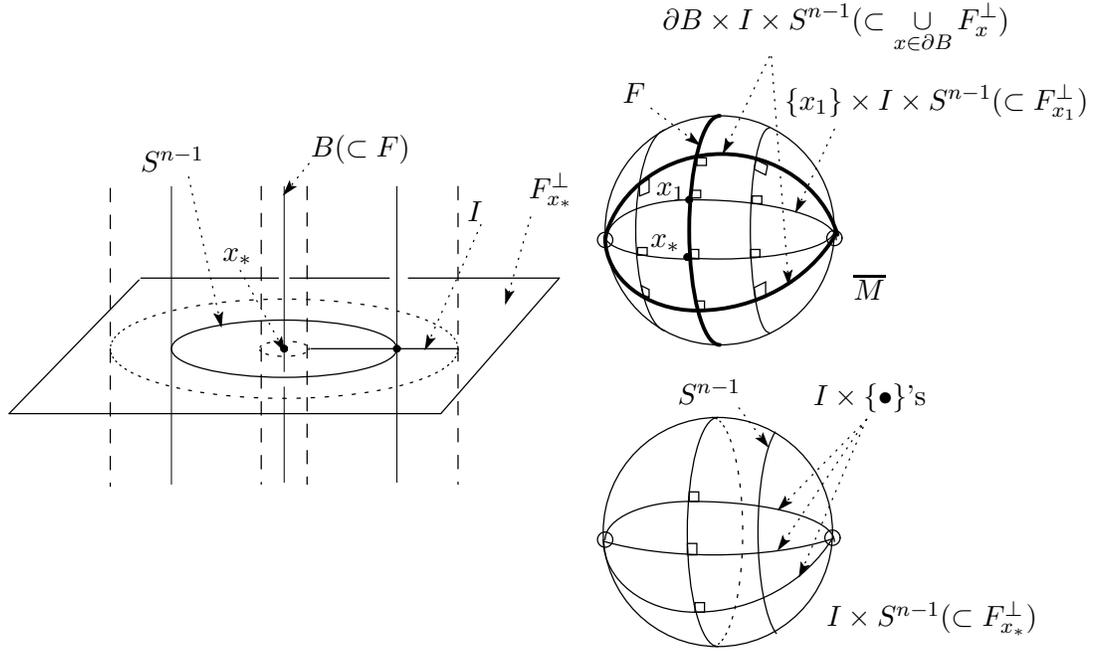
\centerline{{\bf Figure 4.}}

\vspace{0.5truecm}

\noindent

\vspace{0.25truecm}

\noindent
If the initial radius function $r_0$ is constant over a collar neighborhood $U$ of $\partial B$, then 
$t_{r_0}(U)$ is of constant mean curvature because the normal umbrellas of $F$ are of rank one by the assumption 
(furthermore, these umbrellas are automatically isometric to one another).  Hence 
$\overline H_t-H_t$ ($t\in[0,T)$) remain to be constant over $t_{r_0}(U)$, that is, 
$r_t$ ($t\in[0,T)$) remain to be constant over $U$ (see Figure 5).  
If $\overline M$ is a rank one symmetric space (other than a sphere and a (real) hyperbolic space) and if $F$ 
is an invariant submanifold, then $R(\cdot,\xi)\xi\vert_{T_xF}$ $(x\in F,\,\,\xi\in T_x^{\perp}F)$ are 
the constant-multiple of the identity transformation ${\rm id}_{T_xF}$ of $T_xF$ and hence 
the condition (C2) automatically holds.  Here ``invariant submanifold'' means that the tangent spaces of 
the submanifold are invariant by $J$, where $J$ denotes the complex structure of $\overline M$ 
in the case where $\overline M$ is the complex projective space or the complex hyperbolic space, 
any complex structure belonging to the quaternionic structure of $\overline M$ in the case where $\overline M$ 
is the quaternionic projective space or the quaternionic hyperbolic space, and any complex structure belonging 
to the Cayley structure of $\overline M$ in the case where $\overline M$ is the Cayley plane 
(i.e., $F_4/Spin(9)$) or the Cayley hyperbolic space (i.e., $F_4^{-20}/Spin(9)$).  
Also, if $F\subset\overline M$ is one of meridians in irreducible rank two symmetric spaces of compact type 
in Table 2 (see Appendix) and if ${\cal D}$ is the corresponding distribution on $F$ as in Table 2, then 
${\cal D}_x$ is a common eigenspace of the family $\{R(\cdot,\xi)\xi\,\vert\,\xi\in T^{\perp}_xF\}$ for any 
$x\in F$.  It is easy to show that, if the initial radius function $r_0$ satisfies 
$({\rm grad}\,r_0)_x\in{\cal D}_x$ ($x\in B$), then $r_t$ satisfies $({\rm grad}\,r_t)_x\in{\cal D}_x$ 
($x\in B$) for all $t\in[0,T)$, that is, the condition (C2) holds.  

\vspace{0.5truecm}

\centerline{
\unitlength 0.1in
\begin{picture}( 29.0100, 22.3000)(  7.6500,-29.4000)
\put(34.8000,-15.4000){\makebox(0,0)[lt]{$f(M)$}}%
%
\special{pn 20}%
\special{sh 1}%
\special{ar 2536 1350 10 10 0  6.28318530717959E+0000}%
\special{sh 1}%
\special{ar 2536 1350 10 10 0  6.28318530717959E+0000}%
%
\special{pn 20}%
\special{sh 1}%
\special{ar 2548 2530 10 10 0  6.28318530717959E+0000}%
\special{sh 1}%
\special{ar 2548 2530 10 10 0  6.28318530717959E+0000}%
%
\special{pn 8}%
\special{pa 1880 1090}%
\special{pa 1344 1564}%
\special{pa 3270 1564}%
\special{pa 3666 1090}%
\special{pa 3666 1090}%
\special{pa 1880 1090}%
\special{fp}%
%
\special{pn 4}%
\special{pa 2530 2940}%
\special{pa 2530 902}%
\special{fp}%
\put(26.0000,-8.8000){\makebox(0,0)[rb]{$F$}}%
\put(10.3500,-21.2100){\makebox(0,0)[rb]{$U$}}%
\put(12.2000,-11.6000){\makebox(0,0)[rb]{$P$}}%
%
\special{pn 8}%
\special{pa 1870 2330}%
\special{pa 1334 2804}%
\special{pa 3260 2804}%
\special{pa 3656 2330}%
\special{pa 3656 2330}%
\special{pa 1870 2330}%
\special{fp}%
%
\special{pn 20}%
\special{pa 2530 1360}%
\special{pa 2530 1630}%
\special{fp}%
%
\special{pn 20}%
\special{pa 2542 2520}%
\special{pa 2542 2260}%
\special{fp}%
%
\special{pn 20}%
\special{ar 2530 1350 410 110  0.0000000 6.2831853}%
%
\special{pn 20}%
\special{ar 2542 2540 410 110  0.0000000 6.2831853}%
%
\special{pn 20}%
\special{pa 2120 1360}%
\special{pa 2120 1630}%
\special{fp}%
%
\special{pn 20}%
\special{pa 2940 1350}%
\special{pa 2940 1620}%
\special{fp}%
%
\special{pn 20}%
\special{pa 2132 2520}%
\special{pa 2132 2260}%
\special{fp}%
%
\special{pn 20}%
\special{pa 2940 2540}%
\special{pa 2940 2280}%
\special{fp}%
%
\special{pn 20}%
\special{ar 1806 1640 302 160  6.2831853 6.2831853}%
\special{ar 1806 1640 302 160  0.0000000 0.9621521}%
%
\special{pn 20}%
\special{ar 2060 1960 180 210  3.1415927 4.2545978}%
%
\special{pn 20}%
\special{ar 1830 2270 300 148  5.3231149 6.2831853}%
%
\special{pn 20}%
\special{ar 2090 1960 214 212  2.0258971 3.1415927}%
%
\special{pn 8}%
\special{ar 3140 2270 142 204  3.1415927 4.0997661}%
%
\special{pn 8}%
\special{ar 3010 1950 68 204  6.2831853 6.2831853}%
\special{ar 3010 1950 68 204  0.0000000 1.1188669}%
%
\special{pn 20}%
\special{sh 1}%
\special{ar 2526 1348 10 10 0  6.28318530717959E+0000}%
\special{sh 1}%
\special{ar 2526 1348 10 10 0  6.28318530717959E+0000}%
%
\special{pn 20}%
\special{sh 1}%
\special{ar 2538 2528 10 10 0  6.28318530717959E+0000}%
\special{sh 1}%
\special{ar 2538 2528 10 10 0  6.28318530717959E+0000}%
%
\special{pn 20}%
\special{pa 2520 1360}%
\special{pa 2520 1638}%
\special{fp}%
%
\special{pn 20}%
\special{pa 2532 2518}%
\special{pa 2532 2256}%
\special{fp}%
%
\special{pn 8}%
\special{ar 2530 1336 474 176  0.0000000 6.2831853}%
%
\special{pn 8}%
\special{pa 2060 1350}%
\special{pa 2060 1628}%
\special{fp}%
%
\special{pn 8}%
\special{pa 3000 1360}%
\special{pa 3000 1680}%
\special{fp}%
%
\special{pn 8}%
\special{pa 2060 2520}%
\special{pa 2060 2260}%
\special{fp}%
%
\special{pn 8}%
\special{pa 3000 2528}%
\special{pa 3000 2266}%
\special{fp}%
%
\special{pn 8}%
\special{ar 2530 2536 474 176  0.0000000 6.2831853}%
%
\special{pn 20}%
\special{ar 3230 1610 302 160  2.1794405 3.1415927}%
%
\special{pn 20}%
\special{ar 2966 1940 198 232  5.1707246 6.2831853}%
%
\special{pn 20}%
\special{ar 3198 2260 258 160  3.1415927 4.1014896}%
%
\special{pn 20}%
\special{ar 2958 1960 216 190  6.2831853 6.2831853}%
\special{ar 2958 1960 216 190  0.0000000 1.1144479}%
%
\special{pn 8}%
\special{ar 1936 1640 126 172  6.2831853 6.2831853}%
\special{ar 1936 1640 126 172  0.0000000 0.9633078}%
%
\special{pn 8}%
\special{ar 2038 1970 68 204  3.1415927 4.2642440}%
%
\special{pn 8}%
\special{ar 1880 2280 180 204  5.3231149 6.2831853}%
%
\special{pn 8}%
\special{ar 2038 1970 68 204  2.0227257 3.1415927}%
%
\special{pn 20}%
\special{ar 2530 1640 410 110  0.0000000 6.2831853}%
%
\special{pn 8}%
\special{ar 2530 1650 474 176  0.0000000 6.2831853}%
%
\special{pn 8}%
\special{ar 2530 2250 474 176  0.0000000 6.2831853}%
%
\special{pn 20}%
\special{ar 2540 2250 410 110  0.0000000 6.2831853}%
%
\special{pn 8}%
\special{pa 1080 2020}%
\special{pa 2520 1590}%
\special{dt 0.045}%
\special{sh 1}%
\special{pa 2520 1590}%
\special{pa 2450 1590}%
\special{pa 2470 1606}%
\special{pa 2462 1628}%
\special{pa 2520 1590}%
\special{fp}%
%
\special{pn 8}%
\special{pa 1080 2020}%
\special{pa 2520 2290}%
\special{dt 0.045}%
\special{sh 1}%
\special{pa 2520 2290}%
\special{pa 2458 2258}%
\special{pa 2468 2280}%
\special{pa 2452 2298}%
\special{pa 2520 2290}%
\special{fp}%
\put(34.8000,-17.4000){\makebox(0,0)[lt]{$f_t(M)$}}%
%
\special{pn 8}%
\special{pa 3450 1670}%
\special{pa 3140 1820}%
\special{dt 0.045}%
\special{sh 1}%
\special{pa 3140 1820}%
\special{pa 3210 1810}%
\special{pa 3188 1798}%
\special{pa 3192 1774}%
\special{pa 3140 1820}%
\special{fp}%
%
\special{pn 8}%
\special{pa 3450 1860}%
\special{pa 3080 1980}%
\special{dt 0.045}%
\special{sh 1}%
\special{pa 3080 1980}%
\special{pa 3150 1978}%
\special{pa 3132 1964}%
\special{pa 3138 1940}%
\special{pa 3080 1980}%
\special{fp}%
%
\special{pn 8}%
\special{pa 1180 1180}%
\special{pa 1690 1450}%
\special{dt 0.045}%
\special{sh 1}%
\special{pa 1690 1450}%
\special{pa 1640 1402}%
\special{pa 1644 1426}%
\special{pa 1622 1436}%
\special{pa 1690 1450}%
\special{fp}%
%
\special{pn 8}%
\special{pa 1180 1180}%
\special{pa 1620 2700}%
\special{dt 0.045}%
\special{sh 1}%
\special{pa 1620 2700}%
\special{pa 1622 2630}%
\special{pa 1606 2650}%
\special{pa 1582 2642}%
\special{pa 1620 2700}%
\special{fp}%
%
\special{pn 8}%
\special{ar 3026 1950 54 182  5.1760366 6.2831853}%
%
\special{pn 8}%
\special{ar 3140 1660 150 160  2.1770340 3.1415927}%
\end{picture}%
\hspace{1.55truecm}}

\vspace{0.5truecm}

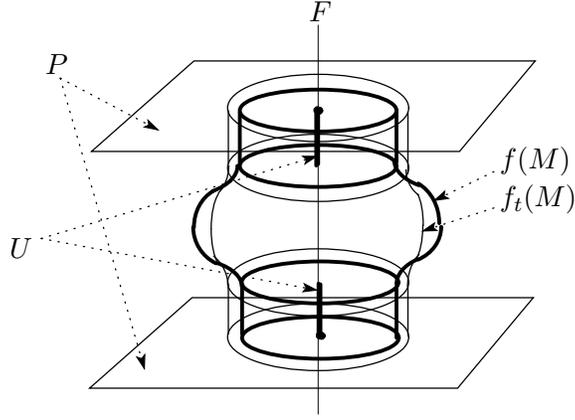
\centerline{{\bf Figure 5.}}

\vspace{0.5truecm}

\noindent
If the condition (C2) holds, then it is shown that there uniquely exists the volume-preserving mean curvature 
flow $f_t:M\hookrightarrow\overline M$ starting from $f$ as in the above setting (S) and satisfying the condition 
(C1) in short time (see Proposition 4.2).  
Under these assumptions, we first derive the evolution equations for the radius functions of the flow and 
some quantities related to the gradients of the functions (see Sections 4 and 5).  
Next, in the following special case, we derive the following preservability theorem for the tubeness along 
the flow by using the evolution equations.  

\vspace{0.5truecm}

\noindent
{\bf Theorem A.} {\sl Let $f$ be as in the above setting (S) and $f_t$ ($t\in[0,T)$) the volume-preserving 
mean curvature flow starting from $f$ and satisfying the above condition (C1).  Assume that $\overline M$ is 
a rank one symmetric space of non-compact type, $F$ is an invariant 
submanifold and that $B$ is a closed geodesic ball of radius $r_B$ centered at $x_{\ast}$ in $F$, where 
the invariantness of $F$ means the totally geodesicness 
in the case where $\overline M$ is a (real) hyperbolic space.  If $r_0$ is radial with respect to $x_{\ast}$ 
(i.e., $r_0$ is constant along each geodesic sphere centered at $x_{\ast}$ in $F$), then 
$M_t$ ($t\in[0,T)$) remain to be tubes over $B$ such that the volume of the closed domain 
surrounded by $M_t$ and $P$ is equal to ${\rm Vol}(D)$.}

\vspace{0.5truecm}

Furthermore, we derive the following results.   

\vspace{0.5truecm}

\noindent
{\bf Theorem B.} {\sl Under the hypothesis of Theorem A, 
one of the following statements {\rm(a)} and {\rm(b)} holds:

{\rm (a)} $M_t:=f_t(M)$ reaches $B$ as $t\to T$, 


{\rm (b)} $T=\infty$ and $M_t$ converges to a tube of constant mean curvature over $B$ 
(in $C^{\infty}$-topology) as $t\to\infty$.}

\vspace{0.5truecm}

\centerline{
\unitlength 0.1in
\begin{picture}( 26.5000, 16.7900)( 15.6000,-27.7000)
%
\special{pn 8}%
\special{ar 2010 1178 438 86  0.0000000 6.2831853}%
%
\special{pn 8}%
\special{ar 2020 1168 450 268  2.3087784 3.1415927}%
%
\special{pn 8}%
\special{ar 1592 1446 178 110  5.3981185 6.2831853}%
%
\special{pn 8}%
\special{pa 1770 1448}%
\special{pa 1770 1488}%
\special{fp}%
%
\special{pn 8}%
\special{ar 2010 1178 450 268  6.2831853 6.2831853}%
\special{ar 2010 1178 450 268  0.0000000 0.8303019}%
%
\special{pn 8}%
\special{ar 2438 1454 178 110  3.1415927 4.0266595}%
%
\special{pn 8}%
\special{pa 2260 1446}%
\special{pa 2260 1486}%
\special{fp}%
%
\special{pn 8}%
\special{ar 2006 1790 444 92  0.0000000 6.2831853}%
%
\special{pn 8}%
\special{ar 2010 1792 450 268  3.1415927 3.9726533}%
%
\special{pn 8}%
\special{ar 1592 1516 178 112  6.2831853 6.2831853}%
\special{ar 1592 1516 178 112  0.0000000 0.8804317}%
%
\special{pn 8}%
\special{pa 1770 1522}%
\special{pa 1770 1480}%
\special{fp}%
%
\special{pn 8}%
\special{pa 2260 1524}%
\special{pa 2260 1484}%
\special{fp}%
%
\special{pn 8}%
\special{ar 3690 1178 516 84  0.0000000 6.2831853}%
%
\special{pn 8}%
\special{ar 3696 1790 516 84  0.0000000 6.2831853}%
%
\special{pn 8}%
\special{pa 3702 1436}%
\special{pa 3702 1530}%
\special{fp}%
%
\special{pn 8}%
\special{ar 2020 1430 240 24  0.0000000 6.2831853}%
%
\special{pn 8}%
\special{ar 2020 1532 252 32  0.0000000 6.2831853}%
%
\special{pn 13}%
\special{pa 2636 1500}%
\special{pa 3086 1500}%
\special{fp}%
\special{sh 1}%
\special{pa 3086 1500}%
\special{pa 3018 1480}%
\special{pa 3032 1500}%
\special{pa 3018 1520}%
\special{pa 3086 1500}%
\special{fp}%
%
\special{pn 8}%
\special{ar 2020 2106 440 88  0.0000000 6.2831853}%
%
\special{pn 8}%
\special{ar 2030 2100 450 268  2.3105321 3.1415927}%
%
\special{pn 8}%
\special{ar 1602 2374 178 112  5.4055539 6.2831853}%
%
\special{pn 8}%
\special{ar 2020 2106 450 268  6.2831853 6.2831853}%
\special{ar 2020 2106 450 268  0.0000000 0.8301444}%
%
\special{pn 8}%
\special{ar 2448 2384 178 110  3.1415927 4.0254582}%
%
\special{pn 8}%
\special{ar 2030 2374 240 24  0.0000000 6.2831853}%
%
\special{pn 13}%
\special{pa 2646 2430}%
\special{pa 3096 2430}%
\special{fp}%
\special{sh 1}%
\special{pa 3096 2430}%
\special{pa 3028 2410}%
\special{pa 3042 2430}%
\special{pa 3028 2450}%
\special{pa 3096 2430}%
\special{fp}%
%
\special{pn 8}%
\special{ar 2016 2670 434 102  0.0000000 6.2831853}%
%
\special{pn 8}%
\special{ar 1968 2660 386 268  3.1415927 3.9740871}%
%
\special{pn 8}%
\special{ar 1592 2374 178 124  6.2831853 6.2831853}%
\special{ar 1592 2374 178 124  0.0000000 0.8819825}%
%
\special{pn 8}%
\special{ar 3694 2122 356 84  0.0000000 6.2831853}%
%
\special{pn 8}%
\special{ar 3694 2680 356 82  0.0000000 6.2831853}%
%
\special{pn 8}%
\special{pa 3340 2120}%
\special{pa 3340 2672}%
\special{fp}%
%
\special{pn 8}%
\special{pa 4048 2118}%
\special{pa 4048 2668}%
\special{fp}%
\put(26.5600,-15.7900){\makebox(0,0)[lt]{$t\to T$}}%
\put(26.7000,-25.3000){\makebox(0,0)[lt]{$t\to\infty$}}%
%
\special{pn 8}%
\special{ar 2070 2680 386 268  5.4506908 6.2831853}%
%
\special{pn 8}%
\special{ar 2448 2380 178 124  2.2596101 3.1415927}%
%
\special{pn 8}%
\special{ar 2010 1800 450 268  5.4521247 6.2831853}%
%
\special{pn 8}%
\special{ar 2438 1520 178 112  2.2611610 3.1415927}%
%
\special{pn 8}%
\special{ar 3700 1190 530 170  1.8009246 3.1415927}%
%
\special{pn 8}%
\special{ar 3560 1442 142 92  4.7885771 6.1588303}%
%
\special{pn 8}%
\special{ar 3744 1180 468 168  6.2831853 6.2831853}%
\special{ar 3744 1180 468 168  0.0000000 1.2657161}%
%
\special{pn 8}%
\special{ar 3890 1420 190 80  3.1796695 4.6486997}%
%
\special{pn 8}%
\special{ar 3690 1790 510 190  3.1415927 4.4112734}%
%
\special{pn 8}%
\special{ar 3520 1520 180 90  6.2831853 6.2831853}%
\special{ar 3520 1520 180 90  0.0000000 1.5374753}%
%
\special{pn 8}%
\special{ar 3680 1790 526 182  5.0837487 6.2831853}%
%
\special{pn 8}%
\special{ar 3930 1530 220 90  1.7330241 3.1415927}%
\end{picture}%
\hspace{0truecm}}

\vspace{0.5truecm}

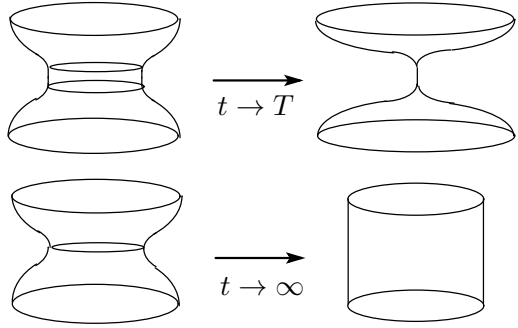
\centerline{{\bf Figure 6.}}

\vspace{0.5truecm}

\noindent
{\bf Theorem C.} {\sl 
Under the hypothesis of Theorem A, assume that 
$${\rm Vol}(M_0)\leq v_{m^H-1}v_{m^V}(\delta_2\circ\delta_1^{-1})
\left(\frac{{\rm Vol}(D)}{v_{m^V}{\rm Vol}(B)}\right),$$
where $m^H:={\rm dim}\,F$, $m^V:={\rm codim}\,F-1$, $v_{m_H-1}$ (resp. $v_{m^V}$) is the volume of the 
$m^H-1$ (resp. $m^V$)-dimensional Euclidean unit sphere and $\delta_i$ ($i=1,2$) are increasing functions over 
${\Bbb R}$ explicitly described (see Section 6).  
Then $T=\infty$ and $M_t$ converges to a tube of constant mean curvature over $B$ (in $C^{\infty}$-topology) 
as $t\to\infty$.}


\vspace{0.5truecm}

\noindent
{\it Remark 1.2.} Let $\overline g$ be the metric of $\overline M$ and $c$ a positive constant.  
As $c\to\infty$, $c\overline g$ approches to a flat metric and $\delta_i$ ($i=1,2$) approches to the identity 
transformation of $[0,\infty)$ and hence 
$\displaystyle{v_{m^V}(\delta_2\circ\delta_1^{-1})\left(\frac{{\rm Vol}(D)}{v_{m^V}{\rm Vol}(B)}\right)}$ 
approaches to $\displaystyle{\frac{{\rm Vol}(D)}{{\rm Vol}(B)}}$.  
Thus, as $c\to\infty$, the condition 
$\displaystyle{{\rm Vol}(M_0)\leq v_{m^H-1}v_{m^V}(\delta_2\circ\delta_1^{-1})
\left(\frac{{\rm Vol}(D)}{v_{m^V}{\rm Vol}(B)}\right)}$ approaches to the condition 
$\displaystyle{{\rm Vol}(M)\leq\frac{{\rm Vol}(D)}{d}}$ in the statement (ii) of Fact 1 in the case of 
${\rm dim}\,F=1$.  

\vspace{0.5truecm}

In the future, we plan to tackle the following problem.  

\vspace{0.5truecm}

\noindent
{\it Problem.} {\sl Under the hypothesis of Theorem C, does $M_t$ converge to a tube of constant radius over 
$B$ (in $C^{\infty}$-topology) as $t\to\infty$?}

\vspace{0.5truecm}

As the first step to solve this problem, I need to classify tubes of constant mean curvature over $F$ other than 
tubes of constant radius over $F$.  

\section{The mean curvature of a tube over a reflective submanifold} 
In this section, we shall calculate the mean curvature of a tube 
over a reflective submanifold in a symmetric space of compact type or non-compact type.  
Let $\overline M=G/K$ be a symmetric space of compact type or non-compact type, where $G$ 
is the identity component of the isometry group of $\overline M$ and $K$ is the isotropy 
group of $G$ at some point $p_0$ of $\overline M$.  
Let $F$ be a reflective submanifold in $\overline M$ such that 
the nomal umbrellas $\Sigma_x$'s ($x\in F$) are symmetric spaces of rank one.  
Denote by $\bar g$ (resp. $g_F$) the Rimannian metric of $\overline M$ (resp. $F$) and 
$\overline{\nabla}$ (resp. $\nabla^F$) the Riemannian connection of $\overline M$ (resp. $F$).  
Denote by $\mathfrak g$ and $\mathfrak k$ the Lie algebras of 
$G$ and $K$, respectively.  Also, let $\theta$ be the Cartan involution of 
$\mathfrak g$ with $({\rm Fix}\,\theta)_0\subset K\subset{\rm Fix}\,\theta$ and 
set $\mathfrak p:={\rm Ker}(\theta+{\rm id})$, which is identified with the tangent space 
of $T_{p_0}\overline M$ of $\overline M$ at $p_0$.  
Without loss of generality, we may assume that $p_0$ belongs to $F$.  
Set $\mathfrak p':=T_{p_0}F$ and ${\mathfrak p'}^{\perp}:=T_{p_0}^{\perp}F$.  
Take a maximal abelian subspace $\mathfrak b$ of ${\mathfrak p'}^{\perp}$ and 
a maximal abelian subspace $\mathfrak a$ of $\mathfrak p$ including $\mathfrak b$.  
Note that the dimension of $\mathfrak b$ is equal to $1$ because the normal umbrellas of 
$F$ is symmetric spaces of rank one by the assumption.  
For each $\alpha\in\mathfrak a^{\ast}$ and $\beta\in\mathfrak b^{\ast}$, we define 
a subspace $\mathfrak p_{\alpha}$ and $\mathfrak p_{\beta}$ of $\mathfrak p$ by 
$$\mathfrak p_{\alpha}:=\{Y\in\mathfrak p\,\vert\,{\rm ad}(X)^2(Y)=-\varepsilon\alpha(X)^2Y
\,\,{\rm for}\,\,{\rm all}\,\,X\in\mathfrak a\}$$
and
$$\mathfrak p_{\beta}:=\{Y\in\mathfrak p\,\vert\,{\rm ad}(X)^2(Y)=-\varepsilon\beta(X)^2Y
\,\,{\rm for}\,\,{\rm all}\,\,X\in\mathfrak b\},$$
respectively, where ${\rm ad}$ is the adjoint representation of $\mathfrak g$, 
$\mathfrak a^{\ast}$ (resp. $\mathfrak b^{\ast}$) is the dual space of 
$\mathfrak a$ (resp. $\mathfrak b$) and 
$\varepsilon$ is given by 
$$\varepsilon:=\left\{
\begin{array}{ll}
1 & ({\rm when}\,\,\overline M\,\,{\rm is}\,\,{\rm of}\,\,{\rm compact}\,\,{\rm type})\\
-1 & ({\rm when}\,\,\overline M\,\,{\rm is}\,\,{\rm of}\,\,
{\rm non-compact}\,\,{\rm type}).
\end{array}
\right.$$
Define a subset $\triangle$ of $\mathfrak a^{\ast}$ by 
$$\triangle:=\{\alpha\in\mathfrak a^{\ast}\,\vert\,\mathfrak p_{\alpha}\not=\{0\}\},$$
and subsets $\triangle'$ and $\triangle'_V$ of $\mathfrak b^{\ast}$ by 
$$\triangle':=\{\beta\in\mathfrak b^{\ast}\,\vert\,\mathfrak p_{\beta}\not=\{0\}\}$$
and 
$$\triangle'_V:=\{\beta\in\mathfrak b^{\ast}\,\vert\,\mathfrak p_{\beta}\cap\mathfrak p'^{\perp}\not=\{0\}\}.$$
The systems $\triangle$ and $\triangle'_V$ are root systems and 
$\triangle'=\{\alpha\vert_{\mathfrak b}\,\vert\,\alpha\in\triangle\}$ holds.  
Let $\triangle_+$ (resp. $(\triangle'_V)_+$) be the positive root system of 
$\triangle$ (resp. $\triangle'_V$) with respect to 
some lexicographic ordering of $\mathfrak a^{\ast}$ (resp. $\mathfrak b^{\ast}$) and 
$\triangle'_+$ be the positive subsystem of $\triangle'$ with respect to the lexicographic ordering of 
$\mathfrak b^{\ast}$, where we take one compatible with the lexicographic ordering of $\mathfrak b^{\ast}$ 
as the lexicographic ordering of $\mathfrak a^{\ast}$.  
Also we have the following root space decomposition:
$$\mathfrak p=\mathfrak a+\sum_{\alpha\in\triangle_+}\mathfrak p_{\alpha}
=\mathfrak z_{\mathfrak p}(\mathfrak b)+\sum_{\beta\in\triangle'_+}\mathfrak p_{\beta},$$
where $\mathfrak z_{\mathfrak p}(\mathfrak b)$ is the centralizer of $\mathfrak b$ in $\mathfrak p$.  
For convenience, we set $\mathfrak p_0:=\mathfrak z_{\mathfrak p}(\mathfrak b)$.  
Since the normal umbrellas of $F$ are symmetric spaces 
of rank one, ${\rm dim}\,\mathfrak b=1$ and this root system $\triangle'_V$ is of 
$(\mathfrak a_1)$-type or $(\mathfrak b\mathfrak d_1)$-type.  
Hence $(\triangle'_V)_+$ is described as 
$$(\triangle'_V)_+=\left\{
\begin{array}{ll}
\{\beta\} & {\rm(}\triangle'_V\,:\,(\mathfrak a_1){\rm -type}{\rm)}\\
\{\beta,2\beta\} & {\rm(}\triangle'_V\,:\,(\mathfrak{bd}_1){\rm -type}{\rm)}
\end{array}\right.$$
for some $\beta(\not=0)\in\mathfrak b^{\ast}$.  
However, in general, we may describe as 
$(\triangle'_V)_+=\{\beta,2\beta\}$ by interpretting as $\mathfrak p_{2\beta}=\{0\}$ 
when $\triangle'_V$ is of $(\mathfrak a_1)$.  
The system $\triangle'_+$ is described as 
$$\triangle'_+=\{k\beta\,\vert\,k\in{\cal K}\}$$ 
for some finite subset ${\cal K}$ of ${\Bbb R}_+$.  
Set $b:=\vert\beta(X_0)\vert$ for a unit vector $X_0$ of $\mathfrak b$.  
Since $F$ is curvature-adapted, $\mathfrak p'$ and ${\mathfrak p'}^{\perp}$ are 
${\rm ad}(X)^2$-invariant for each $X\in\mathfrak b$.  
Hence we have the following direct sum decompositions:
$$\mathfrak p'=\mathfrak p_0\cap\mathfrak p'
+\sum_{k\in{\cal K}}(\mathfrak p_{k\beta}\cap\mathfrak p')$$
and 
$$(\mathfrak p')^{\perp}=\mathfrak b
+\sum_{k=1}^2(\mathfrak p_{k\beta}\cap(\mathfrak p')^{\perp}).$$
For each $p\in\overline M$, we choose a shortest geodesic $\gamma_{p_0p}:[0,1]\to\overline M$ 
with $\gamma_{p_0p}(0)=p_0$ and $\gamma_{p_0p}(1)=p$, 
where the choice of $\gamma_{p_0p}$ is not unique in the case where $p$ belongs to the cut locus 
of $p_0$.  Denote by $\tau_p$ the parallel translation along $\gamma_{p_0p}$.  
For $w\in T_p\overline M$, we define linear transformations ${\rm ad}(w)^2,\,D^{id}_w,\,D^{co}_w$ and 
$D^{si}_w$ of $T_p\overline M$ by 
$$\begin{array}{l}
\displaystyle{{\rm ad}(w)^2:=\tau_p\circ{\rm ad}(\tau_p^{-1}w)^2\circ \tau_p^{-1},}\\
\displaystyle{D^{co}_w:=\tau_p\circ\cos({\bf i}{\rm ad}(\tau_p^{-1}w))\circ \tau_p^{-1}}
\end{array}$$
and 
$$D^{si}_w:=\tau_p{\ast}\circ\frac{\sin({\bf i}{\rm ad}(\tau_p^{-1}w))}
{{\bf i}{\rm ad}(\tau_p^{-1}w)}\circ \tau_p^{-1},$$
respectively, where ${\bf i}$ is the imaginary unit.  
Note that, if $\tau_p^{-1}(X)\in\mathfrak p_{k\beta}$, then 
$$\begin{array}{c}
\displaystyle{{\rm ad}(w)^2(X)=-(\sqrt{\varepsilon}kb\vert\vert X\vert\vert)^2X,\quad\,\,
D^{co}_w(X)=\cos(\sqrt{\varepsilon}kb\vert\vert X\vert\vert)(X)}\\
\displaystyle{{\rm and}\quad\,\,D^{si}_w=\frac{\sin(\sqrt{\varepsilon}kb\vert\vert X\vert\vert)}
{\sqrt{\varepsilon}kb\vert\vert X\vert\vert}X.}
\end{array}
\leqno{(2.1)}$$
Let $r$ be a positive-valued function over $F$ with $r<r_F$ and 
$B$ a compact closed domain in $F$  Set $M:=t_r(B)$ and $f:=\exp^{\perp}\vert_{t_r(B)}$.  
Then $f$ is an embedding.  Denote by $g$ the induced metric on $M$, 
$N$ the outward unit normal vector field of $M$ and $A$ the shape operator of 
$M$ with respect to $-N$.  
Fix $x\in B$ and $\xi\in M\cap T^{\perp}_xB$.  
Without loss of generality, we may assume that $\tau_x^{-1}\xi\in\mathfrak b$.  
Denote by $\gamma_{\xi}$ the normal geodesic of $B$ whose initial vector is equal to $\xi$ 
and $\tau_{\gamma_{\xi}}$ the parallel translation along $\gamma_{\xi}\vert_{[0,1]}$.  
The vertical subspace ${\cal V}_{\xi}$ and the horizontal subspace ${\cal H}_{\xi}$ at 
$\xi$ are defined by 
${\cal V}_{\xi}:=T_{\xi}(M\cap T^{\perp}_xB)$ and 
${\cal H}_{\xi}:=\{\widetilde X_{\xi}\,\vert\,X\in T_xB\}$, respectively, where 
$\widetilde X_{\xi}$ is the natural lift of $X$ to $\xi$ (see [Ko1] about 
the definition of the natural lift).  
Take $v\in{\cal V}_{\xi}$.  
Let $J_{\xi,v}$ be the Jacobi field along $\gamma_{\xi}$ with $J_{\xi,v}(0)=0$ and $J'_{\xi,v}(0)=v$, where 
$v$ is regarded as an element of $T_x^{\perp}B$ under the identification of $T_{\xi}(T_x^{\perp}B)$ and $T_xB$.  
According to $(1.2)$ in [Ko1], we obtain 
$$
J_{\xi,v}(s)=\tau_{\gamma_{\xi}\vert_{[0,s]}}
(D^{co}_{s\xi}(J_{\xi,v}(0))+sD^{si}_{s\xi}(J'_{\xi,v}(0)))
=\tau_{\gamma_{\xi}\vert_{[0,s]}}(sD^{si}_{s\xi}(v)).
$$
Hence we obtain 
$$\exp^{\perp}_{\ast}(v)=J_{\xi,v}(1)=\tau_{\gamma_{\xi}\vert_{[0,1]}}(D^{si}_{\xi}(v)).\leqno{(2.2)}$$
Take $X\in T_xB$.  
Let $J_{\xi,X}$ be the strongly $B$-Jacobi field along $\gamma_{\xi}$ with $J_{\xi,X}(0)=X$, where 
``strongly'' means that $J'_{X,\xi}(0)\in T_xB$.  Since $F$ is reflective (hence totally geodesic), 
we have $J'_{X,\xi}(0)=0$.  
Hence, according to $(1.2)$ in [Ko1], we obtain 
$$
J_{\xi,X}=\tau_{\gamma_{\xi}\vert_{[0,s]}}
(D^{co}_{s\xi}(J_{\xi,X}(0))+sD^{si}_{s\xi}(J'_{\xi,X}(0)))
=\tau_{\gamma_{\xi}\vert_{[0,s]}}(D^{co}_{s\xi}(X)).
$$
On the other hand, according to $(1.1)$ in [Ko1], we have 
$$\exp^{\perp}_{\ast}(\widetilde X_{\xi})=J_{\xi,X}(1)+\frac{Xr}{r(x)}\gamma'_{\xi}(1).$$
Therefore, we can derive 
$$\exp^{\perp}_{\ast}(\widetilde X_{\xi})=\tau_{\gamma_{\xi}\vert_{[0,1]}}(D^{co}_{\xi}(X))
+\frac{Xr}{r(x)}\gamma'_{\xi}(1).\leqno{(2.3)}$$
Take $v\in{\cal V}_{\xi}$ with 
$\tau_x^{-1}(v)\in\mathfrak p_{k\beta}$ and 
$X\in T_xB$ with $\tau_x^{-1}X\in\mathfrak p_{k\beta}$, where $k\in{\cal K}\cup\{0\}$.  
According to (i) of Theorem A in [Ko1], we have 
$$\begin{array}{l}
\displaystyle{Av=\frac{1}{\sqrt{1+\vert\vert (D^{co}_{\xi})^{-1}({\rm grad}\,r)_x\vert\vert^2}}}\\
\hspace{0.9truecm}\displaystyle{\times\left(\frac{\sqrt{\varepsilon}kb}
{\tan(\sqrt{\varepsilon}kbr(x))}v
+((\tau_{\gamma_{\xi}}\circ\tau_x)(Z_{v,\xi}(1)))_T\right),}
\end{array}
\leqno{(2.4)}$$
where ${\rm grad}\,r$ is the gradient vector field of $r$ (with respect to the induced metric on 
$B$), $(\cdot)_T$ is the $TM$-component of $(\cdot)$ and 
$Z_{v,\xi}(:{\Bbb R}\to\mathfrak p)$ is the solution of the following differential equation:
$$\begin{array}{l}
\displaystyle{Z''(s)={\rm ad}(\tau_x^{-1}\xi)^2(Z(s))
-2[[\tau_x^{-1}\xi,\,\tau_x^{-1}((D^{co}_{s\xi}\circ(D^{co}_{\xi})^{-2})(({\rm grad}\,r)_x)],\,
\tau_x^{-1}(D^{co}_{s\xi}(v_0)]}\\
\hspace{1.55truecm}\displaystyle{-2\left[[\tau_x^{-1}\xi,s\tau_x^{-1}(D^{si}_{s\xi}(v_0))],\,
\tau_x^{-1}\left(\frac{dD^{co}_{s\xi}}{ds}\circ(D^{co}_{\xi})^{-2}\right)(({\rm grad}\,r)_x)\right]}
\end{array}\leqno{(2.5)}$$
satisfying the following initial condition:
$$Z(0)=\tau_x^{-1}(((D^{co})^{-2})_{\ast\xi}v)({\rm grad}\,r)_x),\leqno{(2.6)}$$
and 
$$Z'(0)=0,\leqno{(2.7)}$$
where 
$v_0$ is the element of $T^{\perp}_xB$ corresponding to $v$ under the identification 
of $T_{\xi}(T^{\perp}_xB)$ and $T^{\perp}_xB$, 
$((D^{co})^{-2})_{\ast\xi}$ is the differential of $(D^{co})^{-2}$ at $\xi$ 
(which is regarded as a map from $T^{\perp}_xB$ to $T^{\ast}_xB\otimes T_xB$) and 
$((D^{co})^{-2})_{\ast\xi}v$ is regarded as an element of 
$T^{\ast}_xB\otimes T_xB$ under the natural identification of 
$T_{(D^{co}_{\xi})^{-2}}(T^{\ast}_xB\otimes T_xB)$ and $T^{\ast}_xB\otimes T_xB$.  
According to (ii) of Theorem A in [Ko1], we have 
$$\begin{array}{l}
\displaystyle{A\widetilde X_{\xi}=
-\frac{\sqrt{\varepsilon}kb\tan(\sqrt{\varepsilon}kbr(x))}
{\sqrt{1+\vert\vert (D^{co}_{\xi})^{-1}({\rm grad}\,r)_x\vert\vert^2}}\widetilde X_{\xi}}\\
\hspace{1.7truecm}\displaystyle{
+\frac{(Xr)\sqrt{\varepsilon}kb\tan(\sqrt{\varepsilon}kbr(x))}
{(1+\vert\vert (D^{co}_{\xi})^{-1}({\rm grad}\,r)_x\vert\vert^2)^{3/2}}
\widetilde{((D^{co}_{\xi})^{-2}({\rm grad}\,r)_x)}_{\xi}}\\
\hspace{1.7truecm}\displaystyle{
-\frac{1}{\sqrt{1+\vert\vert (D^{co}_{\xi})^{-1}({\rm grad}\,r)_x\vert\vert^2}}
((\tau_{\gamma_{\xi}}\circ \tau_x)(Z_{X,\xi}(1)))_T,}
\end{array}\leqno{(2.8)}$$
where $Z_{X,\xi}(:{\Bbb R}\to\mathfrak p)$ is the solution of the following differential 
equation:
$$\begin{array}{l}
\hspace{0.5truecm}\displaystyle{Z''(s)}\\
\displaystyle{={\rm ad}(\tau_x^{-1}\xi)^2(Z(s))
+\frac{2Xr}{r(x)}{\rm ad}(\tau_x^{-1}\xi)^2(\tau_x^{-1}
((D^{co}_{s\xi}\circ(D^{co}_{\xi})^{-2})(({\rm grad}\,r)_x)))}\\
\hspace{0.5truecm}\displaystyle{-2\left[[\tau_x^{-1}\xi,\tau_x^{-1}
((D^{co}_{s\xi}\circ(D^{co}_{\xi})^{-2})(({\rm grad}\,r)_x))],\,
\tau_x^{-1}\frac{dD^{co}_{s\xi}}{ds}X\right]}\\
\hspace{0.5truecm}\displaystyle{-2\left[[\tau_x^{-1}\xi,\tau_x^{-1}
(D^{co}_{s\xi}X)],\,
\tau_x^{-1}\left(\frac{dD^{co}_{s\xi}}{ds}\circ(D^{co}_{\xi})^{-2}\right)
(({\rm grad}\,r)_x)\right]}
\end{array}\leqno{(2.9)}$$
satisfying the following initial condition:
$$Z(0)=\tau_x^{-1}(\nabla^F)^{\pi\vert_M}_{\widetilde X_{\xi}}
((D^{co}_{\cdot})^{-2}({\rm grad}\,r)_x)
\leqno{(2.10)}$$
and 
$$Z'(0)=-[[\tau_x^{-1}\xi,\tau_x^{-1}X],\,\tau_x^{-1}(D^{co}_{\xi})^{-2}
({\rm grad}\,r)_x],
\leqno{(2.11)}$$
where $(\nabla^F)^{\pi\vert_M}$ is the pull-back connection of the Riemannian connection $\nabla^F$ of $F$ 
by $\pi\vert_M$.  
Here we give the table of the correspondence between the above notations and the notations 
in [Ko1].  

\vspace{0.2truecm}

$$\begin{tabular}{|c|c|c|}
\hline
The notation in [Ko1] & The notations in this paper & Remark\\
\hline
$M$ & $B$ & \\
\hline
$\varepsilon$ & $r$ & \\
\hline
$t_{\varepsilon}(M)$ & $M(=t_r(B))$ & \\
\hline
$A^{\varepsilon}_E$ & $A$ & \\
\hline
$B_{\xi}$ & $D^{co}_{\xi}$ & by the reflectivity of $F$\\
\hline
$g_{\ast}$ & $\tau_x$ & \\
\hline
$\mu(g_{\ast}^{-1}\xi)$ & $\sqrt{\varepsilon}kbr(x)$ & \\
\hline
\end{tabular}$$

\vspace{0.25truecm}

\centerline{{\bf Table 1.}}

\vspace{0.5truecm}

\noindent
{\bf Assumption.} Assume that there exists $k_0\in\{0,1,2\}$ such that 
$\tau_x^{-1}(({\rm grad}\,r)_x)\in\mathfrak p_{k_0\beta}\cap\mathfrak p'$ holds for any 
$\xi\in M\cap T^{\perp}_xB$, 
where we note that $\mathfrak p_{k_0\beta}$ depends on the choice of $\xi$.  

\vspace{0.4truecm}

\noindent
Denote by $(\cdot)_k$ the $\tau_x(\mathfrak p_{k\beta})$-component of 
$(\cdot)\in T_xB$, where $k\in{\cal K}\cup\{0\}$.  Then we have 
$$
(D^{co}_{\xi})^{-j}(({\rm grad}\,r)_x)=\frac{1}{\cos^j(\sqrt{\varepsilon}k_0br(x))}({\rm grad}\,r)_x\,\,\,
(j=1,2)\leqno{(2.12)}$$
and 
$$((D^{co})^{-2}_{\ast\xi}(v))(({\rm grad}\,r)_x)=\frac{2\sqrt{\varepsilon}k_0br(x)
\sin(\sqrt{\varepsilon}k_0br(x))}{\cos^3(\sqrt{\varepsilon}k_0br(x))}({\rm grad}\,r)_x.\leqno{(2.13)}$$
Since $F$ is reflective, it is an orbit of a Hermann action (i.e., the action of 
a symmetric subgroup of $G$) and hence it is homogeneous.  
Let $c:I\to B$ be the (homogeneous) geodesic in $B$ with $c'(0)=X$ and 
$\widehat{\xi}$ the normal vector field of $B$ along $c$ 
such that $\widehat{\xi}(0)=\xi$, 
$\frac{\widehat{\xi}}{\vert\vert\widehat{\xi}\vert\vert}$ is 
parallel (with respect to the normal connection) and that 
$\vert\vert\widehat{\xi}(t)\vert\vert=r(\pi(c(t)))$ for all $t$ 
in the domain of $\widehat{\xi}$.  
This curve $\widehat{\xi}$ is regarded as a curve in $M$ with 
$\widetilde{\xi}'(0)=\widetilde X_{\xi}$.  
Since $c$ is a homogeneous curve, it is described as 
$c(t)=\widehat a(t)x\,\,(t\in I)$ for some curve $\widehat a:I\to G$.  
Then, since $\frac{\widehat{\xi}}{\vert\vert\widehat{\xi}\vert\vert}$ is 
parallel and $F$ is a submanifold with section (i.e., with Lie triple systematic 
normal bundle in the sense of [Ko1]) by the reflectivity of $F$, we can show 
$${\rm Span}\{\widehat{\xi}(t)\}=\widehat a(t)_{\ast}({\rm Span}\{\xi\})
(=(\widehat a(t)_{\ast}\circ\tau_x)(\mathfrak b)\,\,\,\,(t\in I)\leqno{(2.14)}$$
(see Theorem 5.5.12 of [PT]) and 
$$T_{c(t)}B=\sum_{k\in{\cal K}\cup\{0\}}(T_{c(t)}B\cap(\widehat a(t)_{\ast}\circ\tau_x)
(\mathfrak p_{k\beta})).\leqno{(2.15)}$$
Denote by $(\cdot)_k$ the 
$(\widehat a(t)_{\ast}\circ\tau_x)(\mathfrak p_{k\beta})$-component of 
$(\cdot)\in T_{c(t)}B$.  
Then we can show 
$$\begin{array}{l}
\displaystyle{(\nabla^F)^{\pi\vert_M}_{\widetilde X_{\xi}}((D^{co}_{\cdot})^{-2}{\rm grad}\,r)
=\left.\frac{\nabla^F}{dt}\right\vert_{t=0}(D^{co}_{\widehat{\xi}(t)})^{-2}
(({\rm grad}\,r)_{c(t)})}\\
\displaystyle{=\left.\frac{\nabla^F}{dt}\right\vert_{t=0}
\left(\frac{1}{\cos^2(\sqrt{\varepsilon}k_0br(c(t)))}({\rm grad}\,r)_{c(t)}\right)}\\
\displaystyle{=\frac{2(Xr)\sqrt{\varepsilon}k_0b\sin(\sqrt{\varepsilon}k_0br(x))}
{\cos^3(\sqrt{\varepsilon}k_0br(x))}({\rm grad}\,r)_x}\\
\hspace{1.4truecm}
\displaystyle{+\frac{1}{\cos^2(\sqrt{\varepsilon}k_0br(x))}
\left.\frac{\nabla^F}{dt}\right\vert_{t=0}({\rm grad}\,r)_{c(t)},}
\end{array}\leqno{(2.16)}$$
where $\frac{\nabla^F}{dt}$ is the covariant derivative along $c$ 
with respect to $\nabla^F$.  
Since $F$ is a submanifold with section, we have 
$$[[\mathfrak p',{\mathfrak p'}^{\perp}],{\mathfrak p'}^{\perp}]\subset\mathfrak p',\,\,
[[{\mathfrak p'}^{\perp},{\mathfrak p'}^{\perp}],\mathfrak p']\subset\mathfrak p'\,\,
{\rm and}\,\,[[{\mathfrak p'}^{\perp},\mathfrak p'],\mathfrak p']\subset
{\mathfrak p'}^{\perp}.\leqno{(2.17)}$$
Hence, since $Z_v$ satisfies $(2.5)$, we have 
$$Z''_v(s)\equiv{\rm ad}(\tau_x^{-1}\xi)^2Z_v(s)\,\,\,\,
({\rm mod}\,\mathfrak p')$$
Also, since $Z_v$ satisfies $(2.6)$, it follows from $(2.13)$ that 
$Z_v(0)\equiv 0\,\,\,\,({\rm mod}\,\mathfrak p')$.  
Furthermore, since $Z_v$ satisfies $(2.7)$, we have 
$Z'_v(0)=0$.  
Hence we can show $Z_v(s)\equiv 0\,\,({\rm mod}\,\mathfrak p')$.  
Therefore, form $(2.4)$, we obtain 
$$Av\equiv\frac{\sqrt{\varepsilon}kb\cos(\sqrt{\varepsilon}k_0br(x))}
{\tan(\sqrt{\varepsilon}kbr(x))
\sqrt{\cos^2(\sqrt{\varepsilon}k_0br(x))+\vert\vert({\rm grad}\,r)_x\vert\vert^2}}v
\quad({\rm mod}\,\mathfrak p').\leqno{(2.18)}$$
On the other hand, from $(2.3)$ and $(2.8)$, we have 
$$\begin{array}{l}
\displaystyle{A\widetilde X_{\xi}\equiv\frac{g(A\widetilde X_{\xi},\widetilde X_{\xi})}
{\vert\vert\widetilde X\vert\vert^2}\widetilde X_{\xi}}\\
\displaystyle{=\left\{
-\frac{\sqrt{\varepsilon}kb\tan(\sqrt{\varepsilon}kbr(x))\cos(\sqrt{\varepsilon}k_0br(x))}
{\sqrt{\cos^2(\sqrt{\varepsilon}k_0br(x))+\vert\vert({\rm grad}\,r)_x\vert\vert^2}}\right.}\\
\hspace{0.5truecm}\displaystyle{
+\frac{(Xr)^2\sqrt{\varepsilon}kb\tan(\sqrt{\varepsilon}kbr(x))\cos(\sqrt{\varepsilon}k_0br(x))}
{\sqrt{\cos^2(\sqrt{\varepsilon}k_0br(x))+\vert\vert({\rm grad}\,r)_x\vert\vert^2}}}\\
\hspace{1truecm}\displaystyle{\times
\frac{(\cos(\sqrt{\varepsilon}k_0br(x))\cos(\sqrt{\varepsilon}kbr(x))+\vert\vert({\rm grad}\,r)_x\vert\vert^2)}
{(\cos^2(\sqrt{\varepsilon}kbr(x))\vert\vert X\vert\vert^2+(Xr)^2)}}\\
\hspace{0.5truecm}\displaystyle{\left.
-\frac{\cos(\sqrt{\varepsilon}k_0br(x))\bar g(((\tau_{\gamma_{\xi}}\circ\tau_x)(Z_{X,\xi}(1)))_T,
\exp^{\perp}_{\ast}(\widetilde X_{\xi}))}
{\sqrt{\cos^2(\sqrt{\varepsilon}k_0br(x))+\vert\vert({\rm grad}\,r)_x\vert\vert^2}}\right.}\\
\hspace{1truecm}\displaystyle{\left.\times
\frac{1}{(\cos^2(\sqrt{\varepsilon}kbr(x))\vert\vert X\vert\vert^2+(Xr)^2)}\right\}\widetilde X_{\xi}}\\
$\,$\\
\hspace{5truecm}\displaystyle{({\rm mod}\,T_{\xi}M\ominus T^{\perp}{\rm span}\{\widetilde X_{\xi}\}).}
\end{array}\leqno{(2.19)}$$
Denote by $(\bullet)_{\mathfrak p'}$ (resp. $(\bullet)_{{\mathfrak p'}^{\perp}}$) 
the ${\mathfrak p'}$-component (resp. the ${\mathfrak p'}^{\perp}$-component) of $(\bullet)$.  
Also, denote by $(\bullet)_{\mathfrak p'_{\bar k}}$ the $({\mathfrak p'}\cap\mathfrak p_{\bar k})$-component of 
$(\bullet)$ ($\bar k\in{\cal K}\cup\{0\}$), and 
$(\bullet)_{{\mathfrak p'}^{\perp}_{\bar k}}$ the $({\mathfrak p'}^{\perp}\cap\mathfrak p_{\bar k})$-component of 
$(\bullet)$ ($\bar k\in\{0,1,2\}$).  
Now we shall calculate $Z_{v,\xi}$.  
Since $Z_{v,\xi}$ satisfies $(2.5)$, it follows from the relations in $(2.17)$ that 
$$(Z_{v,\xi})''_{\mathfrak p'_{\bar k}}(s)=\left\{
\begin{array}{cl}
\displaystyle{
\begin{array}{l}
\displaystyle{-(\sqrt{\varepsilon}k_0br(x))^2(Z_{v,\xi})_{\mathfrak p'_{\bar k}}(s)}\\
\displaystyle{+m\cos(s\sqrt{\varepsilon}kbr(x))\cos(s\sqrt{\varepsilon}k_0br(x)){\bf C}^1_{v,\xi}}\\
\displaystyle{\pm\sin(s\sqrt{\varepsilon}kbr(x))\sin(s\sqrt{\varepsilon}k_0br(x)){\bf C}^2_{v,\xi}}\\
\end{array}} & (\bar k=\vert k_0\pm k\vert)\\
$\,$ & \\
\displaystyle{-(\sqrt{\varepsilon}\bar kbr(x))^2(Z_{v,\xi})_{\mathfrak p'_{\bar k}}(s)} & 
(\bar k\not=\vert k_0\pm k\vert)
\end{array}\right.\leqno{(2.20)}$$
and 
$$(Z_{v,\xi})''_{{\mathfrak p'}^{\perp}_{\bar k}}(s)
=-(\sqrt{\varepsilon}\bar kbr(x))^2(Z_{v,\xi})_{{\mathfrak p'}^{\perp}_{\bar k}}(s),
\leqno{(2.21)}$$
where ${\bf C}^i_{v,\xi}$ ($i=1,2$) are given by 
$$
{\bf C}^1_{v,\xi}:=-\frac{2}{\cos^2(\sqrt{\varepsilon}k_0br(x))}
[[\tau_x^{-1}\xi,\tau_x^{-1}({\rm grad}\,r)_x)],\tau_x^{-1}v]
$$
and 
$${\bf C}^2_{v,\xi}:=\frac{k_0}{k\cos^2(\sqrt{\varepsilon}kbr(x))}
[[\tau_x^{-1}\xi,\tau_x^{-1}v],\tau_x^{-1}({\rm grad}\,r)_x)].
$$
Note that ${\bf C}^1_{v,\xi}$ and ${\bf C}^2_{v,\xi}$ belong to 
${\mathfrak p'}_{k_0+k}+{\mathfrak p'}_{\vert k_0-k\vert}$ because of 
$[[\mathfrak p_0,\mathfrak p_{k_0}],\mathfrak p_k],\,[[\mathfrak p_0,\mathfrak p_k],\mathfrak p_{k_0}]$\newline
$\in\mathfrak p_{k_0+k}+\mathfrak p_{\vert k_0-k\vert}$ and $(2.17)$.  
Since $Z_{v,\xi}$ satisfies $(2.6)$, it follows from $(2.13)$ that 
$$(Z_{v,\xi})_{\mathfrak p'_{\bar k}}(0)=
\left\{
\begin{array}{cl}
\displaystyle{\frac{2\sqrt{\varepsilon}k_0br(x)\sin(\sqrt{\varepsilon}k_0br(x))}
{\cos^3(\sqrt{\varepsilon}k_0br(x))}(\tau_x^{-1}({\rm grad}\,r)_x))_{\mathfrak p'_{\bar k}}} & (\bar k=k_0)\\
0 & (\bar k\not=k_0)
\end{array}\right.\leqno{(2.22)}$$
and 
$$(Z_{v,\xi})_{{\mathfrak p'}_{\bar k}^{\perp}}(0)=0.\leqno{(2.23)}$$
Also, since $Z_{v,\xi}$ satisfies $(2.7)$, we have 
$$(Z_{v,\xi})'_{\mathfrak p'_{\bar k}}(0)=0\leqno{(2.24)}$$
and 
$$(Z_{v,\xi})'_{{\mathfrak p'}_{\bar k}^{\perp}}(0)=0.\leqno{(2.25)}$$
By solving $(2.20)$ under the initial conditions $(2.22)$ and $(2.24)$, we can derive 
$$\begin{array}{l}
\displaystyle{(Z_{v,\xi})_{\mathfrak p'_{\bar k}}(s)}\\
\displaystyle{=\left\{
\begin{array}{cl}
\displaystyle{
\begin{array}{l}
\displaystyle{\pm\frac{1}{4}\left(\frac{\cos(s\sqrt{\varepsilon}k_0br(x))\cos(s\sqrt{\varepsilon}kbr(x))}
{k_0k(\sqrt{\varepsilon}br(x))^2}\right.}\\
\hspace{0.5truecm}\displaystyle{+\frac{s\sin(s\sqrt{\varepsilon}(k_0\pm k)br(x))}
{\sqrt{\varepsilon}(k_0\pm k)br(x)}}\\
\hspace{0.5truecm}\displaystyle{\left.-\frac{\cos(s\sqrt{\varepsilon}(k_0\pm k)br(x))}
{k_0k(\vert\sqrt{\varepsilon}br(x))^2}\right)
({\bf C}^1_{v,\xi})_{\mathfrak p'_{\vert k_0\pm k\vert}}}\\
\hspace{0.5truecm}\displaystyle{\pm\frac{1}{4}
\left(\frac{\cos(s\sqrt{\varepsilon}k_0br(x))\cos(s\sqrt{\varepsilon}kbr(x))}
{k_0k(\sqrt{\varepsilon}br(x))^2}\right.}\\
\hspace{0.5truecm}\displaystyle{-\frac{s\sin(s\sqrt{\varepsilon}(k_0\pm k)br(x))}
{\sqrt{\varepsilon}(k_0\pm k)br(x)}}\\
\hspace{0.5truecm}\displaystyle{\left.-\frac{\cos(s\sqrt{\varepsilon}(k_0\pm k)br(x))}
{k_0k(\vert\sqrt{\varepsilon}br(x))^2}\right)
({\bf C}^2_{v,\xi})_{\mathfrak p'_{\vert k_0\pm k\vert}}}
\end{array}} & (\bar k=\vert k_0\pm k\vert\not=0)\\
$\,$ & \\
\displaystyle{
\begin{array}{l}
\displaystyle{\pm\frac{\sin(s\sqrt{\varepsilon}k_0br(x))\sin(s\sqrt{\varepsilon}kbr(x))}
{\sqrt{\varepsilon}(k_0\pm k)br(x)}({\bf C}^1_{v,\xi})_{\mathfrak p'_0}}\\
\hspace{0truecm}\displaystyle{-\frac{s\sin(s\sqrt{\varepsilon}(k_0\pm k)br(x))}
{\sqrt{\varepsilon}(k_0\pm k)br(x)}({\bf C}^2_{v,\xi})_{\mathfrak p'_0}}
\end{array}} & (\bar k=\vert k_0-k\vert=0)\\
$\,$ & \\
\displaystyle{
\begin{array}{l}
\displaystyle{\frac{2\sqrt{\varepsilon}k_0br(x)\sin(\sqrt{\varepsilon}k_0br(x))\cos(s\sqrt{\varepsilon}kbr(x))}
{\cos^3(\sqrt{\varepsilon}(k_0\pm k)br(x))}}\\
\displaystyle{\times\tau_x^{-1}(({\rm grad}\,r)_x)}
\end{array}} & (\bar k=k_0)\\
$\,$ & \\
0 & (\bar k\not=k_0,\vert k_0\pm k\vert).
\end{array}\right.}
\end{array}
\leqno{(2.26)}$$
Also, by solving $(2.21)$ under the initial conditions $(2.23)$ and $(2.25)$, we can derive 
$$(Z_{v,\xi})_{{\mathfrak p'}_{\bar k}^{\perp}}\equiv0.\leqno{(2.27)}$$
According to Lemma 3.2 in [Ko1], we have 
$$
N_{\xi}=\frac{\cos(\sqrt{\varepsilon}k_0br(x))\vert\vert\xi\vert\vert\tau_{\gamma_{\xi}}(\xi)
-\tau_{\gamma_{\xi}}(({\rm grad}\,r)_x)}
{\sqrt{\cos^2(\sqrt{\varepsilon}k_0br(x))+\vert\vert({\rm grad}\,r)_x\vert\vert^2}}
\leqno{(2.28)}$$
Also, according to $(2.26)$,  $(Z_{v,\xi})_{\mathfrak p'}(s)$ is descirbed as 
$$(Z_{v,\xi})_{\mathfrak p'}(s)=\sum_{i=1}^2
\left(a^{i,+}_{v,\xi}(s)({\bf C}^i_{v,\xi})_{\mathfrak p'_{k_0+k}}
+a^{i,-}_{v,\xi}(s)({\bf C}^i_{v,\xi})_{\mathfrak p'_{\vert k_0-k\vert}}\right)
+\bar a_{v,\xi}(s)\tau_x^{-1}(({\rm grad}\,r)_x)$$
in terms of some explicitly described functions $a^{i,\pm}_{v,\xi}$ ($i=1,2$).  
From $(2.26),(2.27),(2.28)$ and the definitions of ${\bf C}^i_{v,\xi}$ ($i=1,2$), we can derive 
$$\begin{array}{l}
\hspace{0.5truecm}\displaystyle{(\tau_{\gamma_{\xi}}\circ\tau_x)((Z_{v,\xi})(1)))_T}\\
\displaystyle{=(\tau_{\gamma_{\xi}}\circ\tau_x)((Z_{v,\xi})_{\mathfrak p'}(1)))
-\bar g((\tau_{\gamma_{\xi}}\circ\tau_x)((Z_{v,\xi})_{\mathfrak p'}(1))),\,N_{\xi})N_{\xi}}\\
\displaystyle{=\sum_{i=1}^2
(\tau_{\gamma_{\xi}}\circ\tau_x)\left(a^{i,+}_{v,\xi}(1)({\bf C}^i_{v,\xi})_{\mathfrak p'_{k_0+k}}
+a^{i,-}_{v,\xi}(1)({\bf C}^i_{v,\xi})_{\mathfrak p'_{\vert k_0-k\vert}}\right)}\\
\hspace{0.5truecm}\displaystyle{+\frac{\bar a_{v,\xi}(1)\cos(\sqrt{\varepsilon}k_0br(x))}
{\cos^2(\sqrt{\varepsilon}k_0br(x))+\vert\vert({\rm grad}\,r)_x\vert\vert^2}}\\
\hspace{1truecm}\displaystyle{\times
\tau_{\gamma_{\xi}}\left(\cos(\sqrt{\varepsilon}k_0br(x))({\rm grad}\,r)_x
+\frac{\vert\vert({\rm grad}\,r)_x\vert\vert^2}{r(x)}\xi\right),}
\end{array}\leqno{(2.29)}$$
where we use 
$$\begin{array}{l}
\hspace{0.5truecm}\displaystyle{\bar g([[\tau_x^{-1}\xi,\tau_x^{-1}(({\rm grad}\,r)_x)],\tau_x^{-1}(v)],
\tau_x^{-1}\xi)}\\
\displaystyle{=\bar g({\rm ad}(\tau_x^{-1}\xi)^2(\tau_x^{-1}(({\rm grad}\,r)_x)),
\tau_x^{-1}(v))=0}
\end{array}\leqno{(2.30)}$$
and 
$$\begin{array}{l}
\hspace{0.5truecm}\displaystyle{\bar g([[\tau_x^{-1}\xi,\tau_x^{-1}(v)],\tau_x^{-1}(({\rm grad}\,r)_x)],
\tau_x^{-1}\xi)}\\
\displaystyle{=\bar g({\rm ad}(\tau_x^{-1}\xi)^2(\tau_x^{-1}(({\rm grad}\,r)_x)),\tau_x^{-1}(v))=0.}
\end{array}\leqno{(2.31)}$$
From $(2.3),(2.4)$ and $(2.29)$, we can derive the following relations:
$$\begin{array}{l}
\displaystyle{Av=\frac{\cos(\sqrt{\varepsilon}k_0br(x))}
{\sqrt{\cos^2(\sqrt{\varepsilon}k_0br(x))+\vert\vert({\rm grad}\,r)_x\vert\vert^2}}
\left\{\frac{\sqrt{\varepsilon}kb}{\tan(\sqrt{\varepsilon}kbr(x))}v\right.}\\
\hspace{1truecm}\displaystyle{+\sum_{i=1}^2
\left(\frac{a^{i,+}_{v,\xi}(1)}{\cos(\sqrt{\varepsilon}(k_0+k)br(x))}
\widetilde{(\tau_x(({\bf C}^i_{v,\xi})_{\mathfrak p'_{k_0+k}}))}_{\xi}\right.}\\
\hspace{2.21truecm}\displaystyle{\left.
+\frac{a^{i,-}_{v,\xi}(1)}{\cos(\sqrt{\varepsilon}(k_0-k)br(x))}
\widetilde{(\tau_x(({\bf C}^i_{v,\xi})_{\mathfrak p'_{\vert k_0-k\vert}}))}_{\xi}\right)}\\
\hspace{2.21truecm}\displaystyle{\left.
+\frac{\bar a_{v,\xi}(1)\cos(\sqrt{\varepsilon}k_0br(x))}
{\cos^2(\sqrt{\varepsilon}k_0br(x))+\vert\vert({\rm grad}\,r)_x\vert\vert^2}
\widetilde{(({\rm grad}\,r)_x)}_{\xi}\right\}.}
\end{array}
\leqno{(2.32)}$$

Next we shall calculate $Z_{X,\xi}X$.  
Since $Z_{X,\xi}$ satisfies $(2.9)$, it follows from the relations in $(2.17)$ that 
$$(Z_{X,\xi})''_{\mathfrak p'_{\bar k}}(s)=\left\{
\begin{array}{ll}
\displaystyle{
\begin{array}{l}
\displaystyle{-(\sqrt{\varepsilon}k_0br(x))^2(Z_{X,\xi})_{\mathfrak p'_{k_0}}(s)}\\
\displaystyle{+\cos(s\sqrt{\varepsilon}k_0br(x)){\bf C}^1_{X,\xi}}
\end{array}} & (\bar k=k_0)\\
$\,$ & \\
\displaystyle{-(\sqrt{\varepsilon}\bar kbr(x))^2(Z_{X,\xi})_{\mathfrak p'_{\bar k}}(s)} & (\bar k\not=k_0)
\end{array}\right.\leqno{(2.33)}$$
and 
$$\begin{array}{l}
\displaystyle{(Z_{X,\xi})''_{{\mathfrak p'}^{\perp}_{\bar k}}(s)}\\
\displaystyle{
=\left\{
\begin{array}{ll}
\displaystyle{
\begin{array}{l}
\displaystyle{-(\sqrt{\varepsilon}\bar kbr(x))^2(Z_{X,\xi})_{{\mathfrak p'}^{\perp}_{\bar k}}(s)}\\
\displaystyle{\pm\sin(s\sqrt{\varepsilon}kbr(x))\cos(s\sqrt{\varepsilon}k_0br(x))
({\bf C}^2_{X,\xi})_{{\mathfrak p'}^{\perp}_{\vert k_0\pm k\vert}}}\\
\displaystyle{+\sin(s\sqrt{\varepsilon}k_0br(x))\cos(s\sqrt{\varepsilon}kbr(x))
({\bf C}^3_{X,\xi})_{{\mathfrak p'}^{\perp}_{\vert k_0\pm k\vert}}}
\end{array}} & (\bar k=\vert k_0\pm k\vert)\\
$\,$ & \\
\displaystyle{-(\sqrt{\varepsilon}\bar kbr(x))^2(Z_{X,\xi})_{{\mathfrak p'}^{\perp}_{\bar k}}(s)} 
& (\bar k\not=\vert k_0\pm k\vert),
\end{array}\right.}
\end{array}
\leqno{(2.34)}$$
where ${\bf C}^i_{X,\xi}$ ($i=1,2,3$) are given by 
$$\begin{array}{l}
\displaystyle{{\bf C}^1_{X,\xi}:=-\frac{2(Xr)(\sqrt{\varepsilon}k_0b)^2r(x)}
{\cos^2(\sqrt{\varepsilon}k_0br(x))}\tau_x^{-1}(({\rm grad}\,r)_x)\,\,\,\,(\in\mathfrak p'_{k_0}),}\\
\displaystyle{{\bf C}^2_{X,\xi}:=\frac{2\sqrt{\varepsilon}kbr(x)}
{\cos^2(\sqrt{\varepsilon}k_0br(x))}[[\tau_x^{-1}\xi,\tau_x^{-1}({\rm grad}\,r)_x)],\tau_x^{-1}X]}
\end{array}$$
and 
$${\bf C}^3_{X,\xi}:=\frac{2\sqrt{\varepsilon}k_0br(x)}
{\cos^2(\sqrt{\varepsilon}k_0br(x))}[[\tau_x^{-1}\xi,\tau_x^{-1}X],\tau_x^{-1}({\rm grad}\,r)_x)].$$
Note that ${\bf C}^2_{X,\xi}$ and ${\bf C}^3_{X,\xi}$ belong to 
${\mathfrak p'}^{\perp}_{k_0+k}+{\mathfrak p'}^{\perp}_{\vert k_0-k\vert}$ because of 
$[[\mathfrak p_0,\mathfrak p_{k_0}],\mathfrak p_k],\,[[\mathfrak p_0,\mathfrak p_{k_0}],\mathfrak p_k]
\in\mathfrak p_{k_0+k}+\mathfrak p_{\vert k_0-k\vert}$ and $(2.17)$.  
Also, since $Z_{X,\xi}$ satisfies $(2.10)$, it follows from $(2.16)$ that 
$$\begin{array}{l}
\displaystyle{(Z_{X,\xi})_{\mathfrak p'_{\bar k}}(0)=(\tau_x^{-1}(\nabla^F)^{\pi\vert_M}_{\widetilde X_{\xi}}
((D^{co}_{\cdot})^{-2}{\rm grad}\,r))_{\mathfrak p'_{\bar k}}}\\
\displaystyle{=\left\{
\begin{array}{ll}
\displaystyle{
\begin{array}{l}
\displaystyle{\frac{2\sqrt{\varepsilon}k_0b(Xr)\sin(\sqrt{\varepsilon}k_0br(x))}
{\cos^3(\sqrt{\varepsilon}k_0br(x))}\tau_x^{-1}(({\rm grad}\,r)_x)}\\
\displaystyle{+\frac{1}{\cos^2(\sqrt{\varepsilon}k_0br(x))}
(\tau_x^{-1}(\nabla^F_X{\rm grad}\,r))_{\mathfrak p'_{k_0}}}
\end{array}} & (\bar k=k_0)\\
$\,$ & \\
\displaystyle{\frac{1}{\cos^2(\sqrt{\varepsilon}k_0br(x))}
(\tau_x^{-1}(\nabla^F_X{\rm grad}\,r))_{\mathfrak p'_{\bar k}}} & (\bar k\not=k_0)
\end{array}\right.}
\end{array}
\leqno{(2.35)}$$
and 
$$(Z_{X,\xi})_{{\mathfrak p'_{\bar k}}^{\perp}}(0)=0.\leqno{(2.36)}$$
Furthermore, since $Z_{X,\xi}'(0)$ satisfies $(2.10)$, it follows from the third relation in $(2.16)$ that 
$$(Z_{X,\xi})'_{\mathfrak p'_{\bar k}}(0)=0\leqno{(2.37)}$$
and 
$$(Z_{X,\xi})'_{{\mathfrak p'}_{\bar k}^{\perp}}(0)
=\left\{
\begin{array}{cl}
\displaystyle{-\frac{1}{2\sqrt{\varepsilon}k_0br(x)}
({\bf C}^3_{X,\xi})_{{\mathfrak p'}^{\perp}_{\vert k_0\pm k\vert}}} & (\bar k=\vert k_0\pm k\vert)\\
\displaystyle{0} & (\bar k\not=\vert k_0\pm k\vert).
\end{array}\right.\leqno{(2.38)}$$
By solving $(2.33)$ under the initial conditions $(2.35)$ and $(2.37)$, we can derive 
$$\begin{array}{l}
\displaystyle{(Z_{X,\xi})_{\mathfrak p'_{\bar k}}(s)}\\
\displaystyle{
=\left\{
\begin{array}{ll}
\displaystyle{
\begin{array}{l}
\displaystyle{\frac{2\sqrt{\varepsilon}k_0b(Xr)\sin(\sqrt{\varepsilon}k_0br(x))}
{\cos^3(\sqrt{\varepsilon}k_0br(x))}}\\
\displaystyle{\times\cos(s\sqrt{\varepsilon}k_0br(x))\tau_x^{-1}(({\rm grad}\,r)_x)}\\
\displaystyle{+\frac{\cos(s\sqrt{\varepsilon}k_0br(x))}{\cos^2(\sqrt{\varepsilon}k_0br(x))}
(\tau_x^{-1}(\nabla^F_X{\rm grad}\,r))_{\mathfrak p'_{k_0}}}\\
\displaystyle{+\frac{s\sin(s\sqrt{\varepsilon}k_0br(x))}{2\sqrt{\varepsilon}k_0br(x)}}{\bf C}^1_{X,\xi}
\end{array}} & (\bar k=k_0)\\
$\,$ & \\
\displaystyle{\frac{\cos(s\sqrt{\varepsilon}\bar kbr(x))}{\cos^2(\sqrt{\varepsilon}k_0br(x))}
(\tau_x^{-1}(\nabla^F_X{\rm grad}\,r))_{\mathfrak p'_{\bar k}}} & (\bar k\not=k_0)
\end{array}\right.}
\end{array}\leqno{(2.39)}$$
Also, by solving $(2.34)$ under the initial conditions $(2.36)$ and $(2.38)$, we can derive 
$$\begin{array}{l}
\displaystyle{(Z_{X,\xi})_{{\mathfrak p'}^{\perp}_{\bar k}}(s)}\\
\displaystyle{=\left\{
\begin{array}{cl}
\begin{array}{l}
\displaystyle{-\frac{1}{4}\left(\frac{\sin(s\sqrt{\varepsilon}k_0br(x))\cos(s\sqrt{\varepsilon}kbr(x))}
{\pm k_0k(\sqrt{\varepsilon}br(x))^2}\right.}\\
\hspace{0.5truecm}\displaystyle{\left.+\frac{s\cos(s\sqrt{\varepsilon}(k_0\pm k)br(x))}
{\vert k_0\pm k\vert\sqrt{\varepsilon}br(x)}\right)
({\bf C}^2_{X,\xi})_{\mathfrak p'^{\perp}_{\vert k_0\pm k\vert}}}\\
\displaystyle{-\frac{1}{4}\left(\frac{\sin(s\sqrt{\varepsilon}kbr(x))\cos(s\sqrt{\varepsilon}k_0br(x))}
{k_0k(\sqrt{\varepsilon}br(x))^2}\right.}\\
\hspace{0.5truecm}\displaystyle{\left.
+\frac{s\cos(s\sqrt{\varepsilon}(k_0+k)br(x))}{\vert k_0\pm k\vert\sqrt{\varepsilon}br(x)}\right)
({\bf C}^3_{X,\xi})_{\mathfrak p'^{\perp}_{\vert k_0\pm k\vert}}}\\
\displaystyle{-\frac{\sin(s\sqrt{\varepsilon}\vert k_0\pm k\vert br(x))}{2(\sqrt{\varepsilon}br(x))^2k_0
\vert k_0\pm k\vert}
({\bf C}^3_{X,\xi})_{\mathfrak p'^{\perp}_{\vert k_0\pm k\vert}}}
\end{array} & (\bar k=\vert k_0\pm k\vert\not=0)\\
$\,$ & \\
\displaystyle{\begin{array}{l}
\displaystyle{-\frac{1}{2}\left(\frac{\sin(s\sqrt{\varepsilon}k_0br(x))\cos(s\sqrt{\varepsilon}k_0br(x))}
{(\sqrt{\varepsilon}k_0br(x))^2}\right.}\\
\hspace{0.5truecm}\displaystyle{\left.+\frac{s}{\sqrt{\varepsilon}k_0br(x)}\right)
({\bf C}^2_{X,\xi})_{\mathfrak p'^{\perp}_0}}
\end{array}} & (\bar k=\vert k_0-k\vert=0)\\
$\,$ & \\
0 & (\bar k\not=\vert k_0\pm k\vert).
\end{array}\right.}
\end{array}
\leqno{(2.40)}$$
On the other hand, according to Lemma 3.3 in [Ko1], we have 
$$\begin{array}{l}
\hspace{0.5truecm}\displaystyle{
((\tau_{\gamma_{\xi}}\circ\tau_x)((Z_{X,\xi})_{\mathfrak p'}(1)))_T}\\
\displaystyle{=\widetilde{((D^{co}_{\xi})^{-1}(\tau_x(Z_{X,\xi})_{\mathfrak p'}(1)))}_{\xi}}\\
\hspace{0.5truecm}\displaystyle{
-\frac{((D^{co}_{\xi})^{-1}(\tau_x(Z_{X,\xi})_{\mathfrak p'}(1)))r}
{1+\vert\vert(D^{co}_{\xi})^{-1}({\rm grad}\,r)_x\vert\vert^2}
\widetilde{((D^{co}_{\xi})^{-2}({\rm grad}\,r)_x)}_{\xi}.}
\end{array}$$
From this relation, $(2.3)$ and $(2.39)$, we obtain 
$$\begin{array}{l}
\hspace{0.5truecm}\displaystyle{
((\tau_{\gamma_{\xi}}\circ\tau_x)((Z_{X,\xi})_{\mathfrak p'}(1)))_T}\\
\displaystyle{=\frac{2\sqrt{\varepsilon}k_0b\tan(\sqrt{\varepsilon}k_0br(x))\cdot(Xr)}
{\cos^2(\sqrt{\varepsilon}k_0br(x))+\vert\vert({\rm grad}\,r)_x\vert\vert^2}
(\widetilde{({\rm grad}\,r)_x})_{\xi}}\\
\hspace{0.5truecm}\displaystyle{-\frac{g_F(\nabla^F_X{\rm grad}\,r,({\rm grad}\,r)_x)}
{(\cos^2(\sqrt{\varepsilon}k_0br(x))+\vert\vert({\rm grad}\,r)_x\vert\vert^2)
\cos^2(\sqrt{\varepsilon}k_0br(x))}(\widetilde{({\rm grad}\,r)_x})_{\xi}}\\
\hspace{0.5truecm}\displaystyle{+\frac{1}{\cos^2(\sqrt{\varepsilon}k_0br(x))}
(\widetilde{\nabla^F_X{\rm grad}\,r})_{\xi}.}
\end{array}\leqno{(2.41)}$$
According to Lemma 3.2 in [Ko1], we have 
$$
N_{\xi}=\frac{\cos(\sqrt{\varepsilon}k_0br(x))\vert\vert\xi\vert\vert^{-1}\tau_{\gamma_{\xi}}(\xi)
-\tau_{\gamma_{\xi}}(({\rm grad}\,r)_x)}
{\sqrt{\cos^2(\sqrt{\varepsilon}k_0br(x))+\vert\vert({\rm grad}\,r)_x\vert\vert^2}}
\leqno{(2.42)}$$
Also, according to $(2.40)$,  $(Z_{X,\xi})_{{\mathfrak p'}^{\perp}}(s)$ is descirbed as 
$$(Z_{X,\xi})_{{\mathfrak p'}^{\perp}}(s)=\sum_{i=2}^3\left(a^{i,+}_{X,\xi}(s)
({\bf C}^i_{X,\xi})_{{\mathfrak p'}^{\perp}_{k_0+k}}
+a^{i,-}_{X,\xi}(s)({\bf C}^i_{X,\xi})_{{\mathfrak p'}^{\perp}_{\vert k_0-k\vert}}\right)\leqno{(2.43)}$$
in terms of 
some explicitly described functions $a^{i,\pm}_{X,\xi}$ ($i=2,3$).  
From $(2.42),(2.43)$ and the definitions of ${\bf C}^i_{X,\xi}$ ($i=1,2$), we can derive 
$$\begin{array}{l}
\hspace{0.5truecm}\displaystyle{(\tau_{\gamma_{\xi}}\circ\tau_x)((Z_{X,\xi})_{\mathfrak p'^{\perp}}(1)))_T}\\
\displaystyle{=(\tau_{\gamma_{\xi}}\circ\tau_x)((Z_{X,\xi})_{\mathfrak p'^{\perp}}(1)))
-\bar g((\tau_{\gamma_{\xi}}\circ\tau_x)((Z_{X,\xi})_{\mathfrak p'^{\perp}}(1))),\,N_{\xi})N_{\xi}}\\
\displaystyle{=\sum_{i=2}^3
(\tau_{\gamma_{\xi}}\circ\tau_x)\left(a^{i,+}_{X,\xi}(1)({\bf C}^i_{X,\xi})_{{\mathfrak p'}^{\perp}_{k_0+k}}
+a^{i,-}_{X,\xi}(1)({\bf C}^i_{X,\xi})_{{\mathfrak p'}^{\perp}_{\vert k_0-k\vert}}\right)}\\
\hspace{0.5truecm}\displaystyle{
+\frac{(\sqrt{\varepsilon}k_0br(x))^2\cos(\sqrt{\varepsilon}k_0br(x))(Xr)}
{\cos^2(\sqrt{\varepsilon}k_0br(x))+\vert\vert({\rm grad}\,r)_x\vert\vert^2}}\\
\hspace{1truecm}\displaystyle{\times\left(\frac{2a^{2,-}_{X,\xi}(1)\sqrt{\varepsilon}kb}
{\cos^2(\sqrt{\varepsilon}k_0br(x))}+\frac{2a^{3,-}_{X,\xi}(1)\sqrt{\varepsilon}k_0b}
{\cos^2(\sqrt{\varepsilon}k_0br(x))}\right)}\\
\hspace{1truecm}\displaystyle{\times
\left(\frac{\cos(\sqrt{\varepsilon}k_0br(x))}{r(x)}\tau_{\gamma_{\xi}}(\xi)-\tau_{\gamma_{\xi}}(({\rm grad}\,r)_x)
\right),}
\end{array}\leqno{(2.44)}$$
where we use 
$$\bar g([[\tau_x^{-1}\xi,\tau_x^{-1}(({\rm grad}\,r)_x)],\tau_x^{-1}(X)]_{{\mathfrak p'}^{\perp}_{k_0+k}},
\tau_x^{-1}\xi)=0,\leqno{(2.45)}$$
$$\begin{array}{l}
\hspace{0.5truecm}\displaystyle{\bar g(
[[\tau_x^{-1}\xi,\tau_x^{-1}(({\rm grad}\,r)_x)],\tau_x^{-1}(X)]_{{\mathfrak p'}^{\perp}_{\vert k_0-k\vert}},
\tau_x^{-1}\xi)}\\
\displaystyle{=\bar g({\rm ad}(\tau_x^{-1}\xi)^2(\tau_x^{-1}(({\rm grad}\,r)_x)),
\tau_x^{-1}(X))}\\
\displaystyle{=-(\sqrt{\varepsilon}k_0br(x))^2(Xr),}
\end{array}\leqno{(2.46)}$$
$$\bar g([[\tau_x^{-1}\xi,\tau_x^{-1}(X)],\tau_x^{-1}(({\rm grad}\,r)_x)]_{{\mathfrak p'}^{\perp}_{k_0+k}},
\tau_x^{-1}\xi)=0\leqno{(2.47)}$$
and 
$$\begin{array}{l}
\hspace{0.5truecm}\displaystyle{\bar g(
[[\tau_x^{-1}\xi,\tau_x^{-1}(X)],\tau_x^{-1}(({\rm grad}\,r)_x)]_{{\mathfrak p'}^{\perp}_{\vert k_0-k\vert}},
\tau_x^{-1}\xi)}\\
\displaystyle{=\bar g({\rm ad}(\tau_x^{-1}\xi)^2(\tau_x^{-1}(({\rm grad}\,r)_x)),\tau_x^{-1}(X))}\\
\displaystyle{=-(\sqrt{\varepsilon}k_0br(x))^2(Xr).}
\end{array}\leqno{(2.48)}$$
From $(2.2),(2.8),(2.41)$ and $(2.44)$, we can derive the following relations:
$$\begin{array}{l}
\displaystyle{A\widetilde X_{\xi}=
-\frac{\sqrt{\varepsilon}kb\tan(\sqrt{\varepsilon}kbr(x))\cos(\sqrt{\varepsilon}k_0br(x))}
{\sqrt{\cos^2(\sqrt{\varepsilon}k_0br(x))+\vert\vert({\rm grad}\,r)_x\vert\vert^2}}\widetilde X_{\xi}}\\
\hspace{1truecm}\displaystyle{
-\frac{(Xr)\sqrt{\varepsilon}kb\tan(\sqrt{\varepsilon}kbr(x))
\cos(\sqrt{\varepsilon}k_0br(x))}
{(\cos^2(\sqrt{\varepsilon}k_0br(x))+\vert\vert({\rm grad}\,r)_x\vert\vert^2)^{3/2}}
\widetilde{(({\rm grad}\,r)_x)}_{\xi}}\\
\hspace{1truecm}\displaystyle{
+\frac{\cos(\sqrt{\varepsilon}k_0br(x))}
{\sqrt{\cos^2(\sqrt{\varepsilon}k_0br(x))+\vert\vert({\rm grad}\,r)_x\vert\vert^2}}}\\
\hspace{1truecm}\displaystyle{\times\left\{
\frac{g_F(\nabla^F_X{\rm grad}\,r,({\rm grad}\,r)_x)(\widetilde{({\rm grad}\,r)_x})_{\xi}}
{(\cos^2(\sqrt{\varepsilon}k_0br(x))+\vert\vert({\rm grad}\,r)_x\vert\vert^2)
\cos^2(\sqrt{\varepsilon}k_0br(x))}\right.}\\
\hspace{1.25truecm}\displaystyle{-\frac{1}{\cos^2(\sqrt{\varepsilon}k_0br(x))}
(\widetilde{\nabla^F_X{\rm grad}\,r})_{\xi}}\\
\hspace{1.25truecm}\displaystyle{-\sum_{i=2}^3
\left(\frac{a^{i,+}_{X,\xi}(1)\sqrt{\varepsilon}(k_0+k)br(x)}{\sin(\sqrt{\varepsilon}(k_0+k)br(x))}
\tau_x(({\bf C}^i_{X,\xi})_{{\mathfrak p'}^{\perp}_{k_0+k}})\right.}\\
\hspace{2.1truecm}\displaystyle{
\left.+\frac{a^{i,-}_{X,\xi}(1)\sqrt{\varepsilon}(k_0-k)br(x)}{\sin(\sqrt{\varepsilon}(k_0-k)br(x))}
\tau_x(({\bf C}^i_{X,\xi})_{{\mathfrak p'}^{\perp}_{\vert k_0-k\vert}})\right)}\\
\hspace{1.25truecm}\displaystyle{
-\frac{(\sqrt{\varepsilon}k_0br(x))^2\cos(\sqrt{\varepsilon}k_0br(x))(Xr)}
{\cos^2(\sqrt{\varepsilon}k_0br(x))+\vert\vert({\rm grad}\,r)_x\vert\vert^2}}\\
\hspace{1.75truecm}\displaystyle{\times\left(\frac{2a^{2,-}_{X,\xi}(1)\sqrt{\varepsilon}kb}
{\cos^2(\sqrt{\varepsilon}k_0br(x))}+\frac{2a^{3,-}_{X,\xi}(1)\sqrt{\varepsilon}k_0b}
{\cos^2(\sqrt{\varepsilon}k_0br(x))}\right)}\\
\hspace{1.75truecm}\displaystyle{\times
\left.\left(\frac{\cos(\sqrt{\varepsilon}k_0br(x))}{r(x)}\tau_{\gamma_{\xi}}(\xi)
-\frac{\sqrt{\varepsilon}k_0br(x)}{\sin(\sqrt{\varepsilon}k_0br(x))}(({\rm grad}\,r)_x)\right)\right\}.}
\end{array}\leqno{(2.49)}$$
Set $m^V:={\rm dim}\,{\cal V},\,\,m^H:={\rm dim}\,{\cal H},\,\,
m_k^V:={\rm dim}(\mathfrak p_{k\beta}\cap\mathfrak {p'}^{\perp})$ ($k=1,2$) and 
$m_k^H:={\rm dim}(\mathfrak p_{k\beta}\cap\mathfrak p')$ ($k=0,1,2$).  

From $(2.32)$ and $(2.49)$, we obtain the following description of the mean curvature vector $H$.  

\vspace{0.5truecm}

\noindent
{\bf Proposition 2.1.} {\sl Under the assumption of $\tau_x^{-1}(({\rm grad}\,r)_x)
\in\mathfrak p_{k_0\beta}\cap\mathfrak p'$, the mean curvature $H_{\xi}$ of $M$ at $\xi$ is 
described as 
$$\begin{array}{l}
\hspace{0truecm}\displaystyle{H_{\xi}
=\frac{\cos(\sqrt{\varepsilon}k_0br(x))}
{\sqrt{\cos^2(\sqrt{\varepsilon}k_0br(x))+\vert\vert({\rm grad}\,r)_{x}\vert\vert^2}}}\\
\hspace{0truecm}\displaystyle{
\times\left\{
\sum_{k=1}^2\frac{m_k^V\sqrt{\varepsilon}kb}{\tan(\sqrt{\varepsilon}kbr(x))}
-\sum_{k\in{\cal K}}m_k^H\sqrt{\varepsilon}kb\tan(\sqrt{\varepsilon}kbr(x))
-\frac{(\triangle_Fr)(x)}{\cos^2(\sqrt{\varepsilon}k_0br(x))}\right.}\\
\hspace{0.7truecm}\displaystyle{
-\frac{\vert\vert({\rm grad}\,r)_x\vert\vert^2\sqrt{\varepsilon}k_0b\tan(\sqrt{\varepsilon}k_0br(x))}
{\cos^2(\sqrt{\varepsilon}k_0br(x))+\vert\vert({\rm grad}\,r)_x\vert\vert^2}}\\
\hspace{0.7truecm}\displaystyle{\left.
+\frac{(\nabla^Fdr)_x(({\rm grad}\,r)_x,({\rm grad}\,r)_x)}
{\cos^2(\sqrt{\varepsilon}k_0br(x))
(\cos^2(\sqrt{\varepsilon}k_0br(x))+\vert\vert({\rm grad}\,r)_{x}\vert\vert^2)}\right\}.}
\end{array}\leqno{(2.50)}$$
}

\vspace{0.5truecm}

\noindent
{\it Proof.} Let $(e_1\,\cdots\,e_{m^H})$ be an orthonormal tangent frame of $B$ at $x$, where we take $e_1$ as 
$e_1=({\rm grad}\,r_t)_{\bar c_t(x_0)}/\vert\vert({\rm grad}\,r_t)_{\bar c_t(x_0)}\vert\vert$ in the case of 
$({\rm grad}\,r_t)_{\bar c_t(x_0)}\not=0$.  Also, let $(\hat e_1,\cdots,\hat e_{m^V})$ be an orthonormal base of 
${\cal V}_{\xi}$.  Then we have 
$$H_{\xi}=\sum_{i=1}^{m_V}\frac{g(A\hat e_i,\hat e_i)}
{\vert\vert\hat e_i\vert\vert^2}+\sum_{i=1}^{m_H}\frac{g(A((\widetilde e_i)_{\xi}),(\widetilde e_i)_{\xi})}
{\vert\vert(\widetilde e_i)_{\xi}\vert\vert^2}.$$
By substituting $(2.32)$ and $(2.49)$ into this relation and using $(2.3)$, we obtain the desired relation, 
where we note that the following relation holds:
$$\sum_{i=1}^{m^H}
\frac{\bar g(\widetilde{(\nabla^F_{e_i}{\rm grad}\,r)}_{\xi},(\widetilde{e_i})_{\xi})}
{\vert\vert(\widetilde{e_i})_{\xi}\vert\vert^2}=(\triangle_Fr)(x).\leqno{(2.51)}$$
\begin{flushright}q.e.d.\end{flushright}

\section{The volume element of a tube over a reflective submanifold}
We shall use the notations in Introduction and the previous section.  
In this section, we shall calculate the volume element of $M$.  
First we recall the description of the Jacobi field in a symmetric space $\overline M$.  
The Jacobi field $J$ along the geodesic $\gamma$ in $\overline M$ is described as 
$$J(t)=\tau_{\gamma\vert_{[0,t]}}\left(D^{co}_{t\gamma'(0)}(J(0))
+tD^{si}_{t\gamma'(0)}(J'(0))\right).\leqno{(3.1)}$$
Let $r_F$ be as in Introduction.  By using $(3.1)$, we can derive 
$$r_{co}(\gamma)=\left\{
\begin{array}{ll}
\displaystyle{\frac{\pi}{2\sqrt{\varepsilon}b}} & ({\rm when}\,\,m^V_2\not=0)\\
\displaystyle{\frac{\pi}{\sqrt{\varepsilon}b}} & ({\rm when}\,\,m^V_2=0)
\end{array}\right.$$
and 
$$r_{fo}(\gamma)=\min\left\{\left.\frac{\pi}{2\sqrt{\varepsilon}kb}\,\right\vert\,k\in{\cal K}\right\}.$$
Hence we obtain 
$$r_F=\left\{
\begin{array}{ll}
\displaystyle{\min\left\{\left.\frac{\pi}{2\sqrt{\varepsilon}kb}\,\right\vert\,k\in{\cal K}\cup\{1\}\right\}} & 
({\rm when}\,\,m^V_2\not=0)\\
\displaystyle{\min\left\{\left.\frac{\pi}{2\sqrt{\varepsilon}kb}\,\right\vert\,k\in{\cal K}\cup\left\{\frac{1}{2}
\right\}\,\right\}} & ({\rm when}\,\,m^V_2=0).
\end{array}\right.\leqno{(3.2)}$$
Fix $\xi\in M\cap T^{\perp}_xB$ and $X\in T_xB$.  
Without loss of generality, we may assume that $\tau_x^{-1}\xi\in\mathfrak b$.  
By the assumption (C2), we may assume that 
$\tau_x^{-1}(({\rm grad}\,r)_x)\in\mathfrak p_{k_0\beta}\cap\mathfrak p'$ for some $k_0\in{\cal K}\cup\{0\}$.  
Let $\widetilde S(x,r(x))$ be the hypersphere of radius $r(x)$ in $T^{\perp}_xB$ centered 
the origin and $S(x,r(x))$ the geodesic hypersphere of radius $r(x)$ 
in $F_x^{\perp}:=\exp^{\perp}(T^{\perp}_xB)$ centered $x$.  
Denote by $dv_{(\cdot)}$ the volume element of the induced metric on $(\cdot)$.  
Take $v\in T_{\xi}\widetilde S(x,r(x))$.  
Define a function $\psi_r$ over $B$ by 
$$\begin{array}{l}
\displaystyle{\psi_r(x):=
\left(\mathop{\Pi}_{k=1}^2
\left(\frac{\sin(\sqrt{\varepsilon}kbr(x))}
{\sqrt{\varepsilon}kb}\right)^{m_k^V}\right)
\left(\mathop{\Pi}_{k\in{\cal K}\setminus\{k_0\}}\cos^{m^H_k}(\sqrt{\varepsilon}kbr(x))\right)}\\
\hspace{1.7truecm}\displaystyle{\times\cos^{m^H_{k_0}-1}(\sqrt{\varepsilon}k_0br(x))
\sqrt{\cos^2(\sqrt{\varepsilon}k_0br(x))+\vert\vert({\rm grad}\,r)_x\vert\vert^2}.}
\end{array}
\leqno{(3.3)}$$
We have $r<r_F$ because $\exp^{\perp}\vert_{t_r(B)}$ is not an immersion when $r(x)=r_F$ for some 
$x\in B$.  Hence $\psi_r$ is positive by $(3.2)$.  
Define a function $\rho_r$ over $B$ by 
$$\begin{array}{l}
\displaystyle{\rho_r(x):=\frac{\cos(\sqrt{\varepsilon}k_0br(x))}
{\sqrt{\cos^2(\sqrt{\varepsilon}k_0br(x))+\vert\vert({\rm grad}\,r)_x\vert\vert^2}}}\\
\hspace{1.25truecm}\displaystyle{\times\left(
\sum_{k=1}^2\frac{m_k^V\sqrt{\varepsilon}kb}{\tan(\sqrt{\varepsilon}kbr(x))}
-\sum_{k\in{\cal K}}m_k^H\sqrt{\varepsilon}kb\tan(\sqrt{\varepsilon}kbr(x))\right.}\\
\hspace{2truecm}\displaystyle{\left.
+\frac{\vert\vert({\rm grad}\,r)_x\vert\vert^2\sqrt{\varepsilon}k_0b\tan(\sqrt{\varepsilon}k_0br(x))}
{\cos^2(\sqrt{\varepsilon}k_0br(x))+\vert\vert({\rm grad}\,r)_x\vert\vert^2}\right).}
\end{array}\leqno{(3.4)}$$
According to $(2.50)$, we have 
$$\begin{array}{l}
\displaystyle{H_{\xi}=\rho_r(x)-\frac{(\triangle_Fr)(x)}
{\cos(\sqrt{\varepsilon}k_0br(x))
\sqrt{\cos^2(\sqrt{\varepsilon}k_0br(x))+\vert\vert({\rm grad}\,r)_x\vert\vert^2}}}\\
\hspace{1.2truecm}\displaystyle{
+\frac{(\nabla^Fdr)_x(({\rm grad}\,r)_x,({\rm grad}\,r)_x)}
{\cos(\sqrt{\varepsilon}k_0br(x))
(\cos^2(\sqrt{\varepsilon}k_0br(x))+\vert\vert({\rm grad}\,r)_x\vert\vert^2)^{3/2}}.}
\end{array}\leqno{(3.5)}$$
From $(2.2)$ and $(2.3)$, we can derive the following relation for the volume element of $M$.  

\vspace{0.5truecm}

\noindent
{\bf Proposition 3.1.} {\sl The volume element $dv_M$ is given by 
$$(dv_M)_{\xi}=\psi_r(x)\left(
((\exp^{\perp}\vert_{\widetilde S(x,r(x))})^{-1})^{\ast}
dv_{\widetilde S(x,r(x))}\wedge(\pi\vert_M)^{\ast}dv_F\right),
\leqno{(3.6)}$$
where $\xi\in M\cap T^{\perp}_xB$.}

\vspace{0.5truecm}

From $(3.6)$, we can derive the following relation for the volume of $M$.  

\vspace{0.5truecm}

\noindent
{\bf Proposition 3.2.} {\sl The volume ${\rm Vol}(M)$ of $M$ and 
the average mean curvature $\overline H$ of $M$ are given by 
$${\rm Vol}(M)=v_{m^V}\int_Br^{m^V}\psi_rdv_F\leqno{(3.7)}$$
and 
$$\begin{array}{l}
\displaystyle{\overline H=\frac{1}{\int_Br^{m^V}\psi_rdv_F}}\\
\hspace{0.3truecm}\displaystyle{\times
\int_{x\in B}r^{m^V}\left(
\rho_r(x)-\frac{(\triangle_Fr)(x)}{\cos(\sqrt{\varepsilon}k_0br(x))
\sqrt{\cos^2(\sqrt{\varepsilon}k_0br(x))+\vert\vert({\rm grad}\,r)_x\vert\vert^2}}\right.}\\
\hspace{1.7truecm}\displaystyle{
\left.+\frac{(\nabla^Fdr)_x(({\rm grad}\,r)_x,({\rm grad}\,r)_x)}
{\cos(\sqrt{\varepsilon}k_0br(x))
(\cos^2(\sqrt{\varepsilon}k_0br(x))+\vert\vert({\rm grad}\,r)_x\vert\vert^2)^{3/2}}\right)
\psi_r\,dv_F},
\end{array}\leqno{(3.8)}$$
respectively, where $\xi\in M\cap T^{\perp}_xB$ and 
$v_{m^V}$ is the volume of the $m^V$-dimensional Euclidean unit sphere.}

\section{The evolution of the radius function} 
Let $F,\,B,\,M=t_{r_0}(B)$ and $f$ be as in Setting (S) of Introduction.  
Assume that the volume-preserving mean curvature flow $f_t$ ($t\in[0,T)$) starting from $f$ satisfies 
the conditions (C1) and (C2).  
We use the notations in Introduction and Sections 1-3.  
Denote by $S^{\perp}B$ the unit normal bundle of of $B$ and 
$S^{\perp}_xB$ the fibre of this bundle over $x\in B$.  
Define a positive-valued function $\widehat r_t:M\to{\Bbb R}$ ($t\in[0,T)$) and 
a map 
$w^1_t:M\to S^{\perp}B$ 
($t\in[0,T)$) by 
$f_t(\xi)=\exp^{\perp}(\widehat r_t(\xi)w^1_t(\xi))$ ($\xi\in M$).  
Also, define a map $c_t:M\to B$ by $c_t(\xi):=\pi(w^1_t(\xi))$ 
($\xi\in M$) and a map $w_t:M\to T^{\perp}B\,\,(t\in[0,T))$ 
by 
$w_t(\xi):=\widehat r_t(\xi)w^1_t(\xi)\,\,(\xi\in M)$.  
Here we note that $c_t$ is surjective by the boundary condition in Theorem A, 
$\widehat r_0(\xi)=r_0(\pi(\xi))$ 
and that $c_0(\xi)=\pi(\xi)$ ($\xi\in M$).  
Define a function $\bar r_t$ over $B$ by 
$\bar r_t(x):=\widehat r_t(\xi)\,\,(x\in B)$ and a map 
$\bar c_t:B\to B$ by $\bar c_t(x):=c_t(\xi)$ $(x\in B)$, where $\xi$ is an 
arbitrary element of $M\cap S^{\perp}_xB$.  It is clear that they are 
well-defined.  This map $\bar c_t$ is not necessarily a diffeomorphism.  In particular, 
if $\bar c_t$ is a diffeomorphism, then $M_t:=f_t(M)$ is equal to the tube 
$\exp^{\perp}(t_{r_t}(B))$, where $r_t:=\bar r_t\circ{\bar c}_t^{-1}$.  
It is easy to show that, if $c_t(\xi_1)=c_t(\xi_2)$, then $\widehat r_t(\xi_1)
=\widehat r_t(\xi_2)$ and $\pi(\xi_1)=\pi(\xi_2)$ hold.  
Also, let $a:M\times[0,T)\to G$ be a smooth map with $a(\xi,t)p_0=c_t(\xi)\,\,
(\xi,t)\in M\times[0,T)$.  
In this section, we shall calculate the evolution equations for the functions $r_t$ and $\widehat r_t$.  
Define $\widetilde f:M\times [0,T)\to \overline M,\,\,
r:B\times [0,T)\to{\Bbb R},\,\,
w^1:M\times [0,T)\to M$ and 
$c:M\times [0,T)\to B$ by 
$\widetilde f(\xi,t):=f_t(\xi),\,\,r(x,t):=r_t(x),\,\,
w^1(\xi,t):=w^1_t(\xi),\,w(\xi,t):=w_t(\xi)$ and 
$c(\xi,t):=c_t(\xi)$, respectively, where 
$\xi\in M,x\in B$ and $t\in[0,T)$ (see Figure 7).  
Fix $(\xi_0,t_0)\in M\times[0,T)$ and set $x_0:=c(\xi_0,t_0)$.  
According to the condition (C2), we may assume that 
$\tau^{-1}_{x_0}(({\rm grad}r_{t_0})_{x_0})$ belong to 
$\mathfrak p_{k_0\beta}$ for some $k_0\in{\cal K}\cup\{0\}$, where ${\cal K}$ and $\mathfrak p_{k_0\beta}$ are 
the quantities defined as in Section 2 for the maximal ableian subspace 
$\mathfrak b:={\rm Span}\{\tau_{x_0}^{-1}(w(\xi_0,t_0))\}$ of $\mathfrak p'^{\perp}$.  
Clearly we have 
$$
\widetilde f_{\ast}\left(\frac{\partial}{\partial t}\right)_{(\xi_0,t_0)}
=\left.\frac{d}{dt}\right\vert_{t=t_0}
\widetilde f(\xi_0,t)
=\left.\frac{d}{dt}\right\vert_{t=t_0}
\exp^{\perp}(w(\xi_0,t)).
\leqno{(4.1)}$$

\vspace{0.25truecm}

\centerline{
\unitlength 0.1in
\begin{picture}( 52.3600, 20.4000)(  8.9500,-32.1000)
%
\special{pn 20}%
\special{sh 1}%
\special{ar 4130 2190 10 10 0  6.28318530717959E+0000}%
\special{sh 1}%
\special{ar 4130 2190 10 10 0  6.28318530717959E+0000}%
%
\special{pn 20}%
\special{sh 1}%
\special{ar 4134 2118 10 10 0  6.28318530717959E+0000}%
\special{sh 1}%
\special{ar 4134 2118 10 10 0  6.28318530717959E+0000}%
%
\special{pn 20}%
\special{sh 1}%
\special{ar 5130 2190 10 10 0  6.28318530717959E+0000}%
\special{sh 1}%
\special{ar 5130 2190 10 10 0  6.28318530717959E+0000}%
%
\special{pn 13}%
\special{pa 4156 2118}%
\special{pa 4986 2124}%
\special{fp}%
\special{sh 1}%
\special{pa 4986 2124}%
\special{pa 4920 2104}%
\special{pa 4934 2124}%
\special{pa 4920 2144}%
\special{pa 4986 2124}%
\special{fp}%
%
\special{pn 8}%
\special{pa 5718 2244}%
\special{pa 4996 2116}%
\special{dt 0.045}%
\special{sh 1}%
\special{pa 4996 2116}%
\special{pa 5058 2148}%
\special{pa 5048 2126}%
\special{pa 5064 2108}%
\special{pa 4996 2116}%
\special{fp}%
%
\special{pn 8}%
\special{pa 5878 1460}%
\special{pa 5056 1728}%
\special{dt 0.045}%
\special{sh 1}%
\special{pa 5056 1728}%
\special{pa 5126 1726}%
\special{pa 5108 1712}%
\special{pa 5114 1688}%
\special{pa 5056 1728}%
\special{fp}%
%
\special{pn 13}%
\special{pa 4130 2190}%
\special{pa 5110 2190}%
\special{fp}%
\special{sh 1}%
\special{pa 5110 2190}%
\special{pa 5044 2170}%
\special{pa 5058 2190}%
\special{pa 5044 2210}%
\special{pa 5110 2190}%
\special{fp}%
\put(54.8000,-23.6000){\makebox(0,0)[lt]{$f(\xi)$}}%
\put(57.6400,-21.9400){\makebox(0,0)[lt]{$f_t(\xi)$}}%
\put(31.3000,-25.5000){\makebox(0,0)[rt]{$c_0(\xi)=\bar c_0(x)=x$}}%
\put(28.7500,-25.0100){\makebox(0,0)[rb]{$c_t(\xi)=\bar c_t(x)$}}%
\put(51.3000,-13.4000){\makebox(0,0)[lb]{$B$}}%
\put(59.3300,-14.6700){\makebox(0,0)[lb]{$f_t(M)$}}%
\put(60.2700,-19.5800){\makebox(0,0)[lb]{$f(M)$}}%
%
\special{pn 8}%
\special{pa 5044 1810}%
\special{pa 5038 1840}%
\special{pa 5034 1872}%
\special{pa 5030 1904}%
\special{pa 5026 1936}%
\special{pa 5020 1968}%
\special{pa 5012 1998}%
\special{pa 5008 2030}%
\special{pa 5004 2062}%
\special{pa 4996 2094}%
\special{pa 4990 2124}%
\special{pa 4982 2156}%
\special{pa 4978 2188}%
\special{pa 4968 2218}%
\special{pa 4960 2248}%
\special{pa 4954 2280}%
\special{pa 4944 2310}%
\special{pa 4934 2342}%
\special{pa 4926 2372}%
\special{pa 4916 2402}%
\special{pa 4904 2432}%
\special{pa 4896 2462}%
\special{pa 4884 2492}%
\special{pa 4872 2522}%
\special{pa 4860 2552}%
\special{pa 4848 2572}%
\special{sp}%
%
\special{pn 8}%
\special{ar 4132 2104 856 112  6.2831853 6.2831853}%
\special{ar 4132 2104 856 112  0.0000000 3.1415927}%
%
\special{pn 8}%
\special{ar 4120 2120 862 124  3.1415927 3.1659581}%
\special{ar 4120 2120 862 124  3.2390546 3.2634201}%
\special{ar 4120 2120 862 124  3.3365165 3.3608820}%
\special{ar 4120 2120 862 124  3.4339784 3.4583439}%
\special{ar 4120 2120 862 124  3.5314404 3.5558059}%
\special{ar 4120 2120 862 124  3.6289023 3.6532678}%
\special{ar 4120 2120 862 124  3.7263642 3.7507297}%
\special{ar 4120 2120 862 124  3.8238262 3.8481916}%
\special{ar 4120 2120 862 124  3.9212881 3.9456536}%
\special{ar 4120 2120 862 124  4.0187500 4.0431155}%
\special{ar 4120 2120 862 124  4.1162119 4.1405774}%
\special{ar 4120 2120 862 124  4.2136739 4.2380394}%
\special{ar 4120 2120 862 124  4.3111358 4.3355013}%
\special{ar 4120 2120 862 124  4.4085977 4.4329632}%
\special{ar 4120 2120 862 124  4.5060597 4.5304251}%
\special{ar 4120 2120 862 124  4.6035216 4.6278871}%
\special{ar 4120 2120 862 124  4.7009835 4.7253490}%
\special{ar 4120 2120 862 124  4.7984454 4.8228109}%
\special{ar 4120 2120 862 124  4.8959074 4.9202729}%
\special{ar 4120 2120 862 124  4.9933693 5.0177348}%
\special{ar 4120 2120 862 124  5.0908312 5.1151967}%
\special{ar 4120 2120 862 124  5.1882932 5.2126586}%
\special{ar 4120 2120 862 124  5.2857551 5.3101206}%
\special{ar 4120 2120 862 124  5.3832170 5.4075825}%
\special{ar 4120 2120 862 124  5.4806789 5.5050444}%
\special{ar 4120 2120 862 124  5.5781409 5.6025064}%
\special{ar 4120 2120 862 124  5.6756028 5.6999683}%
\special{ar 4120 2120 862 124  5.7730647 5.7974302}%
\special{ar 4120 2120 862 124  5.8705267 5.8948921}%
\special{ar 4120 2120 862 124  5.9679886 5.9923541}%
\special{ar 4120 2120 862 124  6.0654505 6.0898160}%
\special{ar 4120 2120 862 124  6.1629125 6.1872779}%
\special{ar 4120 2120 862 124  6.2603744 6.2831853}%
%
\special{pn 8}%
\special{pa 3200 1810}%
\special{pa 3208 1840}%
\special{pa 3210 1872}%
\special{pa 3216 1904}%
\special{pa 3220 1936}%
\special{pa 3226 1968}%
\special{pa 3232 1998}%
\special{pa 3236 2030}%
\special{pa 3240 2062}%
\special{pa 3248 2092}%
\special{pa 3254 2124}%
\special{pa 3260 2156}%
\special{pa 3268 2188}%
\special{pa 3276 2218}%
\special{pa 3284 2248}%
\special{pa 3290 2280}%
\special{pa 3300 2310}%
\special{pa 3310 2340}%
\special{pa 3318 2372}%
\special{pa 3330 2402}%
\special{pa 3340 2432}%
\special{pa 3348 2462}%
\special{pa 3360 2492}%
\special{pa 3372 2522}%
\special{pa 3386 2550}%
\special{pa 3396 2572}%
\special{sp}%
%
\special{pn 8}%
\special{pa 3460 3060}%
\special{pa 3460 3028}%
\special{pa 3460 2996}%
\special{pa 3462 2964}%
\special{pa 3460 2932}%
\special{pa 3460 2900}%
\special{pa 3460 2868}%
\special{pa 3454 2838}%
\special{pa 3452 2806}%
\special{pa 3450 2774}%
\special{pa 3444 2742}%
\special{pa 3436 2712}%
\special{pa 3428 2680}%
\special{pa 3422 2650}%
\special{pa 3414 2618}%
\special{pa 3406 2588}%
\special{pa 3390 2562}%
\special{pa 3390 2562}%
\special{sp}%
%
\special{pn 20}%
\special{sh 1}%
\special{ar 4998 2108 10 10 0  6.28318530717959E+0000}%
\special{sh 1}%
\special{ar 4998 2108 10 10 0  6.28318530717959E+0000}%
%
\special{pn 8}%
\special{pa 3160 2570}%
\special{pa 4118 2198}%
\special{dt 0.045}%
\special{sh 1}%
\special{pa 4118 2198}%
\special{pa 4048 2204}%
\special{pa 4068 2216}%
\special{pa 4062 2240}%
\special{pa 4118 2198}%
\special{fp}%
%
\special{pn 8}%
\special{ar 3468 1652 278 536  2.8443036 3.1415927}%
%
\special{pn 8}%
\special{ar 4788 1648 278 536  6.2831853 6.2831853}%
\special{ar 4788 1648 278 536  0.0000000 0.2972890}%
%
\special{pn 8}%
\special{ar 4150 3070 458 80  6.2831853 6.2831853}%
\special{ar 4150 3070 458 80  0.0000000 3.1415927}%
%
\special{pn 8}%
\special{ar 4148 3082 454 66  3.1415927 3.1878354}%
\special{ar 4148 3082 454 66  3.3265638 3.3728065}%
\special{ar 4148 3082 454 66  3.5115349 3.5577776}%
\special{ar 4148 3082 454 66  3.6965059 3.7427487}%
\special{ar 4148 3082 454 66  3.8814770 3.9277198}%
\special{ar 4148 3082 454 66  4.0664481 4.1126909}%
\special{ar 4148 3082 454 66  4.2514192 4.2976620}%
\special{ar 4148 3082 454 66  4.4363903 4.4826331}%
\special{ar 4148 3082 454 66  4.6213614 4.6676042}%
\special{ar 4148 3082 454 66  4.8063325 4.8525753}%
\special{ar 4148 3082 454 66  4.9913036 5.0375464}%
\special{ar 4148 3082 454 66  5.1762747 5.2225175}%
\special{ar 4148 3082 454 66  5.3612458 5.4074886}%
\special{ar 4148 3082 454 66  5.5462169 5.5924597}%
\special{ar 4148 3082 454 66  5.7311880 5.7774308}%
\special{ar 4148 3082 454 66  5.9161591 5.9624019}%
\special{ar 4148 3082 454 66  6.1011302 6.1473730}%
%
\special{pn 8}%
\special{ar 4126 3082 680 128  6.2831853 6.2831853}%
\special{ar 4126 3082 680 128  0.0000000 3.1415927}%
%
\special{pn 8}%
\special{ar 4114 3088 682 128  3.1415927 3.1712589}%
\special{ar 4114 3088 682 128  3.2602577 3.2899239}%
\special{ar 4114 3088 682 128  3.3789227 3.4085889}%
\special{ar 4114 3088 682 128  3.4975877 3.5272540}%
\special{ar 4114 3088 682 128  3.6162527 3.6459190}%
\special{ar 4114 3088 682 128  3.7349177 3.7645840}%
\special{ar 4114 3088 682 128  3.8535828 3.8832490}%
\special{ar 4114 3088 682 128  3.9722478 4.0019140}%
\special{ar 4114 3088 682 128  4.0909128 4.1205791}%
\special{ar 4114 3088 682 128  4.2095778 4.2392441}%
\special{ar 4114 3088 682 128  4.3282428 4.3579091}%
\special{ar 4114 3088 682 128  4.4469079 4.4765741}%
\special{ar 4114 3088 682 128  4.5655729 4.5952391}%
\special{ar 4114 3088 682 128  4.6842379 4.7139041}%
\special{ar 4114 3088 682 128  4.8029029 4.8325692}%
\special{ar 4114 3088 682 128  4.9215679 4.9512342}%
\special{ar 4114 3088 682 128  5.0402330 5.0698992}%
\special{ar 4114 3088 682 128  5.1588980 5.1885642}%
\special{ar 4114 3088 682 128  5.2775630 5.3072292}%
\special{ar 4114 3088 682 128  5.3962280 5.4258943}%
\special{ar 4114 3088 682 128  5.5148930 5.5445593}%
\special{ar 4114 3088 682 128  5.6335580 5.6632243}%
\special{ar 4114 3088 682 128  5.7522231 5.7818893}%
\special{ar 4114 3088 682 128  5.8708881 5.9005543}%
\special{ar 4114 3088 682 128  5.9895531 6.0192194}%
\special{ar 4114 3088 682 128  6.1082181 6.1378844}%
\special{ar 4114 3088 682 128  6.2268831 6.2565494}%
%
\special{pn 8}%
\special{ar 4116 1638 928 122  6.2831853 6.2831853}%
\special{ar 4116 1638 928 122  0.0000000 3.1415927}%
%
\special{pn 8}%
\special{ar 4128 1646 928 122  3.1415927 6.2831853}%
%
\special{pn 8}%
\special{ar 4140 1624 1358 242  6.2831853 6.2831853}%
\special{ar 4140 1624 1358 242  0.0000000 3.1415927}%
%
\special{pn 8}%
\special{ar 4140 1630 1358 228  3.1415927 6.2831853}%
%
\special{pn 8}%
\special{pa 5090 1300}%
\special{pa 4128 1686}%
\special{dt 0.045}%
\special{sh 1}%
\special{pa 4128 1686}%
\special{pa 4196 1680}%
\special{pa 4178 1666}%
\special{pa 4182 1644}%
\special{pa 4128 1686}%
\special{fp}%
%
\special{pn 8}%
\special{pa 2922 2402}%
\special{pa 4116 2116}%
\special{dt 0.045}%
\special{sh 1}%
\special{pa 4116 2116}%
\special{pa 4048 2112}%
\special{pa 4064 2128}%
\special{pa 4056 2150}%
\special{pa 4116 2116}%
\special{fp}%
%
\special{pn 8}%
\special{pa 6070 1624}%
\special{pa 4656 2116}%
\special{dt 0.045}%
\special{sh 1}%
\special{pa 4656 2116}%
\special{pa 4726 2114}%
\special{pa 4706 2098}%
\special{pa 4712 2076}%
\special{pa 4656 2116}%
\special{fp}%
\put(61.3100,-16.4700){\makebox(0,0)[lb]{$w_t(\xi)$}}%
\put(54.6500,-26.3800){\makebox(0,0)[lt]{$w_0(\xi)=\xi$}}%
%
\special{pn 8}%
\special{pa 4804 3090}%
\special{pa 4802 3058}%
\special{pa 4804 3026}%
\special{pa 4802 2994}%
\special{pa 4804 2962}%
\special{pa 4804 2930}%
\special{pa 4808 2900}%
\special{pa 4808 2868}%
\special{pa 4812 2836}%
\special{pa 4814 2804}%
\special{pa 4816 2772}%
\special{pa 4820 2740}%
\special{pa 4826 2708}%
\special{pa 4832 2678}%
\special{pa 4836 2646}%
\special{pa 4844 2614}%
\special{pa 4850 2584}%
\special{pa 4860 2558}%
\special{sp}%
%
\special{pn 8}%
\special{ar 2530 3070 1160 1050  5.2228136 6.2831853}%
%
\special{pn 8}%
\special{ar 4170 1630 1390 830  2.4232081 3.1325503}%
%
\special{pn 8}%
\special{ar 5750 3082 1160 1050  3.1415927 4.2019644}%
%
\special{pn 8}%
\special{ar 4120 1632 1390 830  0.0090423 0.7183846}%
%
\special{pn 8}%
\special{ar 4150 2200 960 100  6.2831853 6.2831853}%
\special{ar 4150 2200 960 100  0.0000000 3.1415927}%
%
\special{pn 8}%
\special{ar 4160 2200 970 110  3.1415927 3.1638149}%
\special{ar 4160 2200 970 110  3.2304815 3.2527038}%
\special{ar 4160 2200 970 110  3.3193704 3.3415927}%
\special{ar 4160 2200 970 110  3.4082593 3.4304815}%
\special{ar 4160 2200 970 110  3.4971482 3.5193704}%
\special{ar 4160 2200 970 110  3.5860371 3.6082593}%
\special{ar 4160 2200 970 110  3.6749260 3.6971482}%
\special{ar 4160 2200 970 110  3.7638149 3.7860371}%
\special{ar 4160 2200 970 110  3.8527038 3.8749260}%
\special{ar 4160 2200 970 110  3.9415927 3.9638149}%
\special{ar 4160 2200 970 110  4.0304815 4.0527038}%
\special{ar 4160 2200 970 110  4.1193704 4.1415927}%
\special{ar 4160 2200 970 110  4.2082593 4.2304815}%
\special{ar 4160 2200 970 110  4.2971482 4.3193704}%
\special{ar 4160 2200 970 110  4.3860371 4.4082593}%
\special{ar 4160 2200 970 110  4.4749260 4.4971482}%
\special{ar 4160 2200 970 110  4.5638149 4.5860371}%
\special{ar 4160 2200 970 110  4.6527038 4.6749260}%
\special{ar 4160 2200 970 110  4.7415927 4.7638149}%
\special{ar 4160 2200 970 110  4.8304815 4.8527038}%
\special{ar 4160 2200 970 110  4.9193704 4.9415927}%
\special{ar 4160 2200 970 110  5.0082593 5.0304815}%
\special{ar 4160 2200 970 110  5.0971482 5.1193704}%
\special{ar 4160 2200 970 110  5.1860371 5.2082593}%
\special{ar 4160 2200 970 110  5.2749260 5.2971482}%
\special{ar 4160 2200 970 110  5.3638149 5.3860371}%
\special{ar 4160 2200 970 110  5.4527038 5.4749260}%
\special{ar 4160 2200 970 110  5.5415927 5.5638149}%
\special{ar 4160 2200 970 110  5.6304815 5.6527038}%
\special{ar 4160 2200 970 110  5.7193704 5.7415927}%
\special{ar 4160 2200 970 110  5.8082593 5.8304815}%
\special{ar 4160 2200 970 110  5.8971482 5.9193704}%
\special{ar 4160 2200 970 110  5.9860371 6.0082593}%
\special{ar 4160 2200 970 110  6.0749260 6.0971482}%
\special{ar 4160 2200 970 110  6.1638149 6.1860371}%
\special{ar 4160 2200 970 110  6.2527038 6.2749260}%
%
\special{pn 8}%
\special{pa 4130 1630}%
\special{pa 4130 3090}%
\special{fp}%
%
\special{pn 8}%
\special{pa 5460 2360}%
\special{pa 5140 2200}%
\special{dt 0.045}%
\special{sh 1}%
\special{pa 5140 2200}%
\special{pa 5192 2248}%
\special{pa 5188 2224}%
\special{pa 5210 2212}%
\special{pa 5140 2200}%
\special{fp}%
%
\special{pn 8}%
\special{pa 6000 1900}%
\special{pa 5390 1960}%
\special{dt 0.045}%
\special{sh 1}%
\special{pa 5390 1960}%
\special{pa 5458 1974}%
\special{pa 5444 1956}%
\special{pa 5454 1934}%
\special{pa 5390 1960}%
\special{fp}%
%
\special{pn 8}%
\special{pa 5430 2670}%
\special{pa 4920 2190}%
\special{dt 0.045}%
\special{sh 1}%
\special{pa 4920 2190}%
\special{pa 4956 2250}%
\special{pa 4960 2228}%
\special{pa 4982 2222}%
\special{pa 4920 2190}%
\special{fp}%
\end{picture}%
\hspace{3.5truecm}}

\vspace{0.25truecm}

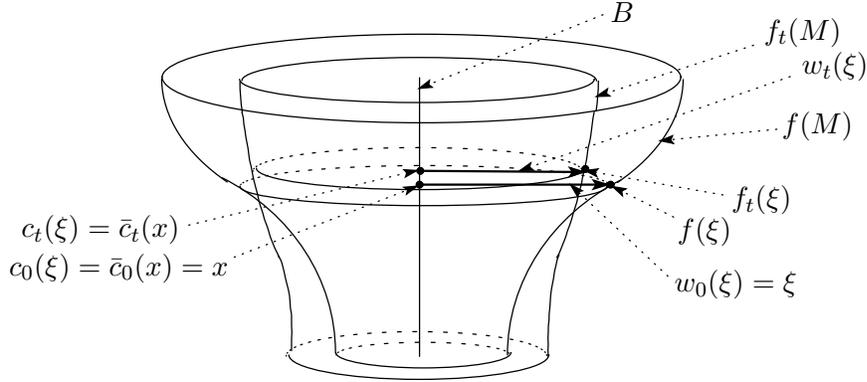
\centerline{{\bf Figure 7.}}

\vspace{0.5truecm}

\noindent
Let $J$ be the Jacobi field along the geodesic 
$\gamma_{w(\xi_0,t_0)}$ (of direction 
$w(\xi_0,t_0)$) with 
$J(0)=\displaystyle{c_{\ast}
\left(\frac{\partial}{\partial t}\right)_{(\xi_0,t_0)}}$ 
and $J'(0)=\displaystyle{
\left.\frac{\overline{\nabla}}{\partial t}\right\vert_{t=t_0}
w(\xi_0,\cdot)}$, where $\displaystyle{\frac{\overline{\nabla}}{\partial t}}$ is the pull-back connection of 
$\overline{\nabla}$ by $t\mapsto c(\xi_0,t)$.  
This Jacobi field $J$ is described as 
$$\begin{array}{l}
\displaystyle{J(s)=\tau_{\gamma_{w(\xi_0,t_0)}\vert_{[0,s]}}
\left(D^{co}_{s w(\xi_0,t_0)}
\left(c_{\ast}\left(\frac{\partial}{\partial t}\right)_{(\xi_0,t_0)}
\right)\right)}\\
\hspace{1.45truecm}\displaystyle{+\tau_{\gamma_{w(\xi_0,t_0)}\vert_{[0,s]}}
\left(sD^{si}_{sw(\xi_0,t_0)}
\left(\left.\frac{\overline{\nabla}}{\partial t}\right\vert_{t=t_0}
w(\xi_0,\cdot)\right)\right).}
\end{array}$$
According to $(4.1)$, we have 
$$\begin{array}{l}
\displaystyle{\widetilde f_{\ast}\left(\frac{\partial}{\partial t}
\right)_{(\xi_0,t_0)}=J(1)}\\
\hspace{2.4truecm}\displaystyle{=\tau_{\gamma_{w(\xi_0,t_0)}
\vert_{[0,1]}}
\left(D^{co}_{w(\xi_0,t_0)}\left(
c_{\ast}\left(\frac{\partial}{\partial t}\right)_{(\xi_0,t_0)}\right)\right)}\\
\hspace{2.9truecm}\displaystyle{+\tau_{\gamma_{w(\xi_0,t_0)}\vert_{[0,1]}}
\left(D^{si}_{w(\xi_0,t_0)}
\left(\left.\frac{\overline{\nabla}}{\partial t}\right\vert_{t=t_0}
w(\xi_0,\cdot)\right)\right).}
\end{array}
\leqno{(4.2)}$$
On the other hand, we have 
$$\begin{array}{l}
\hspace{0.5truecm}\displaystyle{
\left.\frac{\overline{\nabla}}{\partial t}\right\vert_{t=t_0}
w(\xi_0,\cdot)}\\
\displaystyle{=\left.\frac{d\widehat r(\xi_0,t)}{dt}\right\vert_{t=t_0}
w^1(\xi_0,t_0)+\widehat r(\xi_0,t_0)\left.\frac{\overline{\nabla}}{\partial t}\right\vert_{t=t_0}
w^1(\xi_0,t).}
\end{array}\leqno{(4.3)}$$
From $(4.2)$ and $(4.3)$, we have 
$$\begin{array}{r}
\displaystyle{\widetilde f_{\ast}\left(\frac{\partial}{\partial t}\right)_{(\xi_0,t_0)}
\equiv
\left.\frac{d\widehat r(\xi_0,t)}{dt}\right\vert_{t=t_0}
\tau_{\gamma_{w(\xi_0,t_0)}\vert_{[0,1]}}(w^1(\xi_0,t_0))}\\
\displaystyle{({\rm mod}\,\,\,\,
{\rm Span}\{\tau_{\gamma_{w(\xi_0,t_0)}\vert_{[0,1]}}
(w^1(\xi_0,t_0))\}^{\perp}).}
\end{array}\leqno{(4.4)}$$

\vspace{0.5truecm}

\noindent
{\bf Notation.} Set 
$$T_1:=\mathop{\sup}\{t'\in[0,T)\,\vert\,\,M_t:=f_t(M)\,\,(0\leq t\leq t')\,:\,{\rm tubes}\,\,{\rm over}\,\,B\}.$$
(Note that $\bar c_t$ ($0\leq t<T_1$) are diffeomorphisms.)

\vspace{0.5truecm}

\noindent
Assume that $t_0<T_1$.  According to Lemma 3.2 in [Ko1], we have 
$$\begin{array}{l}
\displaystyle{N_{(\xi_0,t_0)}
\equiv
\frac{\cos(\sqrt{\varepsilon}k_0br_{t_0}(x_0))}
{\sqrt{\cos^2(\sqrt{\varepsilon}k_0br_{t_0}(x_0))
+\vert\vert({\rm grad}\,r_{t_0})_{x_0}\vert\vert^2}}}\\
\hspace{1.85truecm}\displaystyle{\times
\tau_{\gamma_{w(\xi_0,t_0)}\vert_{[0,1]}}(w^1(\xi_0,t_0))}\\
\hspace{3.6truecm}\displaystyle{({\rm mod}\,\,\,\,
{\rm Span}\{\tau_{\gamma_{\widetilde w(\xi_0,t_0)}\vert_{[0,1]}}
(w^1(\xi_0,t_0))\}^{\perp}).}
\end{array}\leqno{(4.5)}$$
From $(1.1),\,(1.2)$ and $(4.5)$, we obtain the following relation 
$$\begin{array}{l}
\hspace{0.5truecm}\displaystyle{\widetilde f_{\ast}\left(\frac{\partial}{\partial t}\right)_{(\xi_0,t_0)}}\\
\displaystyle{\equiv\frac{\cos(\sqrt{\varepsilon}k_0br_{t_0}(x_0))
(\overline H_{t_0}-\rho_{r_{t_0}}(x_0))}
{\sqrt{\cos^2(\sqrt{\varepsilon}k_0br_{t_0}(x_0))
+\vert\vert({\rm grad}\,r_{t_0})_{x_0}\vert\vert^2}}}\\
\hspace{0.5truecm}\displaystyle{
+\frac{\cos(\sqrt{\varepsilon}k_0br_{t_0}(x_0))}
{\sqrt{\cos^2(\sqrt{\varepsilon}k_0br_{t_0}(x_0))
+\vert\vert({\rm grad}\,r_{t_0})_{x_0}\vert\vert^2}}}\\
\hspace{1truecm}\displaystyle{
\times\left(\frac{(\triangle_Fr_{t_0})(x_0)}
{\cos(\sqrt{\varepsilon}k_0br_{t_0}(x_0))
\sqrt{\cos^2(\sqrt{\varepsilon}k_0br_{t_0}(x_0))+\vert\vert({\rm grad}_{t_0}\,r_{t_0})_{x_0}\vert\vert^2}}
\right.}\\
\hspace{1.5truecm}\displaystyle{\left.
-\frac{(\nabla^Fdr_{t_0})_{x_0}(({\rm grad}_{t_0}r_{t_0})_{x_0},({\rm grad}_{t_0}\,r_{t_0})_{x_0})}
{\cos(\sqrt{\varepsilon}k_0br_{t_0}(x_0))
(\cos^2(\sqrt{\varepsilon}k_0br{t_0}(x_0))+\vert\vert({\rm grad}\,r_{t_0})_{x_0}\vert\vert^2)^{3/2}}\right)}\\
\hspace{1truecm}\displaystyle{
\times\tau_{\gamma_{w(\xi_0,t_0)}\vert_{[0,1]}}(w^1(\xi_0,t_0))}\\
\hspace{4.5truecm}\displaystyle{({\rm mod}\,\,
{\rm Span}\{\tau_{\gamma_{w(\xi_0,t_0)}\vert_{[0,1]}}
(w^1(\xi_0,t_0))\}^{\perp}).}
\end{array}
\leqno{(4.6)}$$
From $(3.5),\,(4.4),\,(4.6)$ and the arbitrariness of $(\xi_0,t_0)$, we can derive the following relation:
$$\begin{array}{l}
\hspace{0.5truecm}\displaystyle{
\frac{\partial\widehat r}{\partial t}(\xi,t)=\frac{\partial r}{\partial t}(c(\xi,t),t)+dr_t
\left(\frac{\partial c}{\partial t}(\xi,t)\right)}\\
\hspace{0truecm}\displaystyle{=
\frac{\cos(\sqrt{\varepsilon}k_0br_t(c(\xi,t)))
(\overline H_t-\rho_{r_t}(c(\xi,t)))}
{\sqrt{\cos^2(\sqrt{\varepsilon}k_0br_t(c(\xi,t)))
+\vert\vert({\rm grad}\,r_t)_{c(\xi,t)}\vert\vert^2}}}\\
\hspace{0.5truecm}\displaystyle{
+\frac{\cos(\sqrt{\varepsilon}k_0br_t(c(\xi,t)))}
{\sqrt{\cos^2(\sqrt{\varepsilon}k_0br_t(c(\xi,t)))
+\vert\vert({\rm grad}\,r_t)_{c(\xi,t)}\vert\vert^2}}}\\
\hspace{1truecm}\displaystyle{
\times\left(\frac{(\triangle_Fr_t)(c(\xi,t))}
{\cos(\sqrt{\varepsilon}k_0br_t(c(\xi,t)))
\sqrt{\cos^2(\sqrt{\varepsilon}k_0br_t(c(\xi,t)))+\vert\vert({\rm grad}\,r_t)_{c(\xi,t)}\vert\vert^2}}\right.}\\
\hspace{1.5truecm}\displaystyle{\left.
-\frac{(\nabla^Fdr_t)_{c(\xi,t)}(({\rm grad}\,r_t)_{c(\xi,t)},({\rm grad}\,r_t)_{c(\xi,t)})}
{\cos(\sqrt{\varepsilon}k_0br_t(c(\xi,t)))
(\cos^2(\sqrt{\varepsilon}k_0br_t(c(\xi,t)))+\vert\vert({\rm grad}\,r_t)_{c(\xi,t)}\vert\vert^2)^{3/2}}\right)}
\end{array}\leqno{(4.7)}$$
$((\xi,t)\in M\times[0,T_1))$.  Next we shall calculate $\displaystyle{\frac{\partial c}{\partial t}}$.  
Denote by $D_t$ the closed domain surrounded by $P$ and $M_t$, 
and $\widetilde D$ the maximal domain in $T^{\perp}B$ containing the $0$-section such that 
$\exp^{\perp}\vert_{\widetilde D}$ is a diffeomorphism into $\overline M$.  
From $c_t(\xi)=(\pi\circ(\exp^{\perp}\vert_{\widetilde D})^{-1})(f_t(\xi))$, we have 
$$
\begin{array}{l}
\displaystyle{
\frac{\partial c}{\partial t}(\xi,t)
=(\pi\circ(\exp^{\perp}\vert_{\widetilde D})^{-1})_{\ast}\left(
\frac{df_t(\xi)}{dt}\right)}\\
\hspace{1.35truecm}\displaystyle{
=(\bar H_t-H_{(\xi,t)})(\pi\circ(\exp^{\perp}\vert_{\widetilde D})^{-1})_{\ast}(N_{(\xi,t)}).}
\end{array}
\leqno{(4.8)}$$
On the other hand, we have 
$$\begin{array}{l}
\displaystyle{N_{(\xi,t)}=\frac{\tau_{\gamma_{w(\xi,t)}}(w^1(\xi,t)-
(D^{co}_{w(\xi,t)})^{-1}({\rm grad}\,r_t)_{c(\xi,t)})}
{\sqrt{1+\vert\vert(D^{co}_{w(\xi,t)})^{-1}
({\rm grad}\,r_t)_{c(\xi,t)}\vert\vert^2}}}\\
\hspace{1truecm}\displaystyle{
=\frac{\tau_{\gamma_{w(\xi,t)}}(\cos(\sqrt{\varepsilon}k_0br_t(c(\xi,t)))w^1(\xi,t)
-({\rm grad}\,r_t)_{c(\xi,t)})}
{\sqrt{\cos^2(\sqrt{\varepsilon}k_0br_t(c(\xi,t)))+\vert\vert({\rm grad}\,r_t)_{c(\xi,t)}\vert\vert^2}}}
\end{array}\leqno{(4.9)}$$
by Lemma 3.2 in [Ko1].  
Let $s\mapsto J(s)\,\,\,(0\leq s<\infty)$ be the Jacobi field along $\gamma_{w^1(\xi,t)}$ with 
$J(0)=({\rm grad}\,r_t)_{c(\xi,t)}$ and $J'(0)=0$.  
Then we have 
$$(\pi\circ(\exp^{\perp}\vert_{\widetilde D})^{-1})_{\ast}(J(s))=J(0)=({\rm grad}\,r_t)_{c(\xi,t)}.$$
On the other hand, according to $(3.1)$, $J(s)$ is described as 
$$\begin{array}{l}
\displaystyle{J(s)=\tau_{\gamma_{w^1(\xi,t)}\vert_{[0,s]}}(D^{co}_{sw^1(\xi,t)}(({\rm grad}\,r_t)_{c(\xi,t)}))}\\
\hspace{0.82truecm}\displaystyle{=\cos(\sqrt{\varepsilon}k_0bs)\tau_{\gamma_{w^1(\xi,t)}\vert_{[0,s]}}
(({\rm grad}\,r_t)_{c(\xi,t)}).}
\end{array}$$
Therefore we obtain 
$$(\pi\circ(\exp^{\perp}\vert_{\widetilde D})^{-1})_{\ast}(\tau_{\gamma_{w^1(\xi,t)}\vert_{[0,s]}}
(({\rm grad}\,r_t)_{c(\xi,t)}))=\frac{1}{\cos(\sqrt{\varepsilon}k_0bs)}({\rm grad}\,r_t)_{c(\xi,t)}.$$
Also, it is clear that 
$(\pi\circ(\exp^{\perp}\vert_{\widetilde D})^{-1})_{\ast}\left(
\tau_{\gamma_{w(\xi,t)}\vert_{[0,s]}}(w^1(\xi,t))\right)=0$.  
Therefore we obtain 
$$\begin{array}{l}
\hspace{0.5truecm}\displaystyle{
(\pi\circ(\exp^{\perp}\vert_{\widetilde D})^{-1})_{\ast}(N_{(\xi,t)})}\\
\displaystyle{
=-\frac{({\rm grad}\,r_t)_{c(\xi,t)}}
{\cos(\sqrt{\varepsilon}k_0br_t(c(\xi,t)))
\sqrt{\cos^2(\sqrt{\varepsilon}k_0br_t(c(\xi,t)))+\vert\vert({\rm grad}\,r_t)_{c(\xi,t)}\vert\vert^2}}.}
\end{array}
\leqno{(4.10)}$$
From $(4.8)$ and $(4.10)$, we can derive 
$$\begin{array}{l}
\hspace{0.5truecm}\displaystyle{\frac{\partial c}{\partial t}(\xi,t)}\\
\displaystyle{=\frac{(H_{(\xi,t)}-\overline H_t)({\rm grad}\,r_t)_{c(\xi,t)}}
{\cos(\sqrt{\varepsilon}k_0br_t(c(\xi,t)))
\sqrt{\cos^2(\sqrt{\varepsilon}k_0br_t(c(\xi,t)))+\vert\vert({\rm grad}\,r_t)_{c(\xi,t)}\vert\vert^2}}}
\end{array}\leqno{(4.11)}$$
and hence 
$$\begin{array}{l}
\hspace{0.5truecm}\displaystyle{dr_t\left(\frac{\partial c}{\partial t}(\xi,t)\right)}\\
\displaystyle{=\frac{(H_{(\xi,t)}-\overline H_t)\vert\vert({\rm grad}\,r_t)_{c(\xi,t)}\vert\vert^2}
{\cos(\sqrt{\varepsilon}k_0br_t(c(\xi,t)))
\sqrt{\cos^2(\sqrt{\varepsilon}k_0br_t(c(\xi,t)))+\vert\vert({\rm grad}\,r_t)_{c(\xi,t)}\vert\vert^2}}.}
\end{array}\leqno{(4.12)}$$
Next we shall calculate the Laplacian $\triangle_t\,\widehat r_t$ of $\widehat r_t$ 
with respect to the metric $g_t$ on $M$ induced by $f_t$.  
Let ${\cal H}^t$ be the horizontal distribution on $t_{r_t}(B)$ and 
$\widehat{\cal H}^t$ be the distribution $M$ with 
$f_{t\ast}(\widehat{\cal H}^t)=\exp^{\perp}_{\ast}({\cal H}^t)$.  
Also, let ${\cal V}^t$ be the vertical distribution on $t_{r_t}(B)$.  
Note that $f_{t\ast}({\cal V}^0)=\exp^{\perp}_{\ast}({\cal V}^t)$.  
For $X\in T_{c_t(\xi)}B$, denote by $\widetilde X^t_{w_t(\xi)}$ the natural lift of $X$ to 
$w_t(\xi)\in t_{r_t}(B)$ and let $\widetilde X^t_{\xi}$ be the element of $T_{\xi}M$ with 
$f_{t\ast}(\widetilde X^t_{\xi})=\exp^{\perp}_{\ast}(\widetilde X^t_{w_t(\xi)})$ (see Figure 8).  
Note that $\{\widetilde X^t_{\xi}\,\vert\,X\in T_{c_t(\xi)}B\}=\widehat{\cal H}^t_{\xi}$.  
Fix $x_0\in B$ and $\xi_0\in M\cap T_{x_0}^{\perp}B$.  
Let $(e_1,\cdots,e_{m^H})$ be an orthonormal tangent frame of $B$ at 
$\bar c_t(x_0)=c_t(\xi_0)$ and 
$\gamma_i$ the geodesic in $B$ with $\gamma_i'(0)=e_i$, where we take $e_1$ as 
$e_1=({\rm grad}\,r_t)_{\bar c_t(x_0)}/\vert\vert({\rm grad}\,r_t)_{\bar c_t(x_0)}\vert\vert$ in the case of 
$({\rm grad}\,r_t)_{\bar c_t(x_0)}\not=0$.  
Since $\tau_{\bar c_t(x_0)}^{-1}(e_1)\in\mathfrak p_{k_0\beta}$, we may assume that 
$\tau_{\bar c_t(x_0)}^{-1}(e_i)\in\mathfrak p_{k_i\beta}$ for some $k_i\in{\cal K}$ ($i=1,\cdots,m^H$).  
Note that $k_1=k_0$.  Also, let ${\widetilde{(\gamma_i)}}^t_{\xi_0}$ be the curve in $M$ starting from $\xi_0$ 
such that $f_t\circ\widetilde{(\gamma_i)}^t_{\xi_0}$ is the natural lift of $\gamma_i$ to $M_t$ 
starting from $f_t(\xi_0)$.  Set 
$$(E_i^t)_s:=\frac{(\widetilde{(\gamma_i)}^t_{\xi_0})'(s)}
{\vert\vert(\widetilde{(\gamma_i)}^t_{\xi_0})'(s)\vert\vert}
\left(
=\frac{(\widetilde{\gamma_i'(s)})^t_{\xi_0}}
{\vert\vert(\widetilde{\gamma_i'(s)})^t_{\xi_0}\vert\vert}
\right)\quad(i=1,\cdots,m^H)$$
(see Figure 9).  
According to $(2.3)$, note that $((E_1^t)_0,\cdots(E_{m^H}^t)_0)$ is an orthonormal base 
of the horizontal subspace $\widehat{\cal H}^t_{\xi_0}$ with respect to $(g_t)_{\xi_0}$ because 
we take $e_1$ as above.  
Also, let $\{E_{m^H+1}^t,\cdots,E_n^t\}$ be a local orthonormal frame field  (with respect to $g_t$) 
of the vertical distribution ${\cal V}$ around $\xi_0$.  
Since $F$ is reflective and $\tau^{-1}_{c(\xi_0,t)}(e_i)\in\mathfrak p_{k_i\beta}$, we we can show that 
${\widetilde{\gamma_i}}^t$ is a pregeodesic in $(M,g_t)$.  Hence we have 
$$\nabla^t_{E_i^t}E_i^t=0\quad\,\,(i=1,\cdots,m^H).\leqno{(4.13)}$$
Also, since $\widehat r_t$ is constant along $\pi^{-1}(\bar c_t^{-1}(x))\cap M$ for each $x\in B$, 
we have 
$$E_j^t\widehat r_t=0\quad\,\,(j=m^H+1,\cdots,n)\leqno{(4.14)}$$

\vspace{0.4truecm}

\centerline{
\unitlength 0.1in
%
\hspace{4.5truecm}}

\vspace{0.4truecm}


\centerline{{\bf Figure 9.}}

\vspace{0.5truecm}

\noindent
From $(4.13)$ and $(4.14)$, we have 
$$
(\triangle_t\widehat r_t)_{\xi_0}=\sum_{i=1}^n
\nabla^t_{(E^t_i)_0}\nabla^t_{(E^t_i)_0}\widehat r_t
=\sum_{i=1}^{m^H}(E^t_i)_0(E_i^t\widehat r_t).\leqno{(4.15)}$$
Denote by $\triangle_F\,r_t$ the Laplacian of $r_t$ with respect to $g_F$.  Then we have 
$(\triangle_Fr_t)(c(\xi_0,t))=\sum\limits_{i=1}^{m^H}e_i(\gamma_i'r_t)$.  
From the definitions of $r_t$ and $\widehat r_t$, we have 
$$
e_i(\gamma_i'r_t)=\left.\frac{d^2}{ds^2}\right\vert_{s=0}r_t(\gamma_i(s))
=\left.\frac{d^2}{ds^2}\right\vert_{s=0}\widehat r_t((\widetilde{\gamma_i})^t_{\xi_0}(s))
=\widetilde{(e_i)}^t_{\xi_0}\left(((\widetilde{\gamma_i})^t_{\xi_0})'\widehat r_t\right).
$$
On the other hand, from $(2.3)$, we have 
$$\begin{array}{l}
\hspace{0.5truecm}\displaystyle{
\widetilde{(e_i)}^t_{\xi_0}\left(((\widetilde{\gamma_i})^t_{\xi_0})'\widehat r_t\right)
=\left(\cos^2(\sqrt{\varepsilon}k_0br_t(x_0))
+\vert\vert({\rm grad}\,r_t)_{x_0}\vert\vert^2\right)(E^t_i)_0(E^t_i\widehat r_t)}\\
\hspace{3.95truecm}\displaystyle{-\frac{\sqrt{\varepsilon}k_0b\sin(\sqrt{\varepsilon}k_0br_t(x_0))
\cos(\sqrt{\varepsilon}k_0br_t(x_0))(e_ir_t)^2}
{\sqrt{\cos^2(\sqrt{\varepsilon}k_0br_t(x_0))+\vert\vert({\rm grad}\,r_t)_{x_0}\vert\vert^2}}}\\
\hspace{3.95truecm}\displaystyle{+\frac{(\nabla^Fdr_t)(e_i,({\rm grad}\,r_t)_{x_0})(e_ir_t)}
{\sqrt{\cos^2(\sqrt{\varepsilon}k_0br_t(x_0))+\vert\vert({\rm grad}\,r_t)_{x_0}\vert\vert^2}}.}
\end{array}$$
Hence we can derive 
$$\begin{array}{l}
\hspace{0.5truecm}\displaystyle{(\triangle_Fr_t)(x_0)
=\left(\cos^2(\sqrt{\varepsilon}k_0br_t(x_0))+\vert\vert({\rm grad}\,r_t)_{x_0}\vert\vert^2\right)
(\triangle_t\widehat r_t)(\xi_0)}\\
\hspace{2.3truecm}\displaystyle{-\frac{\sqrt{\varepsilon}k_0b\sin(\sqrt{\varepsilon}k_0br_t(x_0))
\cos(\sqrt{\varepsilon}k_0br_t(x_0))\vert\vert({\rm grad}\,r_t)_{x_0}\vert\vert^2}
{\cos^2(\sqrt{\varepsilon}k_0br_t(x_0))+\vert\vert({\rm grad}\,r_t)_{x_0}\vert\vert^2}}\\
\hspace{2.3truecm}\displaystyle{+\frac{(\nabla^Fdr_t)_{x_0}(({\rm grad}\,r_t)_{x_0},({\rm grad}\,r_t)_{x_0})}
{\cos^2(\sqrt{\varepsilon}k_0br_t(x_0))+\vert\vert({\rm grad}\,r_t)_{x_0}\vert\vert^2}.}
\end{array}\leqno{(4.16)}$$
From $(3.5),(4.7),(4.12),(4.16)$ and the definition of $\rho_{r_t}$, 
the following evolution equations are derived.  

\vspace{0.5truecm}

\noindent
{\bf Lemma 4.1.} {\sl {\rm (i)} The radius functions $r_t$'s satisfies the following evolution 
equation:
$$\begin{array}{l}
\hspace{0.5truecm}\displaystyle{\frac{\partial r}{\partial t}(x,t)
-\frac{(\triangle_Fr_t)(x)}{\cos^2(\sqrt{\varepsilon}k_0br_t(x))}}\\
\displaystyle{=\frac{\sqrt{\cos^2(\sqrt{\varepsilon}k_0br_t(x))+\vert\vert({\rm grad}\,r_t)_x\vert\vert^2}}
{\cos(\sqrt{\varepsilon}k_0br_t(x))}\cdot(\overline H_t-\rho_{r_t}(x))}\\
\hspace{0.5truecm}\displaystyle{-\frac{(\nabla^Fdr_t)_x(({\rm grad}\,r_t)_x,({\rm grad}\,r_t)_x)}
{\cos^2(\sqrt{\varepsilon}k_0br_t(x))(\cos^2(\sqrt{\varepsilon}k_0br_t(x))
+\vert\vert({\rm grad}\,r_t)_x\vert\vert^2)}}
\end{array}\leqno{(4.17)}$$
$((x,t)\in M\times[0,T_1))$.  

{\rm(ii)} The radius functions $\widehat r_t$'s satisfies the following evolution 
equation:
$$\begin{array}{l}
\hspace{0.5truecm}\displaystyle{\frac{\partial\widehat r}{\partial t}(\xi,t)
-(\triangle_t\widehat r_t)(\xi)}\\
\displaystyle{=\frac{\cos(\sqrt{\varepsilon}k_0b\widehat r_t(\xi))
(\overline H_t-\rho_{r_t}(c(\xi,t)))}
{\sqrt{\cos^2(\sqrt{\varepsilon}k_0b\widehat r_t(\xi))+\vert\vert({\rm grad}\,r_t)_{c(\xi,t)}\vert\vert^2}}}\\
\hspace{0.5truecm}\displaystyle{-\frac{\sqrt{\varepsilon}k_0b\sin(\sqrt{\varepsilon}k_0b\widehat r_t(\xi))
\cos(\sqrt{\varepsilon}k_0b\widehat r_t(\xi))\vert\vert({\rm grad}\,r_t)_{c(\xi,t)}\vert\vert^2}
{(\cos^2(\sqrt{\varepsilon}k_0b\widehat r_t(\xi))+\vert\vert({\rm grad}\,r_t)_{c(\xi,t)}\vert\vert^2)^2}}
\end{array}\leqno{(4.18)}$$
$((\xi,t)\in M\times[0,T_1))$.
}

\vspace{0.5truecm}

Replacing $\overline H$ in $(4.17)$ to any $C^{1,\alpha/2}$ real-valued function $\phi$ such that 
$\phi(0)=\overline H(0)$, we obtain a parabolic equation, which has a unique solution $r_t$ such that 
${\rm grad}\,r_t=0$ along $\partial B$ in short time for any initial data $r_0$ such that ${\rm grad}\,r_0=0$ 
holds along $\partial B$.  By using a routin fixed point argument (see [M]), 
we can establish the short time existence and uniqueness also for $(4.17)$ with the same boundary condition.  
From this fact, we can derive the following statement.  

\vspace{0.5truecm}

\noindent
{\bf Proposition 4.2.} {\sl Under Setting (S), 
assume that $\overline M$ is a rank one symmetric space 
and that $F$ is an invariant submanifold.  Then there uniquely exists the volume-preserving mean curvature flow 
$f_t:M\hookrightarrow\overline M$ starting from $f$ and satisfying the condition (C1) in short time.}

\vspace{0.5truecm}

\noindent
{\it Proof.} Since $\overline M$ is of rank one and $F$ is invariant, the condition (C2) holds automatically.  
Hence the radius functions $r_t$'s of $f_t$ satisfy $(4.17)$.  Hence the statement is derived from the above 
fact.  \hspace{7.65truecm}q.e.d.

\vspace{0.5truecm}

Denote by ${\rm grad}_t\widehat r_t$ the gradient vector field of $\widehat r_t$ with respect to $g_t$.  
The following relation holds between ${\rm grad}\,r_t$ and ${\rm grad}_t\widehat r_t$.  

\vspace{0.5truecm}

\noindent
{\bf Lemma 4.3.} {\sl The norm of the gradient vector $({\rm grad}_t\widehat r_t)_{\xi}$ is described as 
$$\vert\vert({\rm grad}_t\widehat r_t)_{\xi}\vert\vert_t
=\frac{\vert\vert({\rm grad}\,r_t)_{c(\xi,t)}\vert\vert}
{\sqrt{\cos^2(\sqrt{\varepsilon}k_0br_t(c(\xi,t)))+
\vert\vert({\rm grad}\,r_t)_{c(\xi,t)}\vert\vert^2}},
\leqno{(4.19)}$$
where $\vert\vert\cdot\vert\vert_t$ is the norm of $(\cdot)$ with respect to $g_t$.}

\vspace{0.5truecm}

\noindent
{\it Proof.} 
It is clear that 
$dr_t(({\rm grad}\,r_t)_{c(\xi,t)})
=d\widehat r_t((\widetilde{({\rm grad}\,r_t)_{c(\xi,t)}})^t_{\xi})$, and that 
$(\widetilde{({\rm grad}\,r_t)_{c(\xi,t)}})^t_{\xi}$ and $({\rm grad}_t\widehat r_t)_{\xi}$ are linearly 
dependent.  Hence we have 
$$\vert\vert({\rm grad}\,r_t)_{c(\xi,t)}\vert\vert^2
=\vert\vert(\widetilde{({\rm grad}\,r_t)_{c(\xi,t)}})^t_{\xi}\vert\vert_t\cdot 
\vert\vert({\rm grad}_t\widehat r_t)_{\xi}\vert\vert_t.$$
On the other hand, it follows from $(2.3)$ that 
$$\vert\vert(\widetilde{({\rm grad}\,r_t)_{c(\xi,t)}})^t_{\xi}\vert\vert_t
=\vert\vert({\rm grad}\,r_t)_{c(\xi,t)}\vert\vert\cdot
\sqrt{\cos^2(\sqrt{\varepsilon}k_0br_t(c(\xi,t)))+
\vert\vert({\rm grad}\,r_t)_{c(\xi,t)}\vert\vert^2}.$$
Therefore, we obtain the desired relation.  
\hspace{6truecm}q.e.d.

\section{The evolution of the gradient of the radius function} 
We use the notations in Introduction and Sections 1-4.  
Let $T_1$ be as in Section 4.  
Define a function $\widehat u_t:M\to{\Bbb R}$ ($t\in[0,T_1)$) by 
$$\widehat u_t(\xi):=\bar g(N_{(\xi,t)},\tau_{\gamma_{w(\xi,t)}\vert_{[0,1]}}
(w^1(\xi,t)))\,\,\,\,(\xi\in M)$$ 
and a map 
$\widehat v_t:M\to{\Bbb R}$ by $\widehat v_t:=\frac{1}{\widehat u_t}$ ($0\leq t<T_1$).  
Define a map $\widehat u:M\times[0,T_1)\to{\Bbb R}$ by 
$\widehat u(\xi,t):=\widehat u_t(\xi)\,\,\,\,((\xi,t)\in M\times[0,T_1))$ and a map 
$\widehat v:M\times[0,T_1)\to{\Bbb R}$ by 
$\widehat v(\xi,t):=\widehat v_t(\xi)\,\,\,((\xi,t)\in M\times[0,T_1))$.  
Define a function $\bar u_t$ (resp. $\bar v_t$) over $B$ by 
$\bar u_t(x):=\widehat u_t(\xi)\,\,(x\in B)$ (resp. $\bar v_t(x):=\widehat v_t(\xi)\,\,(x\in B)$, 
where $\xi$ is an arbitrary element of $M\cap S^{\perp}_xB$.  
It is clear that these functions are well-defined.  
Set $u_t:=\bar u_t\circ{\bar c}_t^{-1}$ and $v_t:=\bar v_t\circ{\bar c}_t^{-1}$.  
We have only to show 
$\inf_{(x,t)\in B\times[0,T_1)}u(x,t)>0$, that is, 
$\sup_{(x,t)\in B\times[0,T_1)}v(x,t)<\infty$.  
In the sequel, assume that $t<T_1$.  
From $(4.9)$, we have 
$$\widehat u_t(\xi)=\frac{\cos(\sqrt{\varepsilon}k_0b\widehat r(\xi,t))}
{\sqrt{\cos^2(\sqrt{\varepsilon}k_0b\widehat r(\xi,t))+\vert\vert({\rm grad}\,r_t)_{c(\xi,t)}\vert\vert^2}}.
\leqno{(5.1)}$$
In order to investigate te evolution of the gradient ${\rm grad}\,r_t$ of the radius function $r_t$, we suffice 
to investigate the evolution of $\widehat u_t$.  
It is easy to show that the outward unit normal vector field $N_t$ satisfies the following 
evolution equation:
$$\frac{\partial N}{\partial t}=f_{t\ast}({\rm grad}\,H_t).\leqno{(5.2)}$$
For simplicity, we set 
${\widehat w}^1(\xi,t):=\tau_{\gamma_{w(\xi,t)}\vert_{[0,1]}}(w^1(\xi,t))$.  
We calculate $\frac{\partial\widehat u}{\partial t}$.  
Define a map $\delta:[0,T_1)\times{\Bbb R}\to\overline M$ by 
$$\delta(t,s):=\exp_{c(\xi,t)}(sw^1(\xi,t)).$$
Then we have 
$\displaystyle{\left.\frac{\partial\delta}{\partial s}\right\vert_{s=\widehat r(\xi,t)}
=\widehat w^1(\xi,t)}$.  
For a fixed $t\in[0,T_1)$, 
$\displaystyle{\widehat Y_t:\,s\,\mapsto\left(\frac{\partial\delta}{\partial t}\right)(t,s)}$ is 
the Jacobi field along the geodesic $\gamma_{w^1(\xi,t)}$.  
Since $\displaystyle{\widehat Y_t(0)=\left(\frac{\partial c}{\partial t}\right)(\xi,t)}$ and 
$$\widehat Y'_t(0)=\left.\left(\overline{\nabla}^{\delta}_{\frac{\partial}{\partial s}}
\frac{\partial\delta}{\partial t}\right)\right\vert_{s=0}
=\left.\left(\overline{\nabla}^{\delta}_{\frac{\partial}{\partial t}}
\frac{\partial\delta}{\partial s}\right)\right\vert_{s=0}
=\overline{\nabla}^{\delta(\cdot,0)}_{\frac{\partial}{\partial t}}w^1(\xi,t)
=\nabla^{\perp_B}_{\widehat Y_t(0)}w^1(\cdot,t)$$
by the reflectivity of $F$, it follows from $(3.1)$ that 
$$\widehat Y_t(s)=\tau_{\gamma_{w^1(\xi,t)}\vert_{[0,s]}}\left(D^{co}_{sw^1(\xi,t)}
\left(\left(\frac{\partial c}{\partial t}\right)(\xi,t)\right)
+sD^{si}_{sw^1(\xi,t)}(\nabla^{\perp_B}_{\widehat Y_t(0)}w^1(\cdot,t))\right),$$
where $\nabla^{\perp_B}$ is the normal connection of $B$.  
This implies together with $(4.11)$ that 
$$\begin{array}{l}
\displaystyle{\widehat Y_t(s)=-\frac{(\overline H_t-H_{(\xi,t)})\cos(\sqrt{\varepsilon}k_0bs)\widehat u(\xi,t)}
{\cos^2(\sqrt{\varepsilon}k_0b\widehat r(\xi,t))}}\\
\hspace{1.9truecm}\displaystyle{\times\tau_{\gamma_{w^1(\xi,t)}\vert_{[0,s]}}(({\rm grad}_tr_t)_{c(\xi,t)})}\\
\hspace{1.4truecm}\displaystyle{
+s\tau_{\gamma_{w^1(\xi,t)}\vert_{[0,s]}}(D^{si}_{sw^1(\xi,t)}(\nabla^{\perp_B}_{\widehat Y_t(0)}w^1(\cdot,t))).}
\end{array}\leqno{(5.3)}$$
In particular, we have 
$$\begin{array}{l}
\displaystyle{\widehat Y_t(\widehat r(\xi,t))
=-\frac{(\overline H_t-H_{(\xi,t)})\widehat u(\xi,t)}{\cos(\sqrt{\varepsilon}k_0b\widehat r(\xi,t))}
\tau_{\gamma_{w(\xi,t)}}(({\rm grad}_tr_t)_{c(\xi,t)})}\\
\hspace{2.1truecm}\displaystyle{+\widehat r(\xi,t)\tau_{\gamma_{w(\xi,t)}\vert_{[0,1]}}
(D^{si}_{w(\xi,t)}(\nabla^{\perp_B}_{\widehat Y_t(0)}w^1(\cdot,t))).}
\end{array}\leqno{(5.4)}$$
Also, we have 
$$\begin{array}{l}
\displaystyle{\left(\frac{\partial\widehat w^1}{\partial t}\right)(\xi,t)
=\left.\left(\overline{\nabla}^{\delta}_{\frac{\partial}{\partial t}}
\frac{\partial\delta}{\partial s}\right)\right\vert_{s=\widehat r(\xi,t)}
=\overline{\nabla}^{\delta}_{\frac{\partial}{\partial s}\vert_{s=\widehat r(\xi,t)}}\widehat Y_t}\\
\hspace{0truecm}\displaystyle{=-\frac{(\overline H_t-H_{(\xi,t)})\sqrt{\varepsilon}k_0b
\sin(\sqrt{\varepsilon}k_0b\widehat r(\xi,t))\widehat u(\xi,t)}
{\cos^2(\sqrt{\varepsilon}k_0b\widehat r(\xi,t))}}\\
\hspace{0.5truecm}\displaystyle{\times\tau_{\gamma_{w(\xi,t)}}(({\rm grad}_tr_t)_{c(\xi,t)})}\\
\hspace{0.5truecm}\displaystyle{+\tau_{\gamma_{w(\xi,t)}\vert_{[0,1]}}
(D^{co}_{w(\xi,t)}(\nabla^{\perp_B}_{\widehat Y_t(0)}w^1(\cdot,t))).}
\end{array}\leqno{(5.5)}$$
From this relation and $(4.9)$, we obtain 
$$\begin{array}{l}
\displaystyle{\bar g\left(N(\xi,t),\left(\frac{\partial{\widehat w}^1}{\partial t}\right)(\xi,t)\right)}\\
\displaystyle{=(\overline H_t-H_{(\xi,t)})
\sqrt{\varepsilon}k_0b\tan(\sqrt{\varepsilon}k_0b\widehat r(\xi,t))(1-\widehat u(\xi,t)^2).}
\end{array}\leqno{(5.6)}$$
From $(5.2)$ and $(5.6)$, we obtain 
$$\begin{array}{l}
\displaystyle{\frac{\partial\widehat u}{\partial t}(\xi,t)
=\bar g(f_{t\ast}(({\rm grad}\,H_t)_{\xi}),{\widehat w}^1(\xi,t))}\\
\hspace{1.9truecm}\displaystyle{
+(\overline H_t-H_{(\xi,t)})
\sqrt{\varepsilon}k_0b\tan(\sqrt{\varepsilon}k_0b\widehat r(\xi,t))(1-\widehat u(\xi,t)^2)}
\end{array}
\leqno{(5.7)}$$
($(\xi,t)\in M\times[0,T_1)$).  
Fix $\xi_0\in M$ and $t_0\in[0,T_1)$.  Set $x_0:=c(\xi_0,t_0)$.  
Let $(e_1,\cdots,e_{m^H}),\gamma_i,\,{\widetilde{(\gamma_i)}}^t_{\xi_0},\,(E_i^t)_s,
\,\,(i=1,\cdots,m^H)$ and $\{E_{m^H+1}^t,\cdots,E_n^t\}$ be as in Section 4.  
We may assume that 
$$(\tau_{\gamma_{w(\xi,t)}\vert_{[0,1]}}\circ\tau_{c(\xi,t)})^{-1}
(f_{t\ast}((E^t_j)_{\xi}))\in\mathfrak p_{k_j\beta}$$
holds for some $k_j\in\{1,2\}$ ($j=m^H+1,\cdots,n$).  
Define a function $\lambda_t$ ($t\in[0,T_1)$) over $M$ by 
$$\begin{array}{l}
\displaystyle{\lambda_t(\xi):=-\widehat u(\xi,t)
\left(\sum_{k\in{\cal K}}m_k^H\sqrt{\varepsilon}kb\tan(\sqrt{\varepsilon}kb\widehat r(\xi,t))
+\frac{(\triangle_Fr_t)(x)}{\cos^2(\sqrt{\varepsilon}k_0b\widehat r(\xi,t))}\right.}\\
\hspace{1.7truecm}\displaystyle{
+\frac{\vert\vert({\rm grad}\,r_t)_x\vert\vert^2\sqrt{\varepsilon}k_0b
\tan(\sqrt{\varepsilon}k_0b\widehat r(\xi,t))}
{\cos^2(\sqrt{\varepsilon}k_0b\widehat r(\xi,t))+\vert\vert({\rm grad}\,r_t)_x\vert\vert^2}}\\
\hspace{1.7truecm}\displaystyle{\left.
-\frac{(\nabla^Fdr_t)_x(({\rm grad}\,r_t)_x,({\rm grad}\,r_t)_x)}
{\cos^2(\sqrt{\varepsilon}k_0b\widehat r(\xi,t))
(\cos^2(\sqrt{\varepsilon}k_0b\widehat r(\xi,t))+\vert\vert({\rm grad}\,r_t)_{x}\vert\vert^2)}\right)}
\end{array}$$
($\xi\in M$), where $x=c(\xi,t)$.  
From $\bar g(f_{t_0\ast}(\widetilde{(e_i)}^{t_0}_{\xi_0}),\widehat w^1(\xi_0,t_0))=0\,\,(i=2,\cdots,m^H)$ and 
$(2.50)$, we have 
$$\begin{array}{l}
\hspace{0.5truecm}\displaystyle{\bar g(f_{t_0\ast}(({\rm grad}\,H_{t_0})_{\xi_0}),{\widehat w}^1(\xi_0,t_0))}\\
\hspace{0truecm}\displaystyle{
=\frac{(\widetilde{(e_1)}^{t_0}_{\xi_0}H_{t_0})\bar g(f_{t_0\ast}(\widetilde{(e_1)}^{t_0}_{\xi_0}),
\widehat w^1(\xi_0,t_0))}{\vert\vert\widetilde{(e_1)}^{t_0}_{\xi_0}\vert\vert^2}}\\
\hspace{0truecm}\displaystyle{=\frac{\sqrt{1-\widehat u(\xi_0,t_0)^2}}
{\vert\vert\widetilde{(e_1)}^{t_0}_{\xi_0}\vert\vert}\cdot\widetilde{(e_1)}^{t_0}_{\xi_0}H_{t_0}}\\
\hspace{0truecm}\displaystyle{
=\frac{\widehat u(\xi_0,t_0)\sqrt{1-\widehat u(\xi_0,t_0)^2}}{\cos(\sqrt{\varepsilon}k_0b\widehat r(\xi_0,t_0))}
\cdot\widetilde{(e_1)}^{t_0}_{\xi_0}\lambda_{t_0}}\\
\hspace{0.5truecm}\displaystyle{+(1-\widehat u(\xi_0,t_0)^2)\left(\lambda_{t_0}(\xi_0)
+\widehat u(\xi_0,t_0)\sum_{k\in{\cal K}}m_k^H\sqrt{\varepsilon}kb
\tan(\sqrt{\varepsilon}kb\widehat r(\xi_0,t_0))\right)}\\
\hspace{1truecm}\displaystyle{\times
\sum_{k=1}^2\frac{m_k^V\sqrt{\varepsilon}kb}{\tan(\sqrt{\varepsilon}kb\widehat r(\xi_0,t_0))}}\\
\hspace{0.5truecm}\displaystyle{-\widehat u(\xi_0,t_0)(1-\widehat u(\xi_0,t_0)^2)
\sum_{k=1}^2\frac{m_k^V(\sqrt{\varepsilon}kb)^2}{\sin^2(\sqrt{\varepsilon}kb\widehat r(\xi_0,t_0))}.}
\end{array}\leqno{(5.8)}$$
From $(2.50),\,(5.7)$ and $(5.8)$, we obtain 
$$\begin{array}{l}
\hspace{0truecm}\displaystyle{\frac{\partial\widehat u}{\partial t}(\xi,t)
=\frac{\widehat u(\xi,t)\sqrt{1-\widehat u(\xi,t)^2}}{\cos(\sqrt{\varepsilon}k_0b\widehat r(\xi,t))}
\cdot\widetilde{(e_1)}^t_{\xi}\lambda_t}\\
\hspace{1.5truecm}\displaystyle{+(1-\widehat u(\xi,t)^2)\left(\lambda_t(\xi)
+\widehat u(\xi,t)\sum_{k\in{\cal K}}m_k^H\sqrt{\varepsilon}kb
\tan(\sqrt{\varepsilon}kb\widehat r(\xi,t))\right)}\\
\hspace{2truecm}\displaystyle{\times
\sum_{k=1}^2\frac{m_k^V\sqrt{\varepsilon}kb}{\tan(\sqrt{\varepsilon}kb\widehat r(\xi,t))}}\\
\hspace{1.5truecm}\displaystyle{-\widehat u(\xi,t)(1-\widehat u(\xi,t)^2)
\sum_{k=1}^2\frac{m_k^V(\sqrt{\varepsilon}kb)^2}{\sin^2(\sqrt{\varepsilon}kb\widehat r(\xi,t))}}\\
\hspace{1.5truecm}\displaystyle{+\overline H_t\sqrt{\varepsilon}k_0b\tan(\sqrt{\varepsilon}k_0b\widehat r(\xi,t))
(1-\widehat u(\xi,t)^2)}\\
\hspace{1.5truecm}\displaystyle{+(1-\widehat u(\xi,t)^2)\sqrt{\varepsilon}k_0b
\tan(\sqrt{\varepsilon}k_0b\widehat r(\xi,t))}\\
\hspace{2truecm}\displaystyle{\times\left(\lambda_t(\xi)-\widehat u(\xi,t)\sum_{k=1}^2
\frac{m_k^V\sqrt{\varepsilon}kb}
{\tan(\sqrt{\varepsilon}kb\widehat r(\xi,t))}\right).}
\end{array}\leqno{(5.9)}$$
Define a map $\delta_i:[0,\infty)\times(-\varepsilon,\varepsilon)\to \overline M$ 
($i=1,\cdots,m^H$) by 
$$\delta_i(\bar t,s):=\exp_{\gamma_i(s)}(\bar tw(\widetilde{(\gamma_i)}^{t_0}_{\xi_0}(s),t_0)),$$
where $\varepsilon$ is a small positive number.  
Set 
$\displaystyle{Y_i^{s_0}:=\left.\frac{\partial\delta_i}{\partial s}
\right\vert_{s=s_0}}$.  
Since $Y_i^{s_0}$ is the Jacobi field along $\gamma_{w(\widetilde{(\gamma_i)}^{t_0}_{\xi_0}(s_0),t_0)}$ with 
$Y_i^{s_0}(0)=\gamma'_i(s_0)$ and 
$({Y_i}^{s_0})'(0)=(\gamma'_i(s_0)r_{t_0})w^1(\widetilde{(\gamma_i)}^{t_0}_{\xi_0}(s_0),t_0)$, 
it is described as 
$$\begin{array}{l}
\hspace{0.5truecm}\displaystyle{Y_i^{s_0}(\bar t)}\\
\displaystyle{=\tau_{\gamma_{w(\widetilde{(\gamma_i)}^{t_0}_{\xi_0}(s_0),t_0)}\vert_{[0,\bar t]}}
\left(D^{co}_{\bar tw(\widetilde{(\gamma_i)}^{t_0}_{\xi_0}(s_0),t_0)}(\gamma_i'(s_0))\right.}\\
\hspace{4truecm}\displaystyle{\left.+\bar tD^{si}_{\bar tw(\widetilde{(\gamma_i)}^{t_0}_{\xi_0}(s_0),t_0)}
((\gamma'_i(s_0)r_{t_0})w^1(\widetilde{(\gamma_i)}^{t_0}_{\xi_0}(s_0),t_0))\right)}\\
\hspace{0truecm}\displaystyle{
=\tau_{\gamma_{w(\widetilde{(\gamma_i)}^{t_0}_{\xi_0}(s_0),t_0)}\vert_{[0,\bar t]}}
\left(D^{co}_{\bar tw(\widetilde{(\gamma_i)}^{t_0}_{\xi_0}(s_0),t_0)}(\gamma_i'(s_0))\right)
+\bar t(\gamma'_i(s_0)r_{t_0})\widehat w^1(\widetilde{(\gamma_i)}^{t_0}_{\xi_0}(s_0),t_0).}
\end{array}$$
Hence we have 
$$\begin{array}{l}
\hspace{0truecm}\displaystyle{\overline{\nabla}^{f_{t_0}}_{(E^{t_0}_i)_{s_0}}{\widehat w}^1_{t_0}}\\
\hspace{0truecm}\displaystyle{
=\frac{(Y_i^{s_0})'(1)}{\sqrt{\cos^2(\sqrt{\varepsilon}k_ib\widehat r({\widetilde{(\gamma_i)}}^{t_0}_{\xi_0}(s_0),
t_0))
+\vert\vert({\rm grad}\,r_{t_0})_{c({\widetilde{(\gamma_i)}}^{t_0}_{\xi_0}(s_0),t_0)}\vert\vert^2\delta_{i1}}}}\\
\displaystyle{=\frac{1}{\sqrt{\cos^2(\sqrt{\varepsilon}k_ib\widehat r(\widetilde{(\gamma_i)}^{t_0}_{\xi_0}(s_0),
t_0))
+\vert\vert({\rm grad}\,r_{t_0})_{c(\widetilde{(\gamma_i)}^{t_0}_{\xi_0}(s_0),t_0)}\vert\vert^2\delta_{i1}}}}\\
\hspace{0truecm}\displaystyle{\times\left\{
-\sqrt{\varepsilon}k_ib\widehat r(\widetilde{(\gamma_i)}^{t_0}_{\xi_0}(s_0),t_0)
\sin(\sqrt{\varepsilon}k_ib\widehat r(\widetilde{(\gamma_i)}^{t_0}_{\xi_0}(s_0),t_0))
\tau_{\gamma_{w(\widetilde{(\gamma_i)}^{t_0}_{\xi_0}(s_0),t_0)}\vert_{[0,1]}}(e_i)\right.}\\
\hspace{0.5truecm}\displaystyle{\left.+(\gamma'_i(s_0)r_{t_0})
\widehat w^1(\widetilde{(\gamma_i)}^{t_0}_{\xi_0}(s_0),t_0)\right\},}
\end{array}\leqno{(5.10)}$$
%
where $\overline{\nabla}^{f_{t_0}}$ is the pull-back connection of $\overline{\nabla}$ by $f_{t_0}$.  
Let $\alpha_j:(-\varepsilon,\varepsilon)\to M$ ($j=m^H+1,\cdots,n$) be 
the integral curve of $E_j^{t_0}$ with $\alpha_j(0)=\xi_0$, where 
$\varepsilon$ is a small positive number.  
Define a map $\widehat{\delta}_j:[0,\infty)\times(-\varepsilon,\varepsilon)\to\overline M$ 
($j=m^H+1,\cdots,n$) by $\widehat{\delta}_j(\bar t,s)
:=\exp_{x_0}(\bar t(f_{t_0}\circ\alpha_j)(s))$.  
Set 
$\displaystyle{\widehat Y_j^{s_0}:=\left.\frac{\partial\widehat\delta_j}
{\partial s}\right\vert_{s=s_0}}$.  
Since $\widehat Y_j^{s_0}$ is the Jacobi field along 
$\gamma_{w(\alpha_j(s_0),t_0)}$ with 
${\widehat Y}_j^{s_0}(0)=0$ and $({\widehat Y}_j^{s_0})'(0)=(\exp_{x_0})_{\ast}^{-1}
(f_{t_0\ast}((E_j^{t_0})_{\alpha_j(s_0)}))$, it is described as 
$$
{\widehat Y}^{s_0}_j(\bar t)=\tau_{\gamma_{w(\alpha_j(s_0),t_0)}\vert_{[0,\bar t]}}
\left(\bar tD^{si}_{\bar tw(\alpha_j(s_0),t_0)}((\exp_{x_0})_{\ast}^{-1}
(f_{t_0\ast}((E_j^{t_0})_{\alpha_j(s_0)})))\right).$$
For simplicity, we denote $(E_j^{t_0})_{\alpha_j(s)}$ by $(E_j^{t_0})_s$.  
Hence we have 
$$\begin{array}{l}
\hspace{0.5truecm}\displaystyle{\overline{\nabla}^{f_{t_0}}_{(E^{t_0}_j)_{s_0}}
{\widehat w}^1_{t_0}=({\widehat Y}_j^{s_0})'(1)}\\
\displaystyle{=\tau_{\gamma_{w(\alpha_j(s_0),t_0)}\vert_{[0,1]}}
\left(D^{co}_{w(\alpha_j(s_0),t_0)}((\exp_{x_0})_{\ast}^{-1}
(f_{t_0\ast}((E_j^{t_0})_{s_0}))\right)}\\
\displaystyle{=\frac{\sqrt{\varepsilon}k_jb\widehat r(\alpha_j(s_0),t_0)}
{\tan(\sqrt{\varepsilon}k_jb\widehat r(\alpha_j(s_0),t_0))}f_{t_0\ast}((E_j^{t_0})_{s_0}).}
\end{array}\leqno{(5.11)}$$
%
Denote by $\triangle_t\widehat u_t$ the rough Laplacian of $\widehat u_t$ with respect to the Riemannian 
connection $\nabla^t$ of $g_t$.  
According to $(2.28)$, we have 
$$\bar g(N(\xi,t),\tau_{\gamma_{w(\xi,t)}\vert_{[0,1]}}(e_1))=-\sqrt{1-\widehat u(\xi,t)^2},$$
that is, 
$$\widehat u(\xi,t)=\sqrt{1-\bar g(N(\xi,t),\tau_{\gamma_{w(\xi,t)}\vert_{[0,1]}}(e_1))^2}.\leqno{(5.12)}$$
By using this relation, we can derive 
$$\begin{array}{l}
\displaystyle{(E_1^{t_0})_s\widehat u_{t_0}
=\frac{\vert\vert({\rm grad}\,r_{t_0})_{c(\widetilde{(\gamma_i)}^{t_0}_{\xi_0}(s),t_0)}\vert\vert}
{\cos(\sqrt{\varepsilon}k_0b\widehat r({\widetilde{(\gamma_1)}}^{t_0}_{\xi_0}(s),t_0))
\vert\vert\widetilde{(\gamma'_1(s))}^{t_0}_{\widetilde{(\gamma_1)}^{t_0}_{\xi_0}(s)}\vert\vert}}\\
\hspace{0.5truecm}\displaystyle{\times\left\{
\bar g\left(f_{t_0\ast}(A^{t_0}(\widetilde{(\gamma'_1(s))}^{t_0}_{\widetilde{(\gamma_1)}^{t_0}_{\xi_0}(s)}),
\tau_{\gamma_{w(\widetilde{(\gamma_1)}^{t_0}_{\xi_0}(s),t_0)}\vert_{[0,1]}}(\gamma'_1(s))\right)\right.}\\
\hspace{1truecm}\displaystyle{\left.
+\bar g\left(N(\widetilde{(\gamma_1)}^{t_0}_{\xi_0}(s),t_0),
\frac{\overline{\nabla}}{ds}\left(\tau_{\gamma_{w(\widetilde{(\gamma_1)}^{t_0}_{\xi_0}(s),t_0)}\vert_{[0,1]}}
(\gamma'_1(s))\right)\right)\right\},}
\end{array}\leqno{(5.13)}$$
where $\displaystyle{\frac{\overline{\nabla}}{ds}}$ denotes the covariant derivative along the curve 
$s\mapsto \exp^{\perp}\circ w(\widetilde{(\gamma_1)}^{t_0}_{\xi_0}(s),t_0)
=f_{t_0}\circ\widetilde{(\gamma_1)}^{t_0}_{\xi_0}(s)$ 
with respect to $\overline{\nabla}$.  
Define a curve $\bar w_{t_0,s}$ in $T^{\perp}B$ by 
$$\bar w_{t_0,s}(\hat s)
:=\frac{\widehat r({\widetilde{(\gamma_1)}}^{t_0}_{\xi_0}(s))}
{\widehat r({\widetilde{(\gamma_1)}}^{t_0}_{\xi_0}(\hat s))}
w(\widetilde{(\gamma_1)}^{t_0}_{\xi_0}(\hat s),t_0),$$
which is a curve in the tube $t_{\widehat r({\widetilde{(\gamma_1)}}^{t_0}_{\xi_0}(s))}(B)$ of $B$ of 
constant radius $\widehat r({\widetilde{(\gamma_1)}}^{t_0}_{\xi_0}(s))$.  
Then we have 
$$\begin{array}{l}
\hspace{0.5truecm}\displaystyle{\left.\frac{\overline{\nabla}}{d\hat s}\right\vert_{\hat s=s}
\left(\tau_{\gamma_{w(\widetilde{(\gamma_1)}^{t_0}_{\xi_0}(\hat s),t_0)}\vert_{[0,1]}}(\gamma'_1(s))\right)}\\
\displaystyle{=\cos(\sqrt{\varepsilon}k_0b\widehat r({\widetilde{(\gamma_1)}}^{t_0}_{\xi_0}(s),t_0))
\left.\frac{\overline{\nabla}}{d\hat s}\right\vert_{\hat s=s}
\left(\tau_{\gamma_{\bar w_{t_0,s}(\hat s)}\vert_{[0,1]}}(\gamma'_1(\hat s))\right)}\\
\displaystyle{=\sqrt{\varepsilon}k_0b
\sin(\sqrt{\varepsilon}k_0b\widehat r({\widetilde{(\gamma_1)}}^{t_0}_{\xi_0}(s),t_0))
w^1(\widetilde{(\gamma_1)}^{t_0}_{\xi_0}(s),t_0),}
\end{array}\leqno{(5.14)}$$
where $\displaystyle{\frac{\overline{\nabla}}{d\hat s}}$ in the left-hand-side (resp. the right-hand side) 
denotes the covariant derivative along the curve 
$\hat s\mapsto\exp^{\perp}\circ w(\widetilde{(\gamma_1)}^{t_0}_{\xi_0}(\hat s),t_0)$ 
(resp. $\hat s\mapsto\exp^{\perp}\circ\bar w_{t_0,s}(\hat s)$) with respect to $\overline{\nabla}$.  
From $(2.49),(4.9)$ and $(5.12)-(5.14)$, we obtain 
$$\begin{array}{l}
\displaystyle{(E_1^{t_0})_s\widehat u_{t_0}
=\sqrt{1-\widehat u(\widetilde{(\gamma_1)}^{t_0}_{\xi_0}(s),t_0)^2}
\left(\lambda_{t_0}(\widetilde{(\gamma_1)}^{t_0}_{\xi_0}(s))\right.}\\
\hspace{0.75truecm}\displaystyle{\left.+\widehat u(\widetilde{(\gamma_1)}^{t_0}_{\xi_0}(s),t_0)
\sum_{k\in{\cal K}}m_k^H\sqrt{\varepsilon}kb
\tan(\sqrt{\varepsilon}kb\widehat r(\widetilde{(\gamma_1)}^{t_0}_{\xi_0}(s),t_0))\right)}
\end{array}\leqno{(5.15)}$$
and hence 
$$\begin{array}{l}
\hspace{0truecm}\displaystyle{\left.\frac{\partial}{\partial s}\right\vert_{s=0}
((E_1^{t_0})_s\widehat u_{t_0})}\\
\displaystyle{=\frac{\widehat u(\xi_0,t_0)\sqrt{1-\widehat u(\xi_0,t_0)^2}}
{\cos(\sqrt{\varepsilon}k_0b\widehat r(\xi_0,t_0))}\widetilde{(e_1)}^{t_0}_{\xi_0}\lambda_{t_0}}\\
\hspace{0.5truecm}\displaystyle{+(1-\widehat u(\xi_0,t_0)^2)
\left(\lambda_{t_0}(\xi_0)+\widehat u(\xi_0,t_0)\sum_{k\in{\cal K}}m_k^H\sqrt{\varepsilon}kb
\tan(\sqrt{\varepsilon}kb\widehat r(\xi_0,t_0))\right)}\\
\hspace{1truecm}\displaystyle{\times\sqrt{\varepsilon}k_0b\tan(\sqrt{\varepsilon}k_0b\widehat r(\xi_0,t_0))}\\
\hspace{0.5truecm}\displaystyle{-\widehat u(\xi_0,t_0)
\left(\lambda_{t_0}(\xi_0)+\widehat u(\xi_0,t_0)\sum_{k\in{\cal K}}m_k^H\sqrt{\varepsilon}kb
\tan(\sqrt{\varepsilon}kb\widehat r(\xi_0,t_0))\right)^2}\\
\hspace{0.5truecm}\displaystyle{-\widehat u(\xi_0,t_0)(1-\widehat u(\xi_0,t_0)^2)
\sum_{k\in{\cal K}}\frac{m_k^H(\sqrt{\varepsilon}kb)^2\vert\vert({\rm grad}\,r_{t_0})_{c(\xi_0,t_0)}\vert\vert}
{\cos^2(\sqrt{\varepsilon}kb\widehat r(\xi_0,t_0))}.}
\end{array}\leqno{(5.16)}$$
Also, it is clear that 
$$\left.\frac{\partial}{\partial s}\right\vert_{s=0}((E_i^{t_0})_s\widehat u_{t_0})=0\quad\,\,
(i=2,\cdots,m^H).\leqno{(5.17)}$$
On the other hand, by using $E^{t_0}_j\widehat u_{t_0}=0$ ($j=m^H+1,\cdots,n$) and $(5.15)$, we have 
$$\begin{array}{l}
\hspace{0.5truecm}\displaystyle{
(\nabla^{t_0}du_{t_0})((E_j^{t_0})_{\xi_0},(E_j^{t_0})_{\xi_0})=-(\nabla^{t_0}_{(E_j^{t_0})_{\xi_0}}E_j^{t_0})
\widehat u_{t_0}}\\
\displaystyle{=-\overline g({\overline{\nabla}}^{f_{t_0}}_{(E_j^{t_0})_{\xi_0}}
f_{t_0\ast}(E_j^{t_0}),\,f_{t_0\ast}((E_1^{t_0})_0))
\cdot(E_1^{t_0})_0\widehat u_{t_0}}\\
\displaystyle{=-\frac{\sqrt{1-\widehat u(\xi_0,t_0)^2}}
{\sqrt{\cos^2(\sqrt{\varepsilon}k_0b\widehat r(\xi_0,t_0))
+\vert\vert({\rm grad}\,r_{t_0})_{c(\xi_0,t_0)}\vert\vert^2}}}\\
\hspace{0.5truecm}\displaystyle{\times
\overline g({\overline{\nabla}}^{f_{t_0}}_{(E_j^{t_0})_{\xi_0}}f_{t_0\ast}(E_j^{t_0}),\,
f_{t_0\ast}(\widetilde{(e_1)}^{t_0}_{\xi_0}))}\\
\hspace{0.5truecm}\displaystyle{\times
\left(\lambda_{t_0}(\xi_0)+\widehat u(\xi_0,t_0)\sum_{k\in{\cal K}}m_k^H\sqrt{\varepsilon}kb
\tan(\sqrt{\varepsilon}kb\widehat r(\xi_0,t_0))\right)}\\
\displaystyle{=-(1-\widehat u(\xi_0,t_0)^2)
\overline g({\overline{\nabla}}^{f_{t_0}}_{(E_j^{t_0})_{\xi_0}}f_{t_0\ast}(E_j^{t_0}),\,
{\widehat w}^1(\xi_0,t_0))}\\
\hspace{0.5truecm}\displaystyle{\times
\left(\lambda_{t_0}(\xi_0)+\widehat u(\xi_0,t_0)\sum_{k\in{\cal K}}m_k^H\sqrt{\varepsilon}kb
\tan(\sqrt{\varepsilon}kb\widehat r(\xi_0,t_0))\right)}\\
\displaystyle{=(1-\widehat u(\xi_0,t_0)^2)
\frac{\sqrt{\varepsilon}k_jb}{\tan(\sqrt{\varepsilon}k_jb\widehat r(\xi_0,t_0))}}\\
\hspace{0.5truecm}\displaystyle{\times
\left(\lambda_{t_0}(\xi_0)+\widehat u(\xi_0,t_0)\sum_{k\in{\cal K}}m_k^H\sqrt{\varepsilon}kb
\tan(\sqrt{\varepsilon}kb\widehat r(\xi_0,t_0))\right),}
\end{array}\leqno{(5.18)}$$
where we use also $E_i^{t_0}\widehat u_{t_0}=0$ ($i=2,\cdots,m^H$) and 
$$\begin{array}{l}
\hspace{0truecm}\displaystyle{
\overline g({\overline{\nabla}}^{f_{t_0}}_{(E_j^{t_0})_{\xi_0}}f_{t_0\ast}(E_j^{t_0}),\,
{\widehat w}^1(\xi_0,t_0))=-\bar g(\widehat A(f_{t_0\ast}((E_j^{t_0})_{\xi_0}),\,(E_j^{t_0})_{\xi_0})}\\
\hspace{5.2truecm}\displaystyle{
=-\frac{\sqrt{\varepsilon}k_jb}{\tan(\sqrt{\varepsilon}k_jb\widehat r(\xi_0,t_0))}}
\end{array}$$
(where $\widehat A$ denotes the shape operator of the geodesic sphere 
$f_{t_0}(M)\cap\exp^{\perp}(T^{\perp}_{c(\xi_0,t_0)}F)$ in the normal umbrella 
$\exp^{\perp}(T^{\perp}_{c(\xi_0,t_0)}F)$ for $-\widehat w^1(\xi_0,t_0)$). 
From $(5.16)-(5.18)$, we obtain 
$$\begin{array}{l}
\hspace{0truecm}\displaystyle{(\triangle_{t_0}\widehat u_{t_0})(\xi_0)
=\sum_{i=1}^{m^H}\left.\frac{\partial}{\partial s}\right\vert_{s=0}((E_i^{t_0})_s\widehat u_{t_0})
-\sum_{j=1}^{m^V}(\nabla^{t_0}_{(E_j^{t_0})_{\xi_0}}E_j^{t_0})\widehat u_{t_0}}\\
\hspace{0truecm}\displaystyle{=\frac{\widehat u(\xi_0,t_0)\sqrt{1-\widehat u(\xi_0,t_0)^2}}
{\cos(\sqrt{\varepsilon}k_0b\widehat r(\xi_0,t_0))}\cdot\widetilde{(e_1)}^{t_0}_{\xi_0}\lambda_{t_0}}\\
\hspace{0.5truecm}\displaystyle{+(1-\widehat u(\xi_0,t_0)^2)}\\
\hspace{1truecm}\displaystyle{\times\left(\lambda_{t_0}(\xi_0)+\widehat u(\xi_0,t_0)
\sum_{k\in{\cal K}}m_k^H\sqrt{\varepsilon}kb\tan(\sqrt{\varepsilon}kb\widehat r(\xi_0,t_0))\right)}\\
\hspace{1truecm}\displaystyle{\times\left(\sqrt{\varepsilon}k_0b\tan(\sqrt{\varepsilon}k_0b\widehat r(\xi_0,t_0))
+\sum_{k=1}^2\frac{m_k^V\sqrt{\varepsilon}kb}{\tan(\sqrt{\varepsilon}kb\widehat r(\xi_0,t_0))}\right)}\\
\hspace{0.5truecm}\displaystyle{-\widehat u(\xi_0,t_0)
\left(\lambda_{t_0}(\xi_0)+\widehat u(\xi_0,t_0)
\sum_{k\in{\cal K}}m_k^H\sqrt{\varepsilon}kb\tan(\sqrt{\varepsilon}kb\widehat r(\xi_0,t_0))\right)^2}\\
\hspace{0.5truecm}\displaystyle{-\widehat u(\xi_0,t_0)(1-\widehat u(\xi_0,t_0)^2)
\sum_{k\in{\cal K}}\frac{m_k^H(\sqrt{\varepsilon}kb)^2\vert\vert({\rm grad}\,r_{t_0})_{c(\xi_0,t_0)}\vert\vert}
{\cos^2(\sqrt{\varepsilon}kb\widehat r(\xi_0,t_0))}.}
\end{array}\leqno{(5.19)}$$
From $(5.9)$ and $(5.19)$, we obatin the following evolution equations for $\widehat u_t$ and $\widehat v_t$.  

\vspace{0.5truecm}

\noindent
{\bf Lemma 5.1.} {\sl The functions $\widehat u_t$'s ($t\in[0,T)$) satisfy the following evolution equation:
$$\begin{array}{l}
\hspace{0.5truecm}\displaystyle{\frac{\partial\widehat u}{\partial t}(\xi,t)-(\triangle_t\widehat u_t)(\xi)}\\
\hspace{0truecm}\displaystyle{=\overline H_t\sqrt{\varepsilon}k_0b\tan(\sqrt{\varepsilon}k_0b\widehat r(\xi,t))
(1-\widehat u(\xi,t)^2)}\\
\hspace{0.5truecm}\displaystyle{-\widehat u(\xi,t)(1-\widehat u(\xi,t)^2)
\sum_{k=1}^2\frac{m_k^V(\sqrt{\varepsilon}kb)^2}{\sin^2(\sqrt{\varepsilon}kb\widehat r(\xi,t))}}\\
\hspace{0.5truecm}\displaystyle{-\widehat u(\xi,t)(1-\widehat u(\xi,t)^2)
\sqrt{\varepsilon}k_0b\tan(\sqrt{\varepsilon}k_0b\widehat r(\xi,t))
\sum_{k=1}^2\frac{m_k^V\sqrt{\varepsilon}kb}{\tan(\sqrt{\varepsilon}kb\widehat r(\xi,t))}}\\
\hspace{0.5truecm}\displaystyle{-\widehat u(\xi,t)(1-\widehat u(\xi,t)^2)
\sqrt{\varepsilon}k_0b\tan(\sqrt{\varepsilon}k_0b\widehat r(\xi,t))}\\
\hspace{1truecm}\displaystyle{\times
\sum_{k\in{\cal K}}m^H_k\sqrt{\varepsilon}kb\tan(\sqrt{\varepsilon}kb\widehat r(\xi,t))}\\
\hspace{0.5truecm}\displaystyle{+\widehat u(\xi,t)
\left(\lambda_{t}(\xi)+\widehat u(\xi,t)
\sum_{k\in{\cal K}}m_k^H\sqrt{\varepsilon}kb\tan(\sqrt{\varepsilon}kb\widehat r(\xi,t))\right)^2}\\
\hspace{0.5truecm}\displaystyle{+\widehat u(\xi,t)(1-\widehat u(\xi,t)^2)
\sum_{k\in{\cal K}}\frac{m_k^H(\sqrt{\varepsilon}kb)^2\vert\vert({\rm grad}\,r_{t})_{c(\xi,t)}\vert\vert}
{\cos^2(\sqrt{\varepsilon}kb\widehat r(\xi,t))}}
\end{array}\leqno{(5.20)}$$
($(\xi,t)\in M\times[0,T_1)$).  
Also, the functions $\widehat v_t$'s ($t\in[0,T)$) satisfy the following evolution equation:
$$\begin{array}{l}
\hspace{0.5truecm}\displaystyle{\frac{\partial\widehat v}{\partial t}(\xi,t)-(\triangle_t\widehat v_t)(\xi)}\\
\hspace{0truecm}\displaystyle{=-\overline H_t\sqrt{\varepsilon}k_0b\tan(\sqrt{\varepsilon}k_0b\widehat r(\xi,t))
(\widehat v(\xi,t)^2-1)}\\
\hspace{0.5truecm}\displaystyle{+\widehat v(\xi,t)\left(1-\frac{1}{\widehat v(\xi,t)^2}\right)
\sum_{k=1}^2\frac{m_k^V(\sqrt{\varepsilon}kb)^2}{\sin^2(\sqrt{\varepsilon}kb\widehat r(\xi,t))}}\\
\hspace{0.5truecm}\displaystyle{+\widehat v(\xi,t)\left(1-\frac{1}{\widehat v(\xi,t)^2}\right)
\sqrt{\varepsilon}k_0b\tan(\sqrt{\varepsilon}k_0b\widehat r(\xi,t))
\sum_{k=1}^2\frac{m_k^V\sqrt{\varepsilon}kb}{\tan(\sqrt{\varepsilon}kb\widehat r(\xi,t))}}\\
\hspace{0.5truecm}\displaystyle{-\widehat v(\xi,t)\left(1-\frac{1}{\widehat v(\xi,t)^2}\right)
\sqrt{\varepsilon}k_0b\tan(\sqrt{\varepsilon}k_0b\widehat r(\xi,t))}\\
\hspace{1truecm}\displaystyle{\times
\sum_{k\in{\cal K}}m^H_k\sqrt{\varepsilon}kb\tan(\sqrt{\varepsilon}kb\widehat r(\xi,t))}\\
\hspace{0.5truecm}\displaystyle{-\widehat v(\xi,t)
\left(\lambda_{t}(\xi)+\frac{1}{\widehat v(\xi,t)}
\sum_{k\in{\cal K}}m_k^H\sqrt{\varepsilon}kb\tan(\sqrt{\varepsilon}kb\widehat r(\xi,t))\right)^2}\\
\hspace{0.5truecm}\displaystyle{-\widehat v(\xi,t)\left(1-\frac{1}{\widehat v(\xi,t)^2}\right)
\sum_{k\in{\cal K}}\frac{m_k^H(\sqrt{\varepsilon}kb)^2\vert\vert({\rm grad}\,r_{t})_{c(\xi,t)}\vert\vert}
{\cos^2(\sqrt{\varepsilon}kb\widehat r(\xi,t))}}\\
\hspace{0.25truecm}\displaystyle{-\frac{2}{\widehat v(\xi,t)}
\vert\vert({\rm grad}_t\widehat v_t)_{\xi}\vert\vert^2}
\end{array}\leqno{(5.21)}$$
($(\xi,t)\in M\times[0,T_1)$).}

\vspace{0.5truecm}

\noindent
{\it Proof.} The evolution equation $(5.20)$ is derived from $(5.9)$ and $(5.19)$ directly.  
The evolution equation $(5.21)$ is derived from $(5.20)$ and the following relations:
$$\frac{\partial\widehat v}{\partial t}=-\frac{1}{\widehat u_t^2}\frac{\partial\widehat u}{\partial t}\quad\,\,
{\rm and}\quad\,\,\triangle_t\widehat v_t=-\frac{1}{\widehat u_t^2}\triangle_t\widehat u_t
+\frac{2}{\widehat u_t^3}\vert\vert{\rm grad}_t\widehat u_t\vert\vert^2.$$
\begin{flushright}q.e.d.\end{flushright}

\section{Estimate of the volume}
We use the notations in Introduction and Sections 1-5.  
Assume that $\overline M$ is a rank one symmetric space, $F$ is invariant, $B$ is a closed geodesic ball of 
radius $r_B$ centered at $x_{\ast}$ in $F$ and that $r_0$ is radial with respect to $x_{\ast}$, where $r_B$ is 
a positive number smaller than the injective radius of $F$ at $x_{\ast}$.  
Note that the invariantness of $F$ means the totally geodesicness 
in the case where $\overline M$ is a sphere or a (real) hyperbolic space.  
Since $F$ is of rank one, each geodesic sphere in $F$ is homogeneous.  Hence, since $r_0$ is radial by the 
assumption, so are also $r_t$.  
For $X\in\widetilde S'(x_{\ast},1)$, denote by $\gamma_X$ the geodesic in $F$ having $X$ as the initial velocity 
vector (i.e., $\gamma_X(z)=\exp^F_{x_{\ast}}(zX)$).  
Then, since $r_t$ is radial, it is described as 
$r_t(\gamma_X(z))=\mathop{r}^{\circ}_t(z)$ ($X\in\widetilde S'(x_{\ast},1),\,\,z\in[0,r_B)$) for some 
function $\mathop{r}^{\circ}_t$ over $[0,r_B)$, where $\widetilde S'(x_{\ast},1)$ denotes the unit geodesic 
sphere in $F$ centered $x_{\ast}$ and $r_B$ denotes the radius of $B$.  
Then we have 
$$\begin{array}{c}
\displaystyle{({\rm grad}\,r_t)_{\gamma_X(z)}=(r^{\circ}_t)'(z)\gamma'_X(z),\,\,\,\,
(\triangle_Fr_t)(\gamma_X(z))=(r^{\circ}_t)''(z)}\\
\displaystyle{{\rm and}\,\,\,\,(\nabla^Fdr_t)
(({\rm grad}\,r_t)_{\gamma_X(z)},({\rm grad}\,r_t)_{\gamma_X(z)})=(r^{\circ}_t)'(z)^2(r^{\circ}_t)''(z).}
\end{array}\leqno{(6.1)}$$
Since $F$ is invariant, we have ${\cal K}=\{1\}$ (hence $k_0=1$) and 
$$r_{fo}(\gamma)=r_{co}(\gamma)=\left\{
\begin{array}{ll}
\displaystyle{\frac{\pi}{2b}} & (\varepsilon=1)\\
\displaystyle{\infty} & (\varepsilon=-1)
\end{array}\right.$$
for any normal geodesic $\gamma$ of $F$ and hence 
$$r_F=\left\{
\begin{array}{ll}
\displaystyle{\frac{\pi}{2b}} & (\varepsilon=1)\\
\displaystyle{\infty} & (\varepsilon=-1),
\end{array}\right.$$
where $r_{fo}(\gamma),r_{co}(\gamma)$ and $r_F$ are as in Introduction.  
Denote by $\nabla^t$ the Levi-Civita connection of $g_t$.  
In the sequel, assume that $t<T_1$.  
In this section, we estimate the volume of $M_t$ from below of in terms of the infimum and the maximum of 
the radius function $r_t$.  Furthermore, we show that 
$r_t$ and the average mean curvature $\overline H_t$ are uniformly bounded in terms of the estimate.  
Denote by $\pi_F^{\perp}$ the bundle projection of the normal bundle $T^{\perp}F$ of $F$.  
Set 
$$\widetilde W:=\{\xi\,\vert\,\xi\in T^{\perp}F\,\,{\rm s.t.}\,\,
\vert\vert\xi\vert\vert<r_F\}$$
and $W:=\exp^{\perp}(\widetilde W)$.  Define a submersion ${\rm pr}_F$ of $W$ onto $F$ by 
${\rm pr}_F(\exp^{\perp}(\xi)):=\pi_F^{\perp}(\xi)$, where $\xi\in\widetilde W$.  
Let $\widetilde r:\overline M\to{\Bbb R}$ be the distance function from $F$, where we note that 
$\widetilde r(\exp^{\perp}(\xi))=\vert\vert\xi\vert\vert$ holds for $\xi\in\widetilde W$.  
Define a function $\overline{\psi}$ over $[0,r_F)$ by 
$$
\overline{\psi}(s):=
\left(\mathop{\Pi}_{k=1}^2
\left(\frac{\sin(\sqrt{\varepsilon}kbs)}
{\sqrt{\varepsilon}kbs}\right)^{m_k^V}\right)
\cos^{m^H}(\sqrt{\varepsilon}bs).
$$
Since $r_F=\frac{\pi}{2\sqrt{\varepsilon}b}$, $\overline{\psi}$ is positive.  
From $(2.2)$ and $(2.3)$, we can give the following explicit descriptions of the volume element 
$dv_{\overline M}$ of $\overline M$ and the volume ${\rm Vol}(D_t)$.  

\vspace{0.5truecm}

\noindent
{\bf Proposition 6.1.} {\sl {\rm (i)} The volume element $dv_{\overline M}$ is given by 
$$\begin{array}{l}
\displaystyle{(dv_{\overline M})_p=(\overline{\psi}\circ\widetilde r)(p)
\left(((\exp^{\perp}\vert_{\widetilde S({\rm pr}_F(p),\widetilde r(p))})^{-1})^{\ast}
dv_{\widetilde S({\rm pr}_F(p),\widetilde r(p))}\right)_p}\\
\hspace{1.7truecm}
\displaystyle{\wedge(d\widetilde r)_p\wedge\left({\rm pr}_F^{\ast}dv_F\right)_p}
\hspace{4truecm}(p\in W)
\end{array}
\leqno{(6.2)}$$

{\rm (ii)} The volume ${\rm Vol}(D_t)$ is given by 
$${\rm Vol}(D_t)=v_{m^V}\int_{x\in B}\left(\int_0^{r_t(x)}
s^{m^V}\overline{\psi}(s)ds\right)dv_F.
\leqno{(6.3)}$$
}

\vspace{0.5truecm}

Denote by $\exp^F_x$ the exponential map of $F$ at $x$ and ${\widetilde S}'(x,a)$ 
the hypersphere of radius $a$ in $T_xF$ centered the origin.  
Since $B$ is star-shaped with respect to $x_{\ast}$, $B$ is described as 
$$B=\{\exp^F_{x_{\ast}}(zX)\,\vert\,X\in{\widetilde S}'(x_{\ast},1),\,\,0\leq z\leq r^B\}.$$
Since $F$ is a symmetric space, we can define the operators corresponding to 
$D^{co}_w$ and $D^{si}_w$ ($w\in T\overline M$) for each $X\in TF$.  
Deonte by $(D_F)^{co}_X$ and $(D_F)^{si}_X$ these operators.  
Since $F$ is of rank one, ${\rm det}((D_F)^{si}_{zX})$ and ${\rm det}((D_F)^{co}_{zX})$ are 
independent of the choice of $X$.  
Define a function $\psi_F$ over ${\Bbb R}$ by 
$\psi_F(z):={\rm det}((D_F)^{si}_{zX})$.  
Then we can describe the volume element $dv_F$ of $F$ and can estimate 
the volume ${\rm Vol}(B)$ of $B$ as follows.  

\vspace{0.5truecm}

\noindent
{\bf Proposition 6.2.} {\sl {\rm (i)} The volume element $dv_F$ is given by 
$$\begin{array}{l}
\displaystyle{(dv_F)_{\exp^F_{x_{\ast}}(zX)}=\psi_F(z)
\left(((\exp^F_{x_{\ast}}\vert_{{\widetilde S}'(x_{\ast},z)})^{-1})^{\ast}
dv_{{\widetilde S}'(x_{\ast},z)}\right)_{\exp^F_{x_{\ast}}(zX)}}\\
\hspace{3truecm}
\displaystyle{\wedge(d\widetilde z)_{\exp^F_{x_{\ast}}(zX)}}
\end{array}
\leqno{(6.4)}$$
for ${\exp^F_{x_{\ast}}(zX)}\in B$, where $\widetilde z$ is the distance function from $x_{\ast}$ in $F$.  

{\rm (ii)} The volume ${\rm Vol}(B)$ is given by 
$${\rm Vol}(B)=v_{m^H-1}\int_0^{r^B}z^{m^H-1}\psi_F(z)dz.
\leqno{(6.5)}$$
}

\vspace{0.5truecm}

\noindent
{\it Proof.} In similar to $(2.2)$, we have 
$$(\exp^F_{x_{\ast}})_{\ast}(Y)=\tau_{\gamma_{zX}}((D_F)^{si}_{zX}(Y)),$$
where $X\in{\widetilde S}'(x_{\ast},1)$ and $Y\in T_{zX}{\widetilde S}'(x_{\ast},z)$.  
From this relation, $(6.4)$ follows directly.  
The relation $(6.5)$ follows directly from $(6.4)$.  
\hspace{5truecm}q.e.d.

\vspace{0.5truecm}

Define a function $\delta_1$ and $\delta_2$ over $[0,r_F)$ by 
$$\delta_1(s):=\int_0^ss^{m^V}\overline{\psi}(s)ds$$
and 
$$\delta_2(s):=\int_0^s\frac{s^{m^V}\overline{\psi}(s)}{\cos(\sqrt{\varepsilon}bs)}ds,$$
respectively.  According to $(6.3)$, we have $\frac{{\rm Vol}(D_t)}{v_{m^V}{\rm Vol}(B)}\in\delta_1([0,r_F))$.  
Since $r_B$ is smaller than the injective radius of $F$ , $\overline{\psi}(s)$ and 
$\frac{\overline{\psi}(s)}{\cos(\sqrt{\varepsilon}bs)}$ are positive over $[0,r_B)$.  
Hence $\delta_i$ ($i=1,2$) are increasing.  
Set 
$$\hat r_1:=\delta_1^{-1}\left(\frac{{\rm Vol}(D)}{v_{m^V}{\rm Vol}(B)}\right)$$
and 
$$\hat r_2:=\delta_2^{-1}\left(\frac{{\rm Vol}(M)}{v_{m^V}}+\delta_2(\hat r_1)\right).$$
Denote by $(r_t)_{\rm max}$ (resp. $(r_t)_{\rm min}$) the maximum (resp. the minimum) 
of $r_t$.  Then we have 
$$\begin{array}{l}
\displaystyle{\delta_1(\widehat r_1)v_{m^V}{\rm Vol}(B)={\rm Vol}(D_t)}\\
\displaystyle{=v_{m^V}\int_{x\in B}\delta_1(r_t(x))\,dv_F
\left\{
\begin{array}{l}
\displaystyle{\geq v_{m^V}{\rm Vol}(B)\delta_1((r_t)_{\rm min})}\\
\displaystyle{\leq v_{m^V}{\rm Vol}(B)\delta_1((r_t)_{\rm max})}
\end{array}
\right.}
\end{array}$$
and hence $\delta_1((r_t)_{\rm min})\leq\delta_1(\widehat r_1)\leq\delta_1((r_t)_{\rm max})$, that is, 
$(r_t)_{\rm min}\leq\widehat r_1\leq(r_t)_{\rm max}$.  
Let $r^{\circ}_t$ be as in Section 5.  From Proposition 3.2 and $(6.4)$, we have 
$$\begin{array}{l}
\displaystyle{{\rm Vol}(M_t)=v_{m^V}\int_{x\in B}\frac{r_t(x)^{m^V}\overline{\psi}(r_t(x))
\sqrt{\cos^2(\sqrt{\varepsilon}br_t(x))+\vert\vert({\rm grad}\,r_t)_x\vert\vert^2}}
{\cos(\sqrt{\varepsilon}br_t(x))}dv_F}\\
\hspace{1.4truecm}\displaystyle{\geq v_{m^V}\int_{x\in B}
\frac{r_t(x)^{m^V}\overline{\psi}(r_t(x))\vert\vert({\rm grad}\,r_t)_x\vert\vert}
{\cos(\sqrt{\varepsilon}br_t(x))}dv_F}\\
\hspace{1.4truecm}\displaystyle{\geq v_{m^V}\int_{X\in\widetilde S'(x_{\ast},1)}\left(\int_0^{r_B}
\frac{r^{\circ}_t(z)^{m^V}\overline{\psi}(r^{\circ}_t(z))}{\cos(\sqrt{\varepsilon}br^{\circ}_t(z))}
\right.}\\
\hspace{4.65truecm}\displaystyle{\left.
\times z^{m^H-1}{\rm det}((D_F)^{si}_{zX})(r^{\circ}_t)'(z)\,dz\right)dv_{\widetilde S'(x_{\ast},1)},}
\end{array}\leqno{(6.6)}$$
where $x=\gamma_X(z)$.  Also, we have 
$$\begin{array}{l}
\displaystyle{z^{m^H-1}{\rm det}((D_F)^{si}_{zX})=\mathop{\Pi}_{k=1}^2
\left(\frac{\sin(\sqrt{\varepsilon}kbz)}{\sqrt{\varepsilon}kb}\right)^{m^H_k}}\\
\hspace{0truecm}\displaystyle{\geq
\left\{
\begin{array}{ll}
\displaystyle{
\mathop{\Pi}_{k=1}^2
\left(\frac{\sin(\sqrt{\varepsilon}kbz)}{\sqrt{\varepsilon}kb}\right)^{m^H_k}} & 
(\overline M\,:\,{\rm of}\,\,{\rm compact}\,\,{\rm type})\\
1 & (\overline M\,:\,{\rm of}\,\,{\rm noncompact}\,\,{\rm type})
\end{array}\right.}
\end{array}\leqno{(6.7)}$$
Set 
$$a_{r_B}:=\left\{
\begin{array}{ll}
\displaystyle{
\mathop{\Pi}_{k=1}^2
\left(\frac{\sin(\sqrt{\varepsilon}kbr_B)}{\sqrt{\varepsilon}kb}\right)^{m^H_k}} & 
(\overline M\,:\,{\rm of}\,\,{\rm compact}\,\,{\rm type})\\
1 & (\overline M\,:\,{\rm of}\,\,{\rm noncompact}\,\,{\rm type}).
\end{array}\right.$$
From $(6.6),\,(6.7)$ and $\max_{z\in[0,r_B]}\,r^{\circ}_t(z)=(r_t)_{\max}$ and 
$\min_{z\in[0,r_B]}\,r^{\circ}_t(z)=(r_t)_{\min}$, we can estimate the volume of $M_t$ from below as follows:
$$\begin{array}{l}
\displaystyle{{\rm Vol}(M_t)\geq 
a_{r_B}v_{m^V}v_{m^H-1}\int_{(r_t)_{\min}}^{(r_t)_{\max}}
\frac{s^{m^V}\overline{\psi}(s)}{\cos(\sqrt{\varepsilon}bs)}\,ds}\\
\hspace{1.425truecm}\displaystyle{=a_{r_B}v_{m^V}v_{m^H-1}(\delta_2((r_t)_{\max})-\delta_2((r_t)_{\min})).}
\end{array}\leqno{(6.8)}$$
On the other hand, since ${\rm Vol}(D_t)$ preserves invariantly along the volume-preserving mean 
curvature flow and ${\rm Vol}(M_t)$ is decreasing along the flow, we have 
${\rm Vol}(D_t)={\rm Vol}(D_0)$ and ${\rm Vol}(M_t)\leq{\rm Vol}(M_0)$.  
Hence we have 
$$\begin{array}{l}
\displaystyle{(r_t)_{\max}\leq\delta_2^{-1}\left(\delta_2(\hat r_1)+\frac{{\rm Vol}(M_0)}{a_{r_B}v_{m^V}v_{m^H-1}}
\right)}\\
\hspace{1.32truecm}\displaystyle{\leq\delta_2^{-1}
\left(\frac{{\rm Vol}(D_0)}{v_{m^V}{\rm Vol}(B)}+
\frac{{\rm Vol}(M_0)}{a_{r_B}v_{m^V}v_{m^H-1}}\right).}
\end{array}$$
Thus we obtain the following result.  

\vspace{0.5truecm}

\noindent
{\bf Proposition 6.3.} {\sl The family $\{r_t\}_{t\in[0,T_1)}$ is uniformly 
bounded as follows:
$$\sup_{(x,t)\in B\times[0,T_1)}\,r_t(x)\leq\delta_2^{-1}\left(\delta_2(\hat r_1)
+\frac{{\rm Vol}(M_0)}{a_{r_B}v_{m^V}v_{m^H-1}}\right).
\leqno{(6.10)}$$}

\vspace{0.5truecm}

For uniformly boundedness of the average mean curvatures $\vert\overline H\vert$, we have the following result.  

\vspace{0.5truecm}

\noindent
{\bf Proposition 6.4.} {\sl 
If $0<a_1\leq r_t(x)\leq a_2<r_F$ holds for all 
$(x,t)\in M\times[0,T_0]$ ($T_0<T_1$), then 
$\displaystyle{\mathop{\max}_{t\in[0,T_0]}\overline H_t\leq C(a_1,a_2)}$ 
holds for some constant $C(a_1,a_2)$ depending only on $a_1$ and $a_2$.}

\vspace{0.5truecm}

\noindent
{\it Proof.} According to $(3.8)$, we have 
$$\begin{array}{l}
\displaystyle{\overline H_t=\frac{1}{\int_Br_t^{m^V}\psi_{r_t}dv_F}\times\int_{x\in B}r_t(x)^{m^V}
\psi_{r_t}(x)}\\
\hspace{1.2truecm}\displaystyle{\left(\rho_{r_t}(x)
-\frac{\cos(\sqrt{\varepsilon}br_t(x))(r^{\circ}_t)''(z)}{(\cos^2(\sqrt{\varepsilon}br_t(x))
+(r^{\circ}_t)'(z)^2)^{3/2}}\right)dv_F}
\end{array}\leqno{(6.11)}$$
($x=\gamma_X(z)$), where $\rho_{r_t}$ and $\psi_{r_t}$ are the functions defined for $r_t$ in similar to 
$\rho_r$ and $\psi_r$ (see $(3.3)$ and $(3.4)$).  
Take any $(\xi,t)\in M\times [0,T_0]$ and let $x:=c(\xi,t)=r_X(z)$.  
According to the definitions of $\rho_{r_t}$ and $\psi_{r_t}$, we have the following estimates 
$$\rho_{r_t}(x)\leq\sum_{k=1}^2\frac{m_k^V\sqrt{\varepsilon}kb}{\tan(\sqrt{\varepsilon}kbr_t(x))}
+\sqrt{\varepsilon}b\tan(\sqrt{\varepsilon}br_t(x)).\leqno{(6.12)}$$
Since $0<a_1\leq r_t(x)\leq a_2<r_F$ by the assumption, 
we have 
$$\sum_{k=1}^2\frac{m_k^V\sqrt{\varepsilon}kb}{\tan(\sqrt{\varepsilon}kbr_t(x))}
+\sqrt{\varepsilon}b\tan(\sqrt{\varepsilon}br_t(x))\leq C_1(a_1,a_2)$$
for some positive constant $C_1(a_1,a_2)$ depending only on $a_1$ and $a_2$.  
Hence we have 
$$\frac{\int_B\rho_{r_t}r_t^{m^V}\psi_{r_t}dv_F}{\int_Br_t^{m^V}\psi_{r_t}dv_F}
\leq C_1(a_1,a_2).\leqno{(6.13)}$$
By a simple calculation, we have 
$$\begin{array}{l}
\hspace{0.5truecm}\displaystyle{\frac{d}{dz}\arctan
\left(\frac{(r^{\circ}_t)'(z)}{\cos(\sqrt{\varepsilon}br_t(x))}\right)}\\
\displaystyle{=\frac{\cos(\sqrt{\varepsilon}br_t(x))(r^{\circ}_t)''(z)}
{\cos^2(\sqrt{\varepsilon}br_t(x))+(r^{\circ}_t)'(z)^2}
+\frac{\sqrt{\varepsilon}b\sin(\sqrt{\varepsilon}br_t(x))(r^{\circ}_t)'(z)^2}
{\cos^2(\sqrt{\varepsilon}br_t(x))+(r^{\circ}_t)'(z)^2}.}
\end{array}\leqno{(6.14)}$$
By using this relation and the definition of $\psi_{r_t}$, we can derive 
$$\begin{array}{l}
\hspace{0.5truecm}\displaystyle{
-\int_{x\in B}\frac{\cos(\sqrt{\varepsilon}br_t(x))(r^{\circ}_t)''(z)
r_t(x)^{m^V}\psi_{r_t}(x)}{(\cos^2(\sqrt{\varepsilon}br_t(x))+(r^{\circ}_t)'(z)^2)^{3/2}}dv_F}\\
\hspace{0truecm}\displaystyle{=-\int_{x\in B}\frac{r_t(x)^{m^V}\psi_{r_t}(x)}
{\sqrt{\cos^2(\sqrt{\varepsilon}br_t(x))+(r^{\circ}_t)'(z)^2}}\cdot\frac{d}{dz}\arctan
\left(\frac{(r^{\circ}_t)'(z)}{\cos(\sqrt{\varepsilon}br_t(x))}\right)dv_F}\\
\hspace{0.5truecm}\displaystyle{+\int_{x\in B}
\frac{\sqrt{\varepsilon}b\sin(\sqrt{\varepsilon}br_t(x))(r^{\circ}_t)'(z)^2r_t(x)^{m^V}\psi_{r_t}(x)}
{(\cos^2(\sqrt{\varepsilon}br_t(x))+(r^{\circ}_t)'(z)^2)^{3/2}}dv_F}\\
\hspace{0truecm}\displaystyle{\leq-\int_{x\in B}
\left(\mathop{\Pi}_{k=1}^2\left(\frac{\sin(\sqrt{\varepsilon}kbr_t(x))}
{\sqrt{\varepsilon}kb}\right)^{m_k^V}\right)\cos^{m^H_{k_0}-1}(\sqrt{\varepsilon}k_0br_t(x))}\\
\hspace{1.7truecm}\displaystyle{\times\frac{d}{dz}\arctan
\left(\frac{(r^{\circ}_t)'(z)}{\cos(\sqrt{\varepsilon}br_t(x))}\right)\cdot r_t(x)^{m^V}dv_F}\\
\hspace{0.5truecm}\displaystyle{+\int_{x\in B}
\sqrt{\varepsilon}b\tan(\sqrt{\varepsilon}br_t(x))r_t(x)^{m^V}\psi_{r_t}(x)dv_F.}
\end{array}\leqno{(6.15)}$$
Since $0<a_1\leq r_t(x)\leq a_2<r_F$ by the assumption, 
we have 
$$\sqrt{\varepsilon}b\tan(\sqrt{\varepsilon}br_t(x))\leq C_2(a_1,a_2)$$
for some positive constant $C_2(a_1,a_2)$ depending only on $a_1$ and $a_2$.  
Hence we have 
$$\frac{\int_{x\in B}\tan(\sqrt{\varepsilon}br_t(x))r_t(x)^{m^V}\psi_{r_t}(x)dv_F}{\int_Br_t^{m^V}\psi_{r_t}dv_F}
\leq C_2(a_1,a_2).\leqno{(6.16)}$$
By integrating by parts and using $(6.4)$ and $(r^{\circ}_t)'(0)=(r^{\circ}_t)'(r_B)=0$ 
(which holds by the boundary condition (C1)), we have 
$$\begin{array}{l}
\hspace{0.5truecm}\displaystyle{-\int_{x\in B}
\left(\mathop{\Pi}_{k=1}^2\left(\frac{\sin(\sqrt{\varepsilon}kbr_t(x))}
{\sqrt{\varepsilon}kb}\right)^{m_k^V}\cdot\cos^{m^H-1}(\sqrt{\varepsilon}br_t(x))\right.}\\
\hspace{2truecm}\displaystyle{\left.\times\frac{d}{dz}\arctan
\left(\frac{(r^{\circ}_t)'(z)}{\cos(\sqrt{\varepsilon}br_t(x))}\right)\cdot r_t(x)^{m^V}\right)dv_F}\\
\hspace{0truecm}\displaystyle{=-\int_{\widetilde S'(x_{\ast},1)}\left(\int_0^{r_B}
\left(\mathop{\Pi}_{k=1}^2\left(\frac{\sin(\sqrt{\varepsilon}kbr_t(x))}
{\sqrt{\varepsilon}kb}\right)^{m_k^V}\cdot\cos^{m^H_{k_0}-1}(\sqrt{\varepsilon}k_0br_t(x))\right.\right.}\\
\hspace{1truecm}\displaystyle{\left.\left.\times\frac{d}{dz}\arctan
\left(\frac{(r^{\circ}_t)'(z)}{\cos(\sqrt{\varepsilon}br_t(x))}\right)\cdot r_t(x)^{m^V}\psi_F(z)z^{m^H-1}\right)
dz\right)dv_{\widetilde S'(x_{\ast},1)}}\\
\hspace{0truecm}\displaystyle{=
\int_{\widetilde S'(x_{\ast},1)}\left(\int_0^{r_B}
\left(\left.\frac{d}{ds}\right\vert_{s=r_t(x)}
\left(\mathop{\Pi}_{k=1}^2\left(\frac{\sin(\sqrt{\varepsilon}kbs)}
{\sqrt{\varepsilon}kb}\right)^{m_k^V}\cdot\cos^{m^H_{k_0}-1}(\sqrt{\varepsilon}k_0bs)s^{m^V}\right)
\right.\right.}\\
\hspace{1.5truecm}\displaystyle{
\left.\times\arctan\left(\frac{(r^{\circ}_t)'(z)}{\cos(\sqrt{\varepsilon}br_t(x))}\right)
(r^{\circ}_t)'(z)\psi_F(z)z^{m^H-1}dz\right)dv_{\widetilde S'(x_{\ast},1)}
}\\
\hspace{0.5truecm}\displaystyle{+\int_{\widetilde S'(x_{\ast},1)}\left(\int_0^{r_B}
\left(\mathop{\Pi}_{k=1}^2\left(\frac{\sin(\sqrt{\varepsilon}kbr_t(x))}
{\sqrt{\varepsilon}kb}\right)^{m_k^V}\cdot\cos^{m^H_{k_0}-1}(\sqrt{\varepsilon}k_0br_t(x))r_t(x)^{m^V}
\right.\right.}\\
\hspace{1.5truecm}\displaystyle{
\left.\left.\times\arctan\left(\frac{(r^{\circ}_t)'(z)}{\cos(\sqrt{\varepsilon}br_t(x))}\right)
\frac{d}{dz}\left(\psi_F(z)z^{m^H-1}\right)\right)dz\right)dv_{\widetilde S'(x_{\ast},1)}.}
\end{array}\leqno{(6.17)}$$
Also, we have 
$$\begin{array}{l}
\hspace{0.5truecm}\displaystyle{\int_{\widetilde S'(x_{\ast},1)}\left(\int_0^{r_B}\left(
\left.\frac{d}{ds}\right\vert_{s=r_t(x)}
\left(\mathop{\Pi}_{k=1}^2\left(\frac{\sin(\sqrt{\varepsilon}kbs)}
{\sqrt{\varepsilon}kb}\right)^{m_k^V}\cdot\cos^{m^H_{k_0}-1}(\sqrt{\varepsilon}k_0bs)s^{m^V}\right)
\right.\right.}\\
\hspace{1.5truecm}\displaystyle{\left.\left.\times
\arctan\left(\frac{(r^{\circ}_t)'(z)}{\cos(\sqrt{\varepsilon}br_t(x))}\right)
(r^{\circ}_t)'(z)\psi_F(z)z^{m^H-1}\right)dz\right)dv_{\widetilde S'(x_{\ast},1)}}\\
\hspace{0truecm}\displaystyle{\leq
\frac{\pi}{2}\int_{\widetilde S'(x_{\ast},1)}\left(\int_0^{r_B}
\left\vert\,\frac{d}{ds}\right\vert_{s=r_t(x)}\log
\left(\mathop{\Pi}_{k=1}^2\left(\frac{\sin(\sqrt{\varepsilon}kbs)}
{\sqrt{\varepsilon}kb}\right)^{m_k^V}\cdot\cos^{m^H_{k_0}-1}(\sqrt{\varepsilon}k_0bs)s^{m^V}\right)\right\vert}\\
\hspace{1.5truecm}\displaystyle{\times\mathop{\Pi}_{k=1}^2\left(\frac{\sin(\sqrt{\varepsilon}kbr_t(x))}
{\sqrt{\varepsilon}kb}\right)^{m_k^V}\cdot\cos^{m^H_{k_0}-1}(\sqrt{\varepsilon}k_0br_t(x))r_t(x)^{m^V}}\\
\hspace{1.5truecm}\displaystyle{\left.
\times\sqrt{\cos^2(\sqrt{\varepsilon}k_0br_t(x))+(r^{\circ}_t)'(z)^2}\psi_F(z)z^{m^H-1}dz\right)
dv_{\widetilde S'(x_{\ast},1)}}\\
\hspace{0truecm}\displaystyle{\leq
\frac{\pi}{2}\int_{x\in B}
\left\vert\,\left.\frac{d}{ds}\right\vert_{s=r_t(x)}\log
\left(\mathop{\Pi}_{k=1}^2\left(\frac{\sin(\sqrt{\varepsilon}kbs)}
{\sqrt{\varepsilon}kb}\right)^{m_k^V}\cdot\cos^{m^H_{k_0}-1}(\sqrt{\varepsilon}k_0bs)s^{m^V}\right)\right\vert}\\
\hspace{1.5truecm}\displaystyle{\times r_t(x)^{m_V}\psi_{r_t}(x)dv_F.}
\end{array}$$
Since $0<a_1\leq r_t(x)\leq a_2<r_F$ by the assumption, 
we have 
$$
\frac{\pi}{2}
\left\vert\left.\frac{d}{ds}\right\vert_{s=r_t(x)}\log
\left(\mathop{\Pi}_{k=1}^2\left(\frac{\sin(\sqrt{\varepsilon}kbs)}
{\sqrt{\varepsilon}kb}\right)^{m_k^V}\cdot\cos^{m^H_{k_0}-1}(\sqrt{\varepsilon}k_0bs)s^{m^V}\right)\right\vert
\leq C_3(a_1,a_2)
$$
for some positive constant $C_3(a_1,a_2)$ depending only on $a_1$ and $a_2$.  
Hence we have 
$$\begin{array}{l}
\hspace{0truecm}\displaystyle{\frac{1}{\int_Br_t^{m^V}\psi_{r_t}dv_F}\times}\\
\hspace{0truecm}\displaystyle{\int_{\widetilde S'(x_{\ast},1)}\left(\int_0^{r_B}\left(
\left.\frac{d}{ds}\right\vert_{s=r_t(x)}
\left(\mathop{\Pi}_{k=1}^2\left(\frac{\sin(\sqrt{\varepsilon}kbs)}
{\sqrt{\varepsilon}kb}\right)^{m_k^V}\cdot\cos^{m^H_{k_0}-1}(\sqrt{\varepsilon}k_0bs)s^{m^V}\right)
\right.\right.}\\
\hspace{1.5truecm}\displaystyle{\left.\left.\times
\arctan\left(\frac{(r^{\circ}_t)'(z)}{\cos(\sqrt{\varepsilon}br_t(x))}\right)
(r^{\circ}_t)'(z)\psi_F(z)z^{m^H-1}\right)dz\right)dv_{\widetilde S'(x_{\ast},1)}}\\
\hspace{0truecm}\displaystyle{\leq C_3(a_1,a_2).}
\end{array}\leqno{(6.18)}$$
Also we have 
$$\begin{array}{l}
\hspace{0.5truecm}\displaystyle{\int_{\widetilde S'(x_{\ast},1)}\left(\int_0^{r_B}
\left(\mathop{\Pi}_{k=1}^2\left(\frac{\sin(\sqrt{\varepsilon}kbr_t(x))}
{\sqrt{\varepsilon}kb}\right)^{m_k^V}\cdot\cos^{m^H-1}(\sqrt{\varepsilon}br_t(x))r_t(x)^{m^V}
\right.\right.}\\
\hspace{1.5truecm}\displaystyle{
\left.\left.\times\arctan\left(\frac{(r^{\circ}_t)'(z)}{\cos(\sqrt{\varepsilon}br_t(x))}\right)
\frac{d}{dz}\left(\psi_F(z)z^{m^H-1}\right)\right)dz\right)dv_{\widetilde S'(x_{\ast},1)}}\\
\hspace{0.5truecm}\displaystyle{\leq\frac{\pi}{2}\int_{\widetilde S'(x_{\ast},1)}\left(\int_0^{r_B}
\left(\frac{\sqrt{\cos^2(\sqrt{\varepsilon}br_t(x))+(r^{\circ}_t)'(z)^2}}{\cos(\sqrt{\varepsilon}br_t(x))}
\right.\right.}\\
\hspace{1.5truecm}\displaystyle{\times\mathop{\Pi}_{k=1}^2\left(\frac{\sin(\sqrt{\varepsilon}kbr_t(x))}
{\sqrt{\varepsilon}kb}\right)^{m_k^V}\cdot\cos^{m^H-1}(\sqrt{\varepsilon}br_t(x))}\\
\hspace{1.5truecm}\displaystyle{\left.\left.\times r_t(x)^{m^V}\cdot\frac{d}{dz}(\psi_F(z)z^{m^H-1})
\right)dz\right)dv_{\widetilde S'(x_{\ast},1)}}\\
\hspace{0.5truecm}\displaystyle{\leq\frac{\pi}{2}\int_{x\in B}
\left(\frac{1}{\cos(\sqrt{\varepsilon}br_t(x))}\cdot\frac{d}{dz}(\psi_F(z)z^{m^H-1})\right)dv_F.}
\end{array}$$
Since $0<a_1\leq r_t(x)\leq a_2<r_F$ by the assumption, 
we have 
$$\frac{1}{\cos(\sqrt{\varepsilon}br_t(x))}\cdot\frac{d}{dz}(\psi_F(z)z^{m^H-1})\leq C_4(a_1,a_2)$$
for some positive constant $C_4(a_1,a_2)$ depending only on $a_1$ and $a_2$.  
Hence we have 
$$\begin{array}{l}
\hspace{0truecm}\displaystyle{\frac{1}{\int_Br_t^{m^V}\psi_{r_t}dv_F}\times}\\
\hspace{0truecm}\displaystyle{\int_{\widetilde S'(x_{\ast},1)}\left(\int_0^{r_B}
\left(\mathop{\Pi}_{k=1}^2\left(\frac{\sin(\sqrt{\varepsilon}kbr_t(x))}
{\sqrt{\varepsilon}kb}\right)^{m_k^V}\cdot\cos^{m^H-1}(\sqrt{\varepsilon}br_t(x))r_t(x)^{m^V}
\right.\right.}\\
\hspace{1truecm}\displaystyle{
\left.\left.\times\arctan\left(\frac{(r^{\circ}_t)'(z)}{\cos(\sqrt{\varepsilon}br_t(x))}\right)
\frac{d}{dz}\left(\psi_F(z)z^{m^H-1}\right)\right)dz\right)dv_{\widetilde S'(x_{\ast},1)}}\\
\hspace{0truecm}\displaystyle{\leq C_4(a_1,a_2).}
\end{array}\leqno{(6.19)}$$
From $(6.15)-(6.19)$, we can derive 
$$\begin{array}{l}
\hspace{0truecm}\displaystyle{
-\frac{1}{\int_Br_t^{m^V}\psi_{r_t}dv_F}\cdot\int_{x\in B}\frac{\cos(\sqrt{\varepsilon}br_t(x))(r^{\circ}_t)''(z)
r_t(x)^{m^V}\psi_{r_t}(x)}{(\cos^2(\sqrt{\varepsilon}br_t(x))+(r^{\circ}_t)'(z)^2)^{3/2}}dv_F}\\
\hspace{0truecm}\displaystyle{\leq\sum_{i=2}^4C_i(a_1,a_2).}
\end{array}\leqno{(6.20)}$$
From $(6.11),(6.13)$ and $(6.20)$, we obtain 
$$\overline H_t\leq\sum_{i=1}^4C_i(a_1,a_2)\quad(t\in[0,T_0]).$$
\begin{flushright}q.e.d.\end{flushright}

\vspace{0.5truecm}

For uniformly positivity of the average mean curvatures $\vert\overline H\vert$, we have the following result.  

\vspace{0.5truecm}

\noindent
{\bf Proposition 6.5.} {\sl 
Assume that $\overline M$ is of non-compact type.  
If $r_t(x)\geq a>0$ holds for all 
$(x,t)\in M\times[0,T_0]$ ($T_0<T_1$), then 
$\displaystyle{\mathop{\min}_{t\in[0,T_0]}\overline H_t\geq\widehat C(a)}$ 
holds for some constant $\widehat C(a)$ depending only on $a$.}

\vspace{0.5truecm}

\noindent
{\it Proof.} Take any $x=\gamma_X(z)\in B$ and any $t\in[0,T_0]$.  
The function $\rho_{r_t}\psi_{r_t}$ is described as follows:
$$\begin{array}{l}
\displaystyle{\rho_{r_t}\psi_{r_t}(x)
=\mathop{\Pi}_{k=1}^2
\left(\frac{\sin(\sqrt{\varepsilon}kbr_t(x))}
{\sqrt{\varepsilon}kb}\right)^{m_k^V}\cdot
cos^{m^H}(\sqrt{\varepsilon}br_t(x))}\\
\hspace{1.75truecm}\displaystyle{\times\left(
\sum_{k=1}^2\frac{m_k^V\sqrt{\varepsilon}kb}{\tan(\sqrt{\varepsilon}kbr_t(x))}
-m^H\sqrt{\varepsilon}kb\tan(\sqrt{\varepsilon}kbr_t(x))\right.}\\
\hspace{2.25truecm}\displaystyle{\left.
+\frac{(r^{\circ}_t)'(z)^2\sqrt{\varepsilon}b\tan(\sqrt{\varepsilon}br_t(x))}
{\cos^2(\sqrt{\varepsilon}br_t(x))+(r^{\circ}_t)'(z)^2}\right),}
\end{array}$$
Since $\overline M$ is of non-compact type (i.e., $\varepsilon=-1$), we have 
$$\begin{array}{c}
\displaystyle{\cos(\sqrt{\varepsilon}kbr_t(x))=\cosh(kbr_t(x)),\,\,\,\,
\frac{\sqrt{\varepsilon}kb}{\tan(\sqrt{\varepsilon}kbr_t(x))}=\frac{kb}{\tanh(kbr_t(x))}}\\
\displaystyle{{\rm and}\,\,\,\,
\sqrt{\varepsilon}kb\tan(\sqrt{\varepsilon}kbr_t(x))=-kb\tanh(kbr_t(x)).}
\end{array}\leqno{(6.21)}$$
Hence, according to the above description of $\rho_{r_t}\psi_{r_t}(x)$, we have 
$\rho_{r_t}\psi_{r_t}\geq0$.  
Furthermore, since $r_t(x)\geq a>0$ by the assumption, 
$$\rho_{r_t}\psi_{r_t}(x)\geq C'(a)\,(>0)\leqno{(6.22)}$$
holds for some positive constant $C'(a)$ (depending only on $a$ (independent of $(x,t)\in B\times[0,T_0]$).  
This fact together with $\int_Br_t^{m^V}\psi_{r_t}dv_F={\rm Vol}(M_t)\leq{\rm Vol}(M_0)$ implies 
$$\frac{\int_Br_t^{m^V}\rho_{r_t}\psi_{r_t}\,dv_F}{\int_Br_t^{m^V}\psi_{r_t}\,dv_F}
\geq\frac{a^{m^V}C'(a){\rm Vol}(B)}{{\rm Vol}(M_0)}.\leqno{(6.23)}$$
Also, we have 
$$\begin{array}{l}
\hspace{0.5truecm}\displaystyle{
-\int_{x\in B}\frac{\cos(\sqrt{\varepsilon}br_t(x))(r^{\circ}_t)''(z)
r_t(x)^{m^V}\psi_{r_t}(x)}{(\cos^2(\sqrt{\varepsilon}br_t(x))+(r^{\circ}_t)'(z)^2)^{3/2}}dv_F}\\
\hspace{0truecm}\displaystyle{=-\int_{x\in B}\frac{r_t(x)^{m^V}\psi_{r_t}(x)}
{\sqrt{\cos^2(\sqrt{\varepsilon}br_t(x))+(r^{\circ}_t)'(z)^2}}\cdot\frac{d}{dz}\arctan
\left(\frac{(r^{\circ}_t)'(z)}{\cos(\sqrt{\varepsilon}br_t(x))}\right)dv_F}\\
\hspace{0.5truecm}\displaystyle{+\int_{x\in B}
\frac{\sin(\sqrt{\varepsilon}br_t(x))(r^{\circ}_t)'(z)^2r_t(x)^{m^V}\psi_{r_t}(x)}
{(\cos^2(\sqrt{\varepsilon}br_t(x))+(r^{\circ}_t)'(z)^2)^{3/2}}dv_F}.
\end{array}$$
By integrating the first term in the right-hand side by parts and using $(6.21)$, we can show 
$$-\int_{x\in B}\frac{\cos(\sqrt{\varepsilon}br_t(x))(r^{\circ}_t)''(z)
r_t(x)^{m^V}\psi_{r_t}(x)}{(\cos^2(\sqrt{\varepsilon}br_t(x))+(r^{\circ}_t)'(z)^2)^{3/2}}dv_F\geq0.
\leqno{(6.24)}$$
From $(6.11),\,(6.23),\,(6.24)$ and the arbitrariness of $t$, we can derive 
$$\mathop{\min}_{t\in[0,T_0]}\overline H_t\geq\frac{a^{m^V}C'(a){\rm Vol}(B)}{{\rm Vol}(M_0)}.$$
\begin{flushright}q.e.d.\end{flushright}

\section{Proof of Theorem A} 
In this section, we shall prove Theorem A in terms of the evolution equation $(5.21)$ of $\widehat v_t$ 
and Propositions 6.4 and 6.5.  

\vspace{0.5truecm}

\noindent
{\it Proof of Theorem A.} 
Suppose that $T_1<T$.  Take any $t_0\in(T_1,T)$.  Set 
$$\beta_1(t_0):=\mathop{\min}_{(x,t)\in B\times[0,t_0]}r_t(x)\,\,(>0)\,\,\,\,
{\rm and}\,\,\,\,\beta_2(t_0):=\mathop{\max}_{(x,t)\in B\times[0,t_0]}r_t(x)\,\,(<\infty).$$
According to Propostions 6.4 and 6.5, we have 
$$0<\widehat C(\beta_1(t_0))<\overline H_t<C(\beta_1(t_0),\beta_2(t_0))\qquad(t\in[0,T_1)),$$
where we note that $r_F=\infty$ because $\overline M$ is of non-compact type.  
By using the evolution equation $(5.21)$ for $\widehat v_t$, this inequality, 
$(1-1/\widehat v_t^2)\widehat v_t\leq\widehat v_t,\,\widehat v_t\geq 1$ and $(6.21)$, 
we can derive 
$$\begin{array}{l}
\hspace{0.5truecm}\displaystyle{\frac{\partial\widehat v}{\partial t}(\xi,t)-(\triangle_t\widehat v_t)(\xi)}\\
\displaystyle{\leq\widehat v(\xi,t)\left(\sum_{k=1}^2\frac{m^V_k(kb)^2}{\sinh^2(kb\widehat r(\xi,t))}
+b\tanh(b\widehat r(\xi,t))\sum_{k=1}^2\frac{m^V_kkb}{\tanh(kb\widehat r(\xi,t))}\right)}\\
\hspace{0.5truecm}\displaystyle{-\widehat v(\xi,t)^2\widehat C(\beta_1(t_0))b\tanh(b\widehat r(\xi,t))
+C(\beta_1(t_0),\beta_2(t_0))}
\end{array}\leqno{(7.1)}$$
($t\in[0,T_1)$).  
For simplicity, set 
$$K_1(\beta_1(t_0),\beta_2(t_0)):=
\sum_{k=1}^2\frac{m^V_k(kb)^2}{\sinh^2(kb\beta_1(t_0))}
+b\tanh(b\beta_2(t_0))\sum_{k=1}^2\frac{m^V_kkb}{\tanh(kb\beta_1(t_0))}$$
and 
$$K_2(\beta_1(t_0)):=\widehat C(\beta_1(t_0))b\tanh(b\beta_1(t_0)).$$
Take any $t_1\in[0,T_1)$.  Let $(\xi_2,t_2)\in M\times[0,t_1]$ be a point attaining the maximum of 
$\widehat v$ over $M\times[0,t_1]$.  Since $\widehat v_{t_2}=1$ along $\partial M$, $(\xi_2,t_2)$ belongs 
to the interior of $M\times[0,t_1]$.  
Then we have $\frac{\partial\widehat v}{\partial t}(\xi_2,t_2)=0$ and 
$(\triangle_{t_2}\widehat v_{t_2})(\xi_2)\leq0$, that is, 
$\frac{\partial \widehat v}{\partial t}(\xi_2,t_2)-(\triangle_{t_2}\widehat v_{t_2})(\xi_2)\geq0$.  
Hence, from $(7.1)$, we can derive 
$$\begin{array}{l}
\hspace{0.5truecm}\displaystyle{\mathop{\max}_{(\xi,t)\in M\times[0,t_1]}\widehat v(\xi,t)
=\widehat v(\xi_2,t_2)}\\
\displaystyle{\leq\frac{K_1(\beta_1(t_0),\beta_2(t_0))+\sqrt{K_1(\beta_1(t_0),\beta_2(t_0))^2
+4K_2(\beta_1(t_0))C(\beta_1(t_0),\beta_2(t_0))}}{2}.}
\end{array}$$
From the arbitrariness of $t_1$, we obtain 
$$\begin{array}{l}
\hspace{0.5truecm}\displaystyle{\mathop{\sup}_{(\xi,t)\in M\times[0,T_1)}\widehat v(\xi,t)}\\
\displaystyle{\leq\frac{K_1(\beta_1(t_0),\beta_2(t_0))+\sqrt{K_1(\beta_1(t_0),\beta_2(t_0))^2
+4K_2(\beta_1(t_0))C(\beta_1(t_0),\beta_2(t_0))}}{2},}
\end{array}$$
which implies that $M_{T_1}$ is a tube over $B$ and that, furthermore, so are also $M_t$'s 
($t\in[T_1,T_1+\varepsilon)$) for a sufficiently small positive number $\varepsilon$ (smaller than $t_0-T_1$).  
This contradicts the definition of $T_1$.  
Therefore we obtain $T_1=T$.  That is, $M_t$ ($t\in[0,T)$) remain to be tubes over $B$.  
\hspace{8.85truecm}q.e.d.

\section{Proofs of Theorem B and C} 
We use the notations in Introduction and Sections 1-7.  Assume that $\overline M$ is a rank one symmetric space 
of non-compact type other than a (real) hyperbolic space and $r_0$ is radial (hence so are lso $r_t$'s 
($0\leq t<T$)).  
Set $r_{\inf}:=\inf_{(x,t)\in B\times[0,T)}r_t(x)$ and $r_{\sup}:=\sup_{(x,t)\in B\times[0,T)}r_t(x)$.  
In this section, we shall prove Theorems B and C.  
%
Since $\overline M$ is a symmetric space, we have $\overline{\nabla}\,\overline R=0$.  Hence, according to 
$(6.1)$ in [CM2], we obtain the following evolution equation for $\vert\vert A_t\vert\vert_t^2$.  

\vspace{0.5truecm}

\noindent
{\bf Lemma 8.1([CM2]).} {\sl The family $\{\vert\vert A_t\vert\vert^2_t\}_{t\in[0,T)}$ satisfies 
$$\begin{array}{l}
\displaystyle{\frac{\partial\vert\vert A_t\vert\vert_t^2}{\partial t}-\triangle_t\vert\vert A_t\vert\vert_t^2}\\
\displaystyle{=-2\vert\vert\nabla^tA_t\vert\vert_t^2+2\vert\vert A_t\vert\vert_t^2(\vert\vert A_t\vert\vert_t^2
+\overline{\rm Ric}(N_t,N_t))}\\
\hspace{0.5truecm}\displaystyle{-2\overline H_t\left({\rm Tr}\,A_t^3
+Tr_{g_t}^{\bullet}\overline R(A_t(\bullet),N_t,N_t,\bullet)\right)}\\
\hspace{0.5truecm}\displaystyle{-4{\rm Tr}_{g_t}^{\bullet}{\rm Tr}_{g_t}^{\cdot}
\overline R(\cdot,A_t(\bullet),A_t(\bullet),\cdot)
+4{\rm Tr}_{g_t}^{\bullet}{\rm Tr}_{g_t}^{\cdot}
\overline R(\cdot,A_t(\bullet),\bullet,A_t(\cdot)).}
\end{array}$$
}

\vspace{0.5truecm}

Set $\kappa:=\frac{1}{2\sup_{(\xi,t)\in M\times[0,T)}\widehat v(\xi,t)^2}$.  
According to Theorem A, we have $\kappa<\infty$.  
Define a function $\varphi$ over $[0,1/\sqrt{\kappa})$ by 
$\displaystyle{\varphi(s):=\frac{s^2}{1-\kappa s^2}}$ and 
a function $\Phi_t$ over $M$ by 
$$\Phi_t(\xi):=(\varphi\circ\widehat v_t)(\xi)\vert\vert(A_t)_{\xi}\vert\vert_t^2\quad(\xi\in M).$$
Also, define a function $\Phi$ over $M\times[0,T)$ by $\Phi(\xi,t):=\Phi_t(\xi)\,\,((\xi,t)\in M\times[0,T))$.  
By a simple calculation, we can derive the following evolution equation for $\Phi_t$'s ($t\in[0,T)$).  

\vspace{0.5truecm}

\noindent
{\bf Lemma 8.2.} {\sl The family $\{\Phi_t\}_{t\in[0,T)}$ satisfies 
$$\begin{array}{l}
\hspace{0.5truecm}\displaystyle{\frac{\partial\,\Phi}{\partial t}-\triangle_t\Phi_t}\\
\displaystyle{=(\varphi'\circ\widehat v_t)\vert\vert A_t\vert\vert_t^2
\left(\frac{\partial\,\widehat v}{\partial t}-\triangle_t\widehat v_t\right)+(\varphi\circ\widehat v_t)
\left(\frac{\partial\,\vert\vert A_t\vert\vert_t^2}{\partial t}-\triangle_t\vert\vert A_t\vert\vert_t^2\right)}\\
\hspace{0.5truecm}\displaystyle{-(\varphi''\circ\widehat v_t)\vert\vert A_t\vert\vert_t^2
\vert\vert{\rm grad}_t\widehat v\vert\vert_t^2-2g_t({\rm grad}_t((\varphi\circ\widehat v_t),\,{\rm grad}_t
\vert\vert A_t\vert\vert_t^2).}
\end{array}$$
}

\vspace{0.5truecm}

By using these lemmas, $(5.21)$, Propositions 6.3, 6.4 and 6.5, we shall prove Theorem B.  

\vspace{0.5truecm}

\noindent
{\it Proof of Theorem B.} Suppose that $\inf_{(\xi,t)\in M\times[0,T)}\widehat r(\xi,t)>0$ and $T<\infty$.  
Then we suffice to show that $T=\infty$ and that $M_t$ converges to a tube of constant mean curvature over $B$ 
(in $C^{\infty}$-topology) as $t\to\infty$.  
Set $a_1:=\inf_{(\xi,t)\in M\times[0,T)}\widehat r(\xi,t)$ and 
$a_2:=\sup_{(\xi,t)\in M\times[0,T)}\widehat r(\xi,t)$.  
According to Proposition 6.3, we have 
$$a_2\leq \delta_2^{-1}\left(\delta_2(\widehat r_1)+\frac{{\rm Vol}(M_0)}{a_{r_B}v_{m^V}v_{m^H-1}}\right).$$
According to Propositions 6.4 and 6.5, we have 
$$0<\widehat C(a_1)\leq \overline H_t\leq C(a_1,a_2)<\infty\qquad(t\in[0,T)).\leqno{(8.1)}$$
According to the first inequality in Page 555 in [EH] (see also $(6.4)$ and $(6.5)$ of P502 in 
[CM2]), we have 
$$\begin{array}{l}
\displaystyle{-2g_t({\rm grad}_t(\varphi\circ\widehat v_t),\,
{\rm grad}_t(\vert\vert A_t\vert\vert_t^2))}\\
\displaystyle{\leq-\frac{1}{\varphi\circ\widehat v_t}g_t({\rm grad}_t\Phi_t,
{\rm grad}_t(\varphi\circ\widehat v_t))+2(\varphi\circ\widehat v_t)\vert\vert\nabla^tA_t\vert\vert_t^2}\\
\hspace{0.5truecm}\displaystyle{+\frac{3}{2(\varphi\circ\widehat v_t)}
\vert\vert A_t\vert\vert_t^2\vert\vert{\rm grad}_t(\varphi\circ\widehat v_t)\vert\vert_t^2.}
\end{array}\leqno{(8.2)}$$
Also, according to Kato's inequality, we have 
$$\vert\vert{\rm grad}_t\vert\vert A_t\vert\vert_t\,\vert\vert_t^2
\leq\vert\vert\nabla^tA_t\vert\vert_t^2.\leqno{(8.3)}$$
By using these relations, $(5.21)$, Lemmas 8.1 and 8.2, we can derive 
$$\begin{array}{l}
\hspace{0.5truecm}\displaystyle{\frac{\partial\,\Phi}{\partial t}-\triangle_t\Phi_t}\\
\displaystyle{\leq(\varphi'\circ\widehat v_t)\vert\vert A_t\vert\vert_t^2
\left\{K_1(a_1,a_2)\widehat v_t-K_2(a_1)\widehat v_t^2-\frac{2}{\widehat v_t}
\vert\vert{\rm grad}_t\widehat v_t\vert\vert_t^2\right.}\\
\hspace{3truecm}\displaystyle{\left.-\widehat v_t\left(\lambda_t+\frac{m^Hb\tanh(b\widehat r_t)}{\widehat v_t}
\right)^2\right\}}\\
\hspace{0.5truecm}\displaystyle{+2(\varphi\circ\widehat v_t)\vert\vert A_t\vert\vert_t^4
+2(\varphi\circ\widehat v_t)\vert\vert A_t\vert\vert_t^2\overline{\rm Ric}(N_t,N_t)}\\
\hspace{0.5truecm}\displaystyle{-2(\varphi\circ\widehat v_t)\overline H_t\cdot{\rm Tr}\,A_t^3
-2(\varphi\circ\widehat v_t)\overline H_t{\rm Tr}_{g_t}^{\bullet}\overline R(A_t(\bullet),N_t,N_t,\bullet)}\\
\hspace{0.5truecm}\displaystyle{-4(\varphi\circ\widehat v_t){\rm Tr}_{g_t}^{\bullet}{\rm Tr}_{g_t}^{\cdot}
\overline R(\cdot,A_t(\bullet),A_t(\bullet),\cdot)}\\
\hspace{0.5truecm}\displaystyle{+4(\varphi\circ\widehat v_t){\rm Tr}_{g_t}^{\bullet}{\rm Tr}_{g_t}^{\cdot}
\overline R(\cdot,A_t(\bullet),\bullet,A_t(\cdot))}\\
\hspace{0.5truecm}\displaystyle{-(\varphi''\circ\widehat v_t)\vert\vert A_t\vert\vert_t^2\cdot
\vert\vert{\rm grad}_t\widehat v_t\vert\vert_t^2-\frac{1}{\varphi\circ\widehat v_t}
g_t({\rm grad}_t\Phi_t,{\rm grad}_t(\varphi\circ\widehat v_t))}\\
\hspace{0.5truecm}\displaystyle{+\frac{3}{2(\varphi\circ\widehat v_t)}
\vert\vert A_t\vert\vert_t^2\cdot\vert\vert{\rm grad}_t(\phi\circ\widehat v_t)\vert\vert_t^2}
\end{array}\leqno{(8.4)}$$
According to $(8.1)$, $\{\tanh(b\widehat r_t)\}_{t\in[0,T)}$ and 
$\{\overline{\rm Ric}(N_t,N_t)\}_{t\in[0,T)}$ are uniformly bounded.  
Hence, by the same discussion as in the second-half part in the proof of Theorem A of [CM2] 
(see Line 15 of Page 502-Line 14 of Page 503), we can derive 
$$\begin{array}{l}
\hspace{0truecm}\displaystyle{\frac{\partial\,\Phi}{\partial t}-\triangle_t\Phi_t
\leq-\frac{\kappa}{2}\Phi_t^2+K_3\Phi_t+K_4\sqrt{\Phi_t}}\\
\hspace{2.55truecm}\displaystyle{-\frac{1}{\varphi\circ\widehat v_t}g_t({\rm grad}_t\Phi_t,\,{\rm grad}_t(\varphi\circ\widehat v_t))}
\end{array}\leqno{(8.5)}$$
for some positive constants $K_3$ and $K_4$.  
Take any $t_0\in[0,T)$.  Let $(\xi_1,t_1)$ be a point attainning $\max_{(\xi,t)\in M\times[0,t_0]}\Phi_t(\xi)$.  
We consider the case where $\xi_1$ belongs to the boundary of $M$.  
Define a radial function $\widetilde r_t$ over the geodesic ball $\widetilde B$ of radius $2r_B$ centered at 
$x_{\ast}$ in $F$ by 
$$\widetilde r_t(\gamma_X(z)):=
\left\{
\begin{array}{ll}
\displaystyle{r^{\circ}_t(z)} & (0\leq z\leq r_B)\\
\displaystyle{r^{\circ}_t(2r_B-z)} & (r_B\leq z\leq 2r_B)
\end{array}\right.$$
Set $\widetilde M:=t_{\widetilde r_0}(\widetilde B)$, which includes $M$.  
Then it is easy to show that the volume-preserving mean curvature flow $\widetilde f_t$ 
($t\in[0,T)$) starting from $\widetilde M$ satisfies 
$\widetilde f_t(\widetilde M)=\exp^{\perp}(t_{\widetilde r_t}(\widetilde B))$.  
Also, it is clear that $(\xi_1,t_1)$ is a point attainning $\max_{(\xi,t)\in\widetilde M\times[0,t_0]}
\widetilde{\Phi}_t(\xi)$ and that $\xi_1$ belongs to the interior of $\widetilde M$, where 
$\widetilde{\Phi}_t$ is the function defined for $\widetilde f_t$ similar to $\Phi_t$.  
Thus we may assume that $\xi_1$ belongs to the interior of $M$ without loss of generality.  
Set $b_1:=\Phi_{t_1}(\xi_1)$.  Assume that $b_1>1$.  
Then, it follows from $(8.5)$ that 
$$0\leq\frac{\partial\Phi}{\partial t}(\xi_1,t_1)-(\triangle_{t_1}\Phi_{t_1})(\xi_1)
\leq-\frac{\kappa b_1^2}{2}+(K_3+K_4)b_1,$$
that is, 
$b_1\leq\frac{2}{\kappa}(K_3+K_4)$.  
Hence we obtain 
$$\max_{(\xi,t)\in M\times[0,t_0]}\Phi_t(\xi)
=b_1\leq\max\left\{1,\frac{2}{\kappa}(K_3+K_4)\right\}.$$
Furthermore, by the arbitariness of $t_0$, we obtain 
$$\sup_{(\xi,t)\in M\times[0,T)}\Phi_t(\xi)\leq\max\left\{1,\frac{2}{\kappa}(K_3+K_4)\right\}.$$
Thus $\{\Phi_t\}_{t\in[0,T)}$ is uniformly bounded.  
On the other hand, we have $\varphi\circ\widehat v_t\geq 2\widehat v_t^2\geq 2$.  
From these facts, it follows that $\{\vert\vert A_t\vert\vert_t^2\}_{t\in[0,T)}$ is uniformly bounded.  
Hence, by the discussion in [Hu1,2], it is shown that 
$\{\vert\vert(\nabla^t)^jA_t\vert\vert_t^2\}_{t\in[0,T)}$ also is uniformly bounded 
for any positive integar $j$, where we use $T<\infty$ and 
$g_t({\rm grad}_t\vert\vert(\nabla^t)^jA_t\vert\vert_t^2,\,N_t)$\newline
$=0$ along $\partial B$ 
(This fact holds because ${\rm grad}_tr_t=0$ along $\partial B$ by the assumption).  
Hence, since $\{f_t\}_{t\in[0,T)}$ is equicontinuous in $C^{\infty}$-norm, there exists a sequence 
$\{t_k\}_{k=1}^{\infty}$ such that $\lim\limits_{k\to\infty}t_k=T$ and that 
$\{f_{t_k}\}_{k=1}^{\infty}$ converges to some $C^{\infty}$-embedding $f_T$ 
in the $C^{\infty}$-topology by Arzel$\acute{\rm a}$-Ascoli's theorem.  
By the standard discussion as in [CM1, P2075-2077], it is shown that $f_T(M)$ is of constant mean curvature.  
On the other hand, according to $(4.17)$, we have 
$$\begin{array}{l}
\displaystyle{\frac{\partial r}{\partial t}
=\frac{\triangle_Fr_t}{\cosh^2(br_t)+\vert\vert{\rm grad}\,r_t\vert\vert^2}
+\frac{\sqrt{\cosh^2(br_t)+\vert\vert{\rm grad}\,r_t\vert\vert^2}}{\cosh(br_t)}(\overline H_t-\rho_{r_t})}\\
\hspace{0.65truecm}\displaystyle{=\frac{\sqrt{\cosh^2(br_t)+\vert\vert{\rm grad}\,r_t\vert\vert^2}}
{\cosh(br_t)}(\overline H_t-H_t).}
\end{array}$$
Hence, for any $\tau_1$ and $\tau_2$ of $[0,T)$ ($\tau_1<\tau_2$), we have
$$\begin{array}{l}
\hspace{0.5truecm}\displaystyle{\mathop{\max}_{x\in B}\vert r_{\tau_1}(x)-r_{\tau_2}(x)\vert
\leq\mathop{\max}_{x\in M}\int_{\tau_1}^{\tau_2}
\left\vert\frac{\partial r}{\partial t}\right\vert\,dt}\\
\hspace{0truecm}\displaystyle{\leq\mathop{\max}_{x\in M}\int_{\tau_1}^{\tau_2}
\vert\overline H_t-\rho_{r_t}\vert\sqrt{\cosh^2(br_t)+\vert\vert{\rm grad}\,r_t\vert\vert^2}\,dt}\\
\hspace{0truecm}\displaystyle{\leq C\varepsilon(t)(\tau_2-\tau_1),}
\end{array}$$
where $t$ is an element of $(\tau_1,\tau_2)$ and $\varepsilon(t)\to 0$ as $t\to T$.  
From this fact, we can show that $r_t$ converges to $r_T$ in the $C^{\infty}$-topology by the discussion as in 
[Hu3], where $r_T$ is the radius function of $f_T$.  
Thus $f_t$ converges to $f_T$ in the $C^{\infty}$-topology as $t\to T$.  
It is clear that $f_T(M)$ also is a tube over $B$.  
Hence, according to Proposition 4.2, there exists the volume-preserving mean curvature flow starting from $f_T$ 
in short time.  Hence the volume-preserving mean curvature flow $f_t$ starting from $M$ is continued after $T$.  
This contradicts the definition of $T$.  
Therefore $T=\infty$ or $\inf_{(\xi,t)\in M\times[0,T)}\widehat r(\xi,t)=0$ holds.  
Furthermore, in the case of $T=\infty$, we can show that $f_t$ converges to a tube of constant mean curvature 
over $B$ in $C^{\infty}$-topology as $t\to\infty$ by imitating the discussion in P2075-2080 of [CM1].  
\hspace{3.75truecm}q.e.d.

\vspace{0.5truecm}

Next we shall prove Theorem C.  

\vspace{0.5truecm}

\noindent
{\it Proof of Theorem C.} 
Suppose that $\inf_{(\xi,t)\in M\times[0,T)}\widehat r(\xi,t)=0$.  
Then there exists a sequence $\{(\xi_k,t_k)\}_{k=1}^{\infty}$ such that 
$\widehat r(\xi_k,t_k)<\frac{1}{k}$ ($k\in{\Bbb N}$).  
By using $(6.8)$ and $(r_t)_{\min}\leq\widehat r_1\leq(r_t)_{\max}$, we have 
$$\begin{array}{l}
\displaystyle{{\rm Vol}(M_0)\geq{\rm Vol}(M_{t_k})\geq a_{r_B}v_{m^V}v_{m^H-1}
\left(\delta_2((r_{t_k})_{\max})-\delta_2((r_{t_k})_{\min})\right)}\\
\hspace{3.525truecm}\displaystyle{>v_{m^V}v_{m^H-1}\left(\delta_2(\widehat r_1)-\delta_2(1/k)\right)}
\end{array}$$
and hence 
$$\begin{array}{l}
\displaystyle{{\rm Vol}(M_0)\geq\lim_{k\to\infty}v_{m^V}v_{m^H-1}
\left(\delta_2(\widehat r_1)-\delta_2(1/k)\right)=v_{m^V}v_{m^H-1}\delta_2(\widehat r_1)}\\
\hspace{1.55truecm}\displaystyle{=v_{m^V}v_{m^H-1}(\delta_2\circ\delta_1^{-1})
\left(\frac{{\rm Vol}(D)}{v_{m^V}{\rm Vol}(B)}\right),}
\end{array}$$
where we note that $a_{r_B}=1$ because $\overline M$ is of non-compact type.  
This contradicts the assumption.  Hence we obtain $\inf_{(\xi,t)\in M\times[0,T)}\widehat r(\xi,t)>0$.  
Therefore, by using Theorem B, we can derive that $T=\infty$ and that $f_t$ converges to a tube of constant 
mean curvature over $B$ in the $C^{\infty}$-topology as $t\to\infty$.  \hspace{3.7truecm}q.e.d.

\section{Appendix} 
Polars and meridians of symmetric spaces of compact type are reflective.  
See [CN] about the classification of polars and meridians of symmetric spaces of compact type.  
If $F\subset\overline M$ is one of meridians in irreducible rank two symmetric spaces of compact type in Table 2 
and if ${\cal D}$ is the corresponding distribution on $F$ as in Table 2, then 
${\cal D}_x$ is a common eigenspace of the family $\{R(\cdot,\xi)\xi\,\vert\,\xi\in T^{\perp}_xF\}$ 
for any $x\in F$.  
In fact, we can confirm this fact as follows.  
First we consider the case of $(1)-(4)$ in Tables 2.  
In these cases, we have ${\cal K}=\{1,2\}$ and $\tau_x^{-1}({\cal D}_x)=\mathfrak p_{2\beta}\cap\mathfrak p'$ for any 
$x\in F$ and any $\xi\in T_x^{\perp}F$, where ${\cal K}$ and $\mathfrak p_{2\beta}$ are the quantities defined 
as in Section 2 for $\mathfrak b={\rm Span}\{\tau_x^{-1}\xi\}$.  
That is, we have 
$$R(({\rm grad}\,\psi)_x,\xi)\xi=-4\sqrt{\varepsilon}^2\beta(\tau_x^{-1}\xi)^2({\rm grad}\,\psi)_x$$
for any positive function $\psi$ over $F$ with $({\rm grad}\,\psi)_x\in {\cal D}_x$ ($x\in F$).  
Next we consider the case of $(5)$ in Table 2.  
In these cases, ${\cal K}=\{1\}$ holds, for one of ${\cal D}=TSp(1)$'s in (5), 
$\tau_x^{-1}({\cal D}_x)=\mathfrak p_0\cap\mathfrak p'$ ($x\in F,\,\,\xi\in T_x^{\perp}F$) holds and, for another 
${\cal D}=TSp(1)$'s in (5), $\tau_x^{-1}({\cal D}_x)=\mathfrak p_{\beta}\cap\mathfrak p'$ ($x\in F,\,\,\xi\in T_x^{\perp}F$) 
holds, where ${\cal K},\,\mathfrak p_0$ and $\mathfrak p_{\beta}$ are the quantities defined as in Section 2 for 
$\mathfrak b={\rm Span}\{\tau_x^{-1}\xi\}$.  
That is, in the first case, we have 
$$R(({\rm grad}\,\psi)_x,\xi)\xi=0$$
for any positive function $\psi$ over $F$ with $({\rm grad}\,\psi)_x\in {\cal D}_x$ ($x\in F$), and in the second case, 
$$R(({\rm grad}\,\psi)_x,\xi)\xi=-\sqrt{\varepsilon}^2\beta(\tau_x^{-1}\xi)^2({\rm grad}\,\psi)_x$$
for any positive function $\psi$ over $F$ with $({\rm grad}\,\psi)_x\in {\cal D}_x$ ($x\in F$).  
Easily we can show the following fact.  

\vspace{0.5truecm}

\noindent
{\bf Fact.} {\sl If the initial radius function $r_0$ satisfies $({\rm grad}\,r_0)_x\in {\cal D}_x$ ($x\in F$), 
then $r_t$ satisfies $({\rm grad}\,r_t)_x\in {\cal D}_x$ ($x\in F$) for all $t\in[0,T)$.}

\vspace{0.5truecm}

$$\begin{tabular}{|c|c|c|c|c|}
\hline
 & $\overline M$ & $F$ & $F^{\perp}_x$ & ${\cal D}$\\
\hline
(1) &$SU(3)/SO(3)$ & $S^1\cdot S^2$ (meridian) & ${\Bbb R}P^2$ (polar) & $TS^1$\\
\hline
(2) & $SU(6)/Sp(3)$ & $S^1\cdot S^5$ (meridian) & ${\Bbb Q}P^2$ (polar) & $TS^1$\\
\hline
(3) & $SU(3)$ & $S^1\cdot S^3$ (meridian) & ${\Bbb C}P^2$ (polar) & $TS^1$\\
\hline
(4) & $E_6/F_4$ & $S^1\cdot S^9$ (meridian) & ${\Bbb O}P^2$ (polar) & $TS^1$\\
\hline
(5) & $Sp(2)$ & $Sp(1)\times Sp(1)$ (meridian) & $S^4$ (polar) & {\rm One}\,\,{\rm of}\,\,$TSp(1)${\rm's}\\
\hline
\end{tabular}$$

\vspace{0.2truecm}

\centerline{($\overline M\,:\,$ an irreducible rank two symmetric space of compact type)}

\vspace{0.5truecm}

\centerline{\bf Table 2.}


\vspace{0.5truecm}

See [He] about the notations in these tables for example.  
Let $B$ be a compact closed domain (in $F$) such that $\pi^{-1}(B)$ is as in Figure 7 
and $r$ be a positive function over $B$ such that the lift $({\rm grad}\,r)^L$ of ${\rm grad}\,r$ by $\pi$ 
is as in Figure 7, where $\pi$ is the covering map of $S^1\times S^k$ onto $S^1\cdot S^k$.  
Then $r$ satisfies both $({\rm grad}\,r)_x=0$ ($x\in\partial B$) and 
$({\rm grad}\,r)_x\in {\cal D}_x:=T_xS^1$ ($x\in B$).  

\vspace{0.5truecm}

\centerline{
\unitlength 0.1in
%
\hspace{2.65truecm}}

\vspace{0.5truecm}

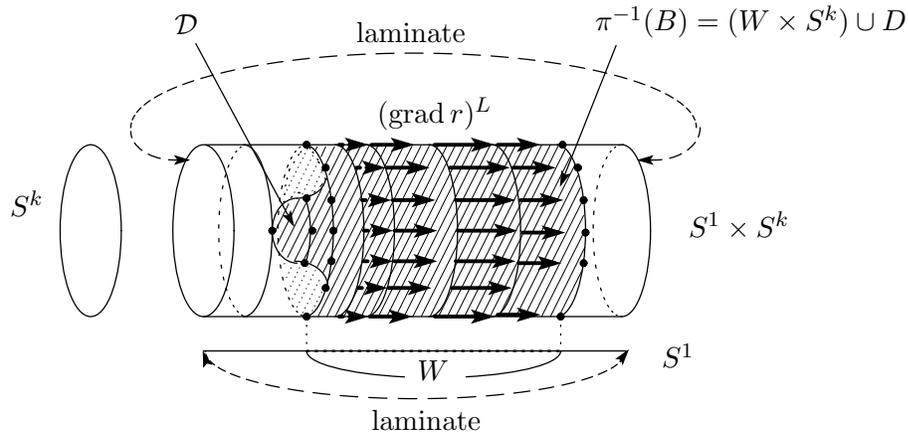
\centerline{{\bf Figure 10.}}

\vspace{0.5truecm}

\vspace{1truecm}

\centerline{{\bf References}}

\vspace{0.5truecm}

{\small 
\noindent
[A1] M. Athanassenas, Volume-preserving mean curvature flow of 
rotationally symme-

tric surfaces, Comment. Math. Helv. {\bf 72} (1997) 52-66.

\noindent
[A2] M. Athanassenas, Behaviour of singularities of the rotationally 
symmetric, volume-

preserving mean curvature flow, Calc. Var. {\bf 17} (2003) 1-16.

\noindent
[BV] J. Berndt and L. Vanhecke, Curvature-adapted submanifolds, Nihonkai Math. J. {\bf 3} 

(1992) 177-185.

\noindent
[BT] J. Berndt and H. Tamaru, Cohomogeneity one actions on noncompact symmetric sp-

aces with a totally geodesic singular orbit, Tohoku Math. J. (2) {\bf 56} (2004), 163--177.

\noindent
[CM1] E. Cabezas-Rivas and V. Miquel, 
Volume preserving mean curvature flow in the 

hyperbolic space, Indiana Univ. Math. J. {\bf 56} 
(2007) 2061-2086.

\noindent
[CM2] E. Cabezas-Rivas and V. Miquel, 
Volume-preserving mean curvature flow of 

revolution hypersurfaces in a rotationally 
symmetric space, Math. Z. {\bf 261} (2009) 

489-510.

\noindent
[CM3] E. Cabezas-Rivas and V. Miquel, Volume-preserving mean curvature flow of 

revolution hypersurfaces between two equidistants, 
Calc. Var. {\bf 43} (2012) 185-210.

\noindent
[CN] B. Y. Chen and T. Nagano, Totally geodesic submanidfold in symmetric spaces II, 

Duke Math. J. {\bf 45} (1978) 405-425.

\noindent
[EH] K. Ecker and G. Huisken, Interior estimates for hypersurfaces moving by mean 

curvature, Invent. Math. {\bf 105} (1991) 547-569.

\noindent
[He] S. Helgason, Differential geometry, Lie groups and symmetric spaces, Academic Press, 

New York, 1978.

\noindent
[Hu1] G. Huisken, Flow by mean curvature of convex surfaces into spheres, J. Differential 

Geom. {\bf 20} (1984) 237-266.

\noindent
[Hu2] G. Huisken, Contracting convex hypersurfaces in Riemannian manifolds by their mean 

curvature, Invent. math. {\bf 84} (1986) 463-480.

\noindent
[Hu3] G. Huisken, Parabolic Monge-Ampere equations on Riemannian manifolds, J. Funct. 

Anal. {\bf 147} (1997) 140-163.

\noindent
[KT] T. Kimura and M. S. Tanaka, Stability of certain minimal submanifolds in compact 

symmetric spaces of rank two, Differential Geom. Appl. {\bf 27} (2009) 23--33.

%

\noindent
[Ko1] N. Koike, Tubes of non-constant radius in symmetric spaces, 
Kyushu J. Math. {\bf 56} 

(2002) 267-291.

\noindent
[Ko2] N. Koike, Collapse of the mean curvature flow for equifocal 
submanifolds, Asian J. 

Math. {\bf 15} (2011) 101-128.

\noindent
[Ko3] N. Koike, Collapse of the mean curvature flow for isoparametric 
submanifolds in 

a symmetric space of non-compact type, Kodai Math. J. {\bf 37} (2014) 355-382.

\noindent
[L1] D. S. P. Leung, Errata: ``On the classification of reflective submanifolds 
of Riemannian 

symmetric spaces'', (Indiana Univ. Math. J. {\bf 24} (1974) 327-339), Indiana Univ. 
Math. 

J. {\bf 24} (1975) 1199.

\noindent
[L2] D. S. P. Leung, Reflective submanifolds. III. Congruenccy of isometric 
reflective subma-

nifolds and corrigenda to the classification of reflective submanifolds, 
J. Differential 

Geom. {\bf 14} (1979) 167-177.  

\noindent
[M] J.A. McCoy, Mixed volume-preserving curvature flows, Calc. Var. Partial Differ. Equ. 

{\bf 24} (2005) 131-154.

\noindent
[PT] R.S. Palais and C.L. Terng, Critical point theory and submanifold 
geometry, Lecture 

Notes in Math. {\bf 1353}, Springer, Berlin, 1988.

\vspace{0.5truecm}

{\small 
\rightline{Department of Mathematics, Faculty of Science}
\rightline{Tokyo University of Science, 1-3 Kagurazaka}
\rightline{Shinjuku-ku, Tokyo 162-8601 Japan}
\rightline{(koike@ma.kagu.tus.ac.jp)}
}

\end{document}